\documentclass{article}
\usepackage[utf8]{inputenc}
\usepackage{amssymb}
\usepackage{amsbsy}
\usepackage{latexsym}
\usepackage{amsmath}
\usepackage{etoolbox}
\usepackage{array}
\usepackage{multirow}
\usepackage{graphicx}
\usepackage{caption}
\usepackage{subcaption}
\usepackage{algorithmic}
\usepackage{float}
\usepackage[hidelinks]{hyperref}
\usepackage{url} 
\usepackage{cleveref}
\usepackage{algorithm}
\usepackage{amsthm}
\usepackage{verbatim}
\usepackage{color}
\usepackage{geometry}  
\usepackage{dsfont}
\usepackage{orcidlink}

\newtheoremstyle{mythmstyle}
  {6pt}
  {6pt}
  {\itshape}
  {}
  {\bfseries}
  {.}
  {0.5em}
  {}

\theoremstyle{mythmstyle}

\newtheorem{theorem}{Theorem}[section]

\newtheorem{lemma}{Lemma}[section]
\newtheorem{remark}{Remark}[section]
\newtheorem{proposition}{Proposition}[section]

\newtheorem{maintheorem}{Theorem}

\AtBeginEnvironment{equation}{\normalfont}
\AtBeginEnvironment{equation*}{\normalfont}
\AtBeginEnvironment{align}{\normalfont}
\AtBeginEnvironment{align*}{\normalfont}
\AtBeginEnvironment{gather}{\normalfont}
\AtBeginEnvironment{gather*}{\normalfont}
\AtBeginEnvironment{multline}{\normalfont}
\AtBeginEnvironment{multline*}{\normalfont}
\AtBeginEnvironment{eqnarray}{\normalfont}
\AtBeginEnvironment{eqnarray*}{\normalfont}

\newcommand{\avg}[1]{\left\{\!\!\left\{#1\right\}\!\!\right\}}
\newcommand{\diff}[1]{\left[\!\left[#1\right]\!\right]}
\newcommand{\avgL}[1]{\left\{\!\!\!\left\{#1\right\}\!\!\!\right\}}
\newcommand{\diffL}[1]{\left[\!\!\left[#1\right]\!\!\right]}
\newcommand{\modL}[1]{\left|#1\right|}
\newcommand{\paraL}[1]{\left(#1\right)}
\newcommand{\sumintK}[1]{\sum_{\sigma \in \mathcal{E}(K) \bigcap \mathcal{E}_{\text{int}}} \frac{|\sigma|}{|K|}}
\newcommand{\sumintKK}[1]{\sum_{\sigma \in \mathcal{E}(K)} \frac{|\sigma|}{|K|}}
\newcommand{\sumK}[1]{\sum_{K \in \mathcal{T}} |K|}
\newcommand{\sumKK}[1]{\sum_{K \in \mathcal{T}}}
\newcommand{\sumintall}[1]{\sum_{\sigma \in \mathcal{E}} |\sigma|}

\newcommand{\dx}[1]{\ \mathrm{d}x}
\newcommand{\dt}[1]{\ \mathrm{d}t}
\newcommand{\divx}[1]{\text{div}_x}
\newcommand{\divh}[1]{\text{div}_h}

\title{Provably fully discrete energy-stable and asymptotic-preserving scheme for barotropic Euler equations\footnote{ \textbf{Funding:} M. A. was funded by the Gutenberg Research College University Mainz. \\
M.L.-M. gratefully acknowledges the support of DFG Project 525800857 funded within Focused Programme SPP 2410 "Hyperbolic Balance Laws: Complexity, Scales and Randomness" and of the Mainz Institute of Multiscale Modeling. \\
This work was also partially  funded by the DAAD DST
(German-India) Project based personnel exchange programme: Development and analysis of higher-order
structure-preserving numerical methods for hyperbolic balance laws.}}
\author{Megala Anandan \orcidlink{0000-0003-4631-3996}\footnote{Institute for Mathematics, Johannes Gutenberg University of Mainz, 55128 Mainz, Germany. E-mail: manandan@uni-mainz.de} and M\'{a}ria Luk\'{a}\v{c}ov\'{a}-Medvid'ov\'{a} \orcidlink{0000-0002-4351-0161} \footnote{Institute for Mathematics, Johannes Gutenberg University of Mainz, 55128 Mainz, Germany. E-mail: lukacova@uni-mainz.de} }
\date{}

\begin{document}

\maketitle

\begin{abstract}
    We develop structure-preserving finite volume schemes for the barotropic Euler equations in the low Mach number regime. Our primary focus lies in ensuring both the asymptotic-preserving (AP) property and the discrete entropy stability. We construct an implicit-explicit (IMEX) method with suitable acoustic/advection splitting including implicit numerical diffusion that is independent of the Mach number. We prove the positivity of density, the entropy stability, and the asymptotic consistency of the fully discrete numerical method rigorously. Numerical experiments for benchmark problems validate the structure-preserving properties of the proposed method.   
\end{abstract}
\section{Introduction}
We consider the compressible barotropic Euler equations given by:
\begin{gather}
\label{Euler mass}
    \partial_t \varrho + \divx{} (\varrho \mathbf{u}) = 0, \\
\label{Euler mom}
    \partial_t (\varrho \mathbf{u}) + \divx{} (\varrho \mathbf{u} \otimes \mathbf{u} ) + \nabla_x p(\varrho) = \mathbf{0},
\end{gather}
on the domain $\Omega \subset \mathbb{R}^d$ and time $t\in \mathbb{R}^+ \cup \{0\}$. The variables $\varrho: \Omega \times \mathbb{R}^+ \cup \{0\}  \rightarrow \mathbb{R}^+$, $\mathbf{u}: \Omega \times \mathbb{R}^+ \cup \{0\} \rightarrow \mathbb{R}^d$, and $p(\varrho)=\kappa \varrho^{\gamma} \in \mathbb{R}^+$ denote the density, the velocity and the pressure of the fluid respectively. Here, $d=1,2,3$ represents the dimension in space, and $\kappa$, $\gamma > 1$ are constants. The system \eqref{Euler mass}-\eqref{Euler mom} is hyperbolic with eigenvalues (in direction $\mathbf{n} \in \mathbb{R}^d$) $\mathbf{u}\cdot \mathbf{n} - c$ and $\mathbf{u}\cdot \mathbf{n} + c$, and $c=\sqrt{\gamma p/\varrho}$ is the sound speed. The conserved quantities are density, $\varrho$ and momentum, $\varrho \mathbf{u}$. The initial conditions required for the system are $\varrho (\cdot,0)=\varrho^0(\cdot)$ and $\mathbf{u}(\cdot,0)=\mathbf{u}^0(\cdot)$, and the boundary is considered to have periodic or no-flux conditions. \\
A specific non-dimensionalization of \eqref{Euler mass} and \eqref{Euler mom} yields the following dimensionless form:
\begin{gather}
    \label{ND Euler mass}
    \partial_t \varrho + \divx{} (\varrho \mathbf{u}) = 0, \\
\label{ND Euler mom}
    \partial_t (\varrho \mathbf{u}) + \divx{} (\varrho \mathbf{u} \otimes \mathbf{u} ) + \frac{1}{\varepsilon^2} \nabla_x p(\varrho) = \mathbf{0},
\end{gather}
where $\varepsilon$ refers to the Mach number. The equations \eqref{ND Euler mass}-\eqref{ND Euler mom} form a hyperbolic system whose eigenvalues in direction $\mathbf{n}$ are $\mathbf{u}\cdot \mathbf{n} - c/\varepsilon$ and $\mathbf{u}\cdot \mathbf{n} + c/\varepsilon$. This system is endowed with the total energy $E = E_{ke} + \frac{1}{\varepsilon^2}P(\varrho)$, where $E_{ke} = \frac{1}{2}\varrho\modL{\mathbf{u}}^2$ is the kinetic energy and $P(\varrho) = \frac{\kappa \varrho^{\gamma}}{\gamma - 1}$ is the pressure potential. The total energy $E$ is a convex functional with respect to $\paraL{\varrho,\varrho\mathbf{u}}$, and it is a mathematical entropy for the system \eqref{ND Euler mass}-\eqref{ND Euler mom}, as the Hessian of the total energy functional symmetrizes the flux Jacobian of the system. The classical solutions of \eqref{ND Euler mass}-\eqref{ND Euler mom} satisfy the following energy identities:
\begin{gather}
     \partial_t P(\varrho) + \divx{} \paraL{P(\varrho) \mathbf{u}} + p(\varrho) \divx{} \mathbf{u} = 0 \text{ (potential energy identity)}, \\
    \partial_t E_{ke} + \divx{} \paraL{E_{ke} \mathbf{u}} + \frac{1}{\varepsilon^2} \nabla_x p \cdot \mathbf{u} = 0 \text{ (kinetic energy identity)}, \\
    \partial_t E + \divx{} \paraL{\paraL{E + \frac{1}{\varepsilon^2}p}\mathbf{u}} = 0 \text{ (total energy identity)}.
\end{gather}
Moreover, the admissible weak solutions of \eqref{ND Euler mass}-\eqref{ND Euler mom} satisfy the energy inequality,
\begin{equation}
\label{ND Euler Energy}
    \partial_t E + \divx{} \paraL{\paraL{E + \frac{1}{\varepsilon^2}p}\mathbf{u}} \leq 0,
\end{equation}
in the sense of distributions. 
\par As the Mach number becomes small (in the limit $\varepsilon \to 0$), the system \eqref{ND Euler mass}-\eqref{ND Euler mom} reduces to an incompressible system as given by the following proposition. 
\begin{proposition}[Asymptotic limit (see Klainerman and Majda \cite{ap_majda1})] 
\everymath{\normalfont}
    Suppose that the initial data for the compressible system \eqref{Euler mass}-\eqref{Euler mom} satisfies the regularity $(\varrho^0(\cdot),\mathbf{u}^0(\cdot)) \in H^s(\Omega;\mathbb{R}^{d+1})$ with $s>\frac{d}{2}+2$, and is well-prepared, i.e.
        \begin{equation}
        \label{Def: Well-prep ID}
            \varrho^0(\cdot) = \varrho_0^0 + \varepsilon^2 \varrho_2^0 (\cdot); \ 
            u^0(\cdot)=u_0^0 (\cdot) + \varepsilon u_1^0 (\cdot) +  \varepsilon^2 u_2^0 (\cdot),
        \end{equation}
        such that $\varrho_0^0$ is constant and $\divx{} u_0^0 = \divx{} u_1^0 = 0$.
    Then, the limit of the system \eqref{ND Euler mass}-\eqref{ND Euler mom} as $\varepsilon \to 0$ is the incompressible barotropic Euler system:
    \begin{gather}
        \label{Lim div free vel}
        \divx{} \mathbf{u_0} = 0, \\
        \label{Lim vel}
        \partial_t \mathbf{u}_0 + \divx{} (\mathbf{u}_0 \otimes \mathbf{u}_0 ) + \frac{\nabla_x p_2}{\varrho_0} = \mathbf{0}.
    \end{gather}
    $p_2$ is the incompressible pressure obtained as the formal limit of $\frac{p(\varrho)-p_0}{\varepsilon^2}$. Further, the following asymptotic expansions hold:
    \begin{equation}
    \label{Lim den vel exp}
        \varrho(\mathbf{x},t) = \varrho_0 + \varepsilon^2 \varrho_2(\mathbf{x},t), \ 
        \divx{} \mathbf{u}(\mathbf{x},t) = \mathcal{O}(\varepsilon^2), \ 
        p(\mathbf{x},t) = p_0 + \varepsilon^2 p_2(\mathbf{x},t). 
    \end{equation}
\end{proposition}
Thus, as described above, the admissible solutions to the barotropic Euler system \eqref{ND Euler mass}-\eqref{ND Euler mom} admits the energy/entropy inequality \eqref{ND Euler Energy}, and also reduces to the incompressible system \eqref{Lim div free vel}-\eqref{Lim vel} with divergence-free velocity as the Mach number goes to zero. Analogous to the continuous setting, it is desirable that a numerical scheme for \eqref{ND Euler mass}-\eqref{ND Euler mom} satisfies the discrete energy/entropy inequality for all values of the Mach number, and also reduces to a consistent discretization of the incompressible system \eqref{Lim div free vel}-\eqref{Lim vel} as the Mach number goes to zero, without having to impose restrictive time-step conditions or large numerical diffusion coefficients that are dependent on $\varepsilon$. Our aim is to propose an implicit-explicit finite volume method that satisfies the following properties:
\begin{enumerate}
    \item Discrete energy/entropy stability property: The numerical scheme satisfies the discrete energy inequality for all values of $\varepsilon$,
    \item Asymptotic-preserving property:
    \begin{itemize}
        \item The asymptotic expansions in \eqref{Lim den vel exp} are satisfied by the numerical scheme at all time levels,
        \item The numerical scheme reduces to a consistent discretization of the incompressible system \eqref{Lim div free vel}-\eqref{Lim vel} as $\varepsilon \to 0$.
    \end{itemize}
\end{enumerate}

\subsection{Related works}
During the last two decades, the entropy-stable methods have become popular for hyperbolic systems after the seminal works of Tadmor \cite{ent_tadmor1,ent_tadmor2}. Several semi-discrete (discrete in space, and continuous in time) entropy-stable methods for shallow water and Euler systems were developed, e.g. \cite{ent_barth,ent_roe,ent_tadmor_ENO,ent_PC,ent_rayPC,ent_kopriva,ent_MA}. Fully discrete entropy stability has also been done numerically up to machine accuracy in \cite{ent_ranocha}, by introducing a relaxation parameter in the Runge-Kutta update. However, provable fully discrete entropy stability is still an ongoing research topic. Despite this gap in the development of entropy-stable schemes, discrete entropy stability is an important characteristic of a numerical scheme as it is considered as a nonlinear stability criterion for hyperbolic systems of conservation laws \cite{ent_harten}, and it is of fundamental importance for analyzing the convergence of numerical schemes for compressible flows \cite{ent_feireisl1,ent_feireisl2,ent_feireisl_book}. 

\par On the other hand, as $\varepsilon \to 0$, the hyperbolic system converges to a mixed hyperbolic-elliptic system \cite{ap_majda1,ap_majda2,ap_schochet}. These systems tend to become stiff with a wide disparity in propagating wave speeds. Explicit numerical methods for such systems require $\varepsilon-$dependent stability condition on the time step, while the fully implicit numerical methods require iterative solvers due to the nonlinearity of the systems. Both these approaches make the computations unreasonably expensive and may even yield wrong results due to the stiffness. The Godunov-type compressible solvers suffer from loss of accuracy as numerical dissipation scales as $\frac{1}{\varepsilon}$ \cite{ap_dellacherie}. Thus, semi-implicit numerical methods with $\varepsilon-$independent time step conditions have become popular for such systems. Such methods are known as \textit{asymptotic-preserving (AP)} methods and were introduced by Jin \cite{ap_jin1,ap_jin2}. Asymptotic-preserving properties in the sense of uniform consistency with respect to the small parameter $\varepsilon$ have been studied for related problems, e.g. \cite{ap_degond1,ap_tang,ap_jin3,ap_lukacova1,ap_lukacova2,ap_lukacova3,ap_dimarco,ap_noelle,ap_lukacova4,ap_chertock,ap_boscarino,ap_arun2,ap_MA1,ap_MA2}. These methods typically involve the splitting of the system into implicit and explicit parts, such that the eigenvalues of the explicit parts are non-stiff.  
\par Since the hyperbolic systems are endowed with entropy inequality for all values of the stiff parameter $\varepsilon$, it is interesting to develop numerical methods that are entropy-stable for all $\varepsilon$ while keeping the time step $\varepsilon-$independent. Such methods have been explored for the barotropic Euler system \eqref{ND Euler mass}-\eqref{ND Euler mom} in our recent work \cite{apes_MA}, where different space discretizations were considered and their discrete entropy stability properties were studied numerically for different values of $\varepsilon$. Our aim in the current work is to propose an asymptotic-preserving method for the barotropic Euler system, and to prove that it satisfies fully discrete entropy/energy stability for all $\varepsilon$ under reasonable conditions on time step and numerical diffusion coefficient that are $\varepsilon-$independent. This is one of the first attempts to prove fully discrete entropy stability and asymptotic preservation for an implicit-explicit method on a collocated grid. In the literature, such fully discrete entropy-stable schemes for all $\varepsilon$ have been proposed either for a fully implicit method \cite{apes_lukacova1} or for a velocity stabilized semi-implicit method on staggered Marker and Cell (MAC) grid \cite{ap_arun1,apes_lukacova2}. Further, unlike the common practice in the AP community that assumes asymptotic expansion \eqref{Lim den vel exp} of variables at all times (evolved by the scheme), inspired by \cite{ap_lukacova3} we also prove that the asymptotic expansion holds for all the variables at any time (evolved by the scheme) if the asymptotic expansion is satisfied at the initial time. This confirms the asymptotic consistency of our method rigorously. 

\subsection{Main contributions}
The main contributions of this paper can be stated in the following informal way. 
\begin{maintheorem}[Positivity of density]
\label{Thm: Positive density-informal}
\everymath{\normalfont}
    Let $\varrho_h^n>0$, $\mathbf{u}_h^n$ be piecewise constant functions on $\Omega$ at time level $n$ and suppose that periodic boundary conditions hold. Then, under a suitable condition on the numerical diffusion coefficient (that will be specified later in \eqref{lam: Positivity}) that is independent of $\varepsilon$, the numerical solution $\varrho_h^{n+1}$ (at time level $n+1$) to the implicit-explicit numerical scheme \eqref{Num mass}-\eqref{Num mom} is positive, for all $0\leq n \leq N-1$.
\end{maintheorem}

\begin{maintheorem}[Discrete total energy inequality]
\label{Thm: TE-informal}
\everymath{\normalfont}
     Let $\varrho_h^n>0$, $\mathbf{u}_h^n$ be piecewise constant functions on $\Omega$ at time level $n$, and $(\varrho_h^{n+1}, \mathbf{u}_h^{n+1})$ be the solution (at time level $n+1$) to the implicit-explicit numerical scheme \eqref{Num mass}-\eqref{Num mom} for $0\leq n \leq N-1$. Then, under a suitable condition on the numerical diffusion coefficient (that will be specified later in \eqref{lambda cond TE}) that is independent of $\varepsilon$, the following discrete total energy inequality holds:
    \begin{equation}
        \int_\Omega \paraL{\frac{1}{2} \varrho_h^{n+1} \modL{\mathbf{u}_h^{n+1}}^2 + \frac{1}{\varepsilon^2} P\paraL{\varrho_h^{n+1}}} \dx{} \leq \int_\Omega \paraL{\frac{1}{2} \varrho_h^{n} \modL{\mathbf{u}_h^{n}}^2 + \frac{1}{\varepsilon^2} P\paraL{\varrho_h^{n}}} \dx{}.
    \end{equation}
\end{maintheorem}

\begin{maintheorem}[Asymptotic-preserving property]
\label{Thm: AP-informal}
\everymath{\normalfont}
    Let the initial data be piecewise constant functions satisfying $\varrho_h^0 = \varrho_0 + \varepsilon^2 \varrho_{h_2}^0$ and $\mathbf{u}_{h}^0=\mathcal{O}(1)$. Then the numerical solution $(\varrho_h^{n}, \mathbf{u}_h^{n})$ to the linear system \eqref{Update mass}-\eqref{Update mom} corresponding to the implicit-explicit numerical scheme satisfies $\varrho_h^n = \varrho_0 + \varepsilon^2 \varrho_{h_2}^n$, $\mathbf{u}_{h}^n=\mathcal{O}(1)$, and $\divh{} \mathbf{u}_{h}^n = \mathcal{O}\paraL{\varepsilon^2}$, for $1\leq n \leq N$. \\
    Further, in the limit as $\varepsilon \to 0$, the momentum update equation in the numerical scheme becomes a consistent approximation of \eqref{Lim vel}. 
\end{maintheorem}
The precise statements and proofs of the above theorems are given in \Cref{Subsec: Positivity}, \Cref{Subsec: Discrete energy}, and \Cref{Subsec: AP} respectively, after introducing the necessary notations and auxiliary lemmas. 

\subsection{Outline}
The remainder of the manuscript is organized as follows: In \Cref{Sec: Preliminaries}, we present the preliminaries, such as the notations involved in the space discretization. Then in \Cref{Sec: Numerical scheme}, we propose the numerical scheme and prove the positivity of density (in \Cref{Subsec: Positivity}), the fully discrete energy stability (in \Cref{Subsec: Discrete energy}), and the asymptotic consistency (in \Cref{Subsec: AP}) properties. Finally, in \Cref{Sec: Num results}, the energy-stable and asymptotically consistent numerical scheme is validated on benchmark problems.

\section{Preliminaries}
\label{Sec: Preliminaries}
In this section, we introduce the main notations and fundamental lemmas associated with the space discretization, which will facilitate a clear presentation of the main results in the next section. 
\subsection{Space discretization}
The physical domain $\Omega$ is discretized as:
\begin{equation}
    \overline{\Omega} = \bigcup_{K \in \mathcal{T}_h} K,
\end{equation}
where $\mathcal{T}_h$ is a set of compact polygons or polyhedrons (if $d=2$ or $d=3$ respectively). For $K, L \in \mathcal{T}_h$, if $K \cap L \neq \emptyset $, $K \neq L$, then $K \cap L$ is either a common face, a common edge, or a common vertex. The set of all faces $\sigma$ of an element $K \in \mathcal{T}_h$ is denoted by $\mathcal{E}(K)$. The set of all faces of the mesh is denoted by $\mathcal{E}$. $\mathcal{E}_{ext} = \mathcal{E} \cap \partial\Omega$ and $\mathcal{E}_{int} = \mathcal{E} \backslash \partial\Omega$ denote the set of all exterior and interior faces, respectively. For periodic boundary conditions, we set $\Omega = \mathbb{T}^d$, $\mathcal{E}_{int} = \mathcal{E}$ and $\mathcal{E}_{ext} = \emptyset$. 
\par $\sigma = K | L$ denotes the face between two neighbors $K,L \in \mathcal{T}_h$ such that $\sigma = \mathcal{E}(K) \cap \mathcal{E}(L)$. For each face $\sigma = K | L$, $\mathbf{n}_{\sigma,K}$ denotes the normal vector of $\sigma$ pointing outwards from K (so that $\mathbf{n}_{\sigma,K} = - \mathbf{n}_{\sigma,L}$). $h_K$ and $h_{\sigma}$ denote the diameters of $K$ and $\sigma$ respectively, and the maximal mesh size is given by $h = \max_{K \in \mathcal{T}_h} h_K$. The mesh $\mathcal{T}_h$ is regular and quasi-uniform, that is, there exist positive real numbers $\theta_0$ and $c_0$ independent of $h$ such that,
\begin{equation}
    \inf_{K \in \mathcal{T}_h} \frac{\xi_K}{h_k} \geq \theta_0 \text{ and } c_0 h \leq h_K,
\end{equation}
where $\xi_K$ denotes the diameter of the largest ball included in $K$. $|K|$ and $|\sigma|$ denote the $d$ and $(d-1)$ dimensional Lebesgue measure of an element $K$ and a face $\sigma$, so that $|K| \approx h^d$, $|\sigma| \approx h^{d-1}$ for any $K \in \mathcal{T}_h$, $\sigma \in \mathcal{E}$. For the family of control points or cell centers $\{ x_K\}_{K\in \mathcal{T}_h}$, any direction vector $\overrightarrow{x_Kx_L}$ is normal to $\sigma$ where $\sigma = K|L$. We set $d_{\sigma} = \modL{x_K - x_L}\approx h$.

\subsection{Discrete function spaces}
Let $\mathcal{P}_n(K)$ denote the set of polynomials of degree at most $n$ on the element $K$. Then,
\begin{equation}
    Q_h = \left\{ \phi_h \in L^1(\Omega) \ | \ \phi_h|_K = \phi_K \in \mathcal{P}_0(K) \text{ for all } K \in \mathcal{T}_h\right\},
\end{equation}
is the function space of piecewise constant functions. We introduce the following projection operator,
\begin{equation}
    \Pi_Q : L^1(\Omega) \to Q_h \text{ with } \ \Pi_Q \phi = \sum_{K \in \mathcal{T}_h} \frac{1_K}{|K|} \int_K \phi \dx{}, \ 1_K (x) = \left\{ \begin{matrix}
        1 & \text{if } x \in K \\ 0 & \text{if } x \notin K
    \end{matrix} \right..
\end{equation} 
Analogously, if $\boldsymbol{\varphi}_h = \paraL{\phi_{1,h},\dots, \phi_{d,h}}$ is a vector-valued function, then $\boldsymbol{\varphi}_h \in \boldsymbol{Q}_h$ means $\phi_{i,h} \in Q_h$ for all $i \in \{ 1,\dots,d\}$, and the projection operator also acts componentwise, i.e., $\Pi_Q \boldsymbol{\varphi} = \paraL{\Pi_Q \phi_1,\dots, \Pi_Q \phi_d}$. 
\par Let $\phi_h \in Q_h$ and let $\phi_h$ be smooth on each element $K \in \mathcal{T}_h$, then for each $\sigma \in \mathcal{E}_{int}$,
\begin{equation}
    \phi_h^{\text{out}} = \lim_{\delta \to 0+} \phi_h(x+\delta\mathbf{n}_{\sigma}), \ \phi_h^{\text{in}} = \lim_{\delta \to 0+} \phi_h(x-\delta\mathbf{n}_{\sigma}), \ x \in \sigma, \ \sigma \in \mathcal{E}_{int}.
\end{equation}
For any face $\sigma \in \mathcal{E}_{ext}$, $\phi_h^{\text{out}}$ must be prescribed from the boundary conditions. If a no-flux boundary condition is used, we simply have $\phi_h^{\text{out}} = \phi_h^{\text{in}}$. For any face $\sigma \in \mathcal{E}$, the following notations are used for discrete jump and average, respectively,
\begin{equation}
    \diff{\phi_h}_{\sigma} = \phi_h^{\text{out}} - \phi_h^{\text{in}}, \ \avg{\phi_h}_{\sigma} = \frac{\phi_h^{\text{out}} + \phi_h^{\text{in}}}{2}.
\end{equation}

\subsection{Discrete difference operators}
We introduce the following discrete gradient, discrete divergence, and discrete Laplace operators, 
\begin{gather}
    \nabla_h \phi_h = \sum_{K \in \mathcal{T}_h} \paraL{\nabla_h \phi_h}_K 1_K, \quad   \paraL{\nabla_h \phi_h}_K = \sumintKK{} \avg{\phi_h}_{\sigma} \mathbf{n}_{\sigma,K}, \\
    \text{div}_h \boldsymbol{\varphi}_h = \sum_{K \in \mathcal{T}_h} \paraL{\text{div}_h \boldsymbol{\varphi}_h}_K 1_K, \quad   \paraL{\text{div}_h \boldsymbol{\varphi}_h}_K = \sumintKK{} \avg{\boldsymbol{\varphi}_h}_{\sigma} \cdot \mathbf{n}_{\sigma,K}, \\
    \Delta_h \phi_h = \sum_{K \in \mathcal{T}_h} \paraL{\Delta_h \phi_h}_K 1_K, \quad   \paraL{\Delta_h \phi_h}_K = \sumintKK{} \frac{\diff{\phi_h}_{\sigma,K}}{d_{\sigma}}, 
\end{gather}
for $\phi_h \in Q_h$ and $\boldsymbol{\varphi}_h \in \boldsymbol{Q}_h$. Here $1_K$ is the characteristic function of $K$. The operators $\nabla_h, \Delta_h$ act component-wisely for vector-valued functions. Note here that these discrete gradient and discrete divergence operators simply represent central differencing. 

\begin{lemma}[Discrete grad-div duality]
\everymath{\normalfont}
\label{Lem: Grad-div duality}
    Let $\phi_h \in Q_h$, $\boldsymbol{\varphi}_h \in \boldsymbol{Q}_h$. Let one of the following conditions holds for all $\sigma \in \mathcal{E}_{ext}$:
    \begin{gather}
        \diff{\phi_h}_{\sigma}=0=\diff{\boldsymbol{\varphi}_h}_{\sigma} \cdot \mathbf{n}_{\sigma} \text{ or } \diff{\phi_h}_{\sigma}=0=\avg{\boldsymbol{\varphi}_h}_{\sigma} \cdot \mathbf{n}_{\sigma} \text{ or } \\
        \avg{\phi_h}_{\sigma}=0=\diff{\boldsymbol{\varphi}_h}_{\sigma} \cdot \mathbf{n}_{\sigma} \text{ or } \avg{\phi_h}_{\sigma}=0=\avg{\boldsymbol{\varphi}_h}_{\sigma} \cdot \mathbf{n}_{\sigma}.
    \end{gather}
    Then,
    \begin{equation}
        \sumK{} \phi_K \paraL{\text{div}_h \boldsymbol{\varphi}_h}_K = - \sumK{} \paraL{\nabla_h \phi_h}_K \cdot \boldsymbol{\varphi}_K. 
    \end{equation} 
\end{lemma}

\begin{lemma}[Discrete inequalities]
\everymath{\normalfont}
    For $f_h,g_h \in Q_h$, the following discrete product rule and algebraic identity hold:
    \begin{gather}
        \diff{f_h g_h} = \avg{f_h} \diff{g_h} + \diff{f_h} \avg{g_h}, \\
        \avg{f_h g_h} = \avg{f_h} \avg{g_h} + \frac{1}{4} \diff{f_h} \diff{g_h}. 
    \end{gather}
    
\end{lemma}

\section{Numerical scheme}
\label{Sec: Numerical scheme}
The numerical scheme is initialized by considering the initial approximations as the projections of the initial data:
\begin{equation}
\label{IC}
    \varrho^0_h|_K = \varrho^0_K = \frac{1}{|K|} \int_K \varrho^0(\cdot) \dx{}, \ \mathbf{u}_h^0|_K = \mathbf{u}_K^0 = \frac{1}{|K|} \int_K \mathbf{u}^0(\cdot) \dx{} \text{ for each } K \in \mathcal{T}_h. 
\end{equation}
The time step size is $\Delta t^n = t^{n+1}-t^n$ with a discretization $0=t^0<t^1<\dots<t^{N}=T$ of the time interval $[0,T]$. 
\par Having the initial data defined as above, we now propose the following semi-implicit finite volume method for \eqref{ND Euler mass}-\eqref{ND Euler mom}: 
\begin{gather}
    \label{Num mass}
    \frac{\varrho_{K}^{n+1}-\varrho_{K}^{n}}{\Delta t^n} +  \paraL{\divh{}\paraL{\varrho_h^{n+1}\mathbf{u}_h^{n+1}}}_K  - \sumintKK{} \lambda_{\sigma,K} \diff{\varrho_h^{n+1}}_{\sigma,K}  = 0,  \\
\label{Num mom}
    \frac{\varrho_{K}^{n+1} \mathbf{u}_{K}^{n+1}-\varrho_{K}^{n}\mathbf{u}_{K}^{n}}{\Delta t^n} + \paraL{ \divh{}\paraL{\varrho_h^{n} \mathbf{u}_h^{n} \otimes \mathbf{u}_h^{n}}}_K  - \sumintKK{} \lambda_{\sigma,K} \diff{\varrho_h^{n+1}\mathbf{u}_h^{n+1}}_{\sigma,K} + \frac{1}{\varepsilon^2}  \paraL{ \nabla_h p\paraL{\varrho_h^{n+1}}}_K  = \mathbf{0},
\end{gather}
where $\lambda_{\sigma,K} \geq 0 $ for $\sigma \in \mathcal{E}(K)$, is the numerical diffusion coefficient. On a uniform quadrilateral mesh, $\divh{} (\cdot)$ and $\nabla_h (\cdot)$ simply represent central discretisations. The method involves central differencing of the pressure gradient and of the mass and momentum convection operators, along with implicit numerical diffusion, 
\begin{gather}
    \sumintKK{} \lambda_{\sigma,K} \diff{\varrho_h^{n+1}}_{\sigma,K} = h_K \paraL{\Delta_{h} \paraL{\lambda_{\sigma} \varrho_h^{n+1}}}_K, \\
    \sumintKK{} \lambda_{\sigma,K} \diff{\varrho_h^{n+1}\mathbf{u}_h^{n+1}}_{\sigma,K} = h_K  \paraL{\Delta_{h} \paraL{\lambda_{\sigma} \varrho_h^{n+1}\mathbf{u}_h^{n+1}}}_K,
\end{gather}
added to both mass and momentum equations.
\par The momentum convection is treated explicitly, while the mass convection and the stiff pressure gradient are treated implicitly in order to facilitate asymptotic preservation. Although several terms in the scheme are implicit, the method can be evolved in only two steps (one step for $\varrho_K^{n+1}$, and the other step for $\mathbf{u}_K^{n+1}$) if the pressure is simply linearized about $p_0$ as will be seen in \Cref{Subsec: Implementation}.    

\subsection{Positivity of density}
In this subsection, we prove that the numerical scheme \eqref{Num mass}-\eqref{Num mom} results in positive density under suitable conditions on the numerical diffusion coefficient. In particular, we present the main lemmas leading to the proof of the formal statement of Theorem~\ref{Thm: Positive density-informal}.  
\label{Subsec: Positivity}
\begin{lemma}[Discrete renormalization identity]
\everymath{\normalfont}
\label{Lem: RE}
    Let $\varrho_h^n \in Q_h$, $\varrho_h^n|_K$ be a solution to \eqref{Num mass} and $\mathbf{u}_h^n \in \boldsymbol{Q}_h$ for $1\leq n \leq N$. Let $\varrho_h^0$ and $\mathbf{u}_h^0$ be the initial data given by \eqref{IC} for each $K\in \mathcal{T}_h$. We suppose that periodic boundary conditions hold. Then, for any function $B \in C^1(\mathbb{R})$, the following holds:
    \begin{multline}
    \label{RE}
        \sumK{} \frac{B\left(\varrho_{K}^{n+1}\right)-B\left(\varrho_{K}^{n}\right)}{\Delta t^n} + \sumK{}  \paraL{\varrho_K^{n+1}B'\paraL{\varrho_K^{n+1}}-B\paraL{\varrho_K^{n+1}}} \paraL{\divh{} \paraL{\mathbf{u}_h^{n+1}}}_K \\ =  - \sumintall{} \lambda_{\sigma} \diff{B'\left(\varrho_h^{n+1}\right)}_{\sigma}\diff{\varrho_h^{n+1}}_{\sigma}  - \sumK{}   \frac{ B\left(\varrho_{K}^{n}\right) -  B\left(\varrho_{K}^{n+1}\right) - B'\left(\varrho_K^{n+1}\right) \left(\varrho_{K}^{n}-\varrho_{K}^{n+1}\right) }{\Delta t^n} \\ + \sumK{} \sumintKK{} \frac{1}{2} \left( \diff{B\left(\varrho_h^{n+1}\right)}_{\sigma,K}  -  B'\left(\varrho_K^{n+1}\right) \diff{\varrho_h^{n+1}}_{\sigma,K}  \right) \mathbf{u}_L^{n+1} \cdot \mathbf{n}_{\sigma,K}  
    \end{multline}
with $\sigma = K|L \in \mathcal{E}(K)$ and $\lambda_{\sigma,K} = \lambda_{\sigma,L} = \lambda_{\sigma}$.
\end{lemma}
\begin{proof}
    See Appendix~\ref{Proof: RE}
\end{proof}

If $B$ is a continuously differentiable convex function, then the following hold:
    \begin{gather}
        - \sumintall{} \lambda_{\sigma} \diff{B'\left(\varrho_h^{n+1}\right)}_{\sigma}\diff{\varrho_h^{n+1}}_{\sigma} \leq 0, \\
        - \sumK{}   \frac{ B\left(\varrho_{K}^{n}\right) -  B\left(\varrho_{K}^{n+1}\right) - B'\left(\varrho_K^{n+1}\right) \left(\varrho_{K}^{n}-\varrho_{K}^{n+1}\right) }{\Delta t^n} \leq 0, \\ 
        \sumK{} \sumintKK{} \frac{1}{2} \left( \diff{B\left(\varrho_h^{n+1}\right)}_{\sigma,K}  -  B'\left(\varrho_K^{n+1}\right) \diff{\varrho_h^{n+1}}_{\sigma,K}  \right) \modL{\mathbf{u}_L^{n+1} \cdot \mathbf{n}_{\sigma,K}} \geq 0. 
    \end{gather}
    For each $\sigma=K|L$, taking $B''\paraL{\varrho^{n+1}_{*,K}}\geq0, B''\paraL{\varrho^{n+1}_{*,L}}\geq0$  with $\varrho^{n+1}_{*,K}, \varrho^{n+1}_{*,L} \in \left[ \varrho_K^{n+1}, \varrho_L^{n+1} \right]$ such that
    \begin{eqnarray}
        B''\paraL{\varrho^{n+1}_{*,L}} \frac{\diff{\varrho_h^{n+1}}_{\sigma}^2}{2} &=& \diff{B\left(\varrho_h^{n+1}\right)}_{\sigma,K}  -  B'\left(\varrho_K^{n+1}\right) \diff{\varrho_h^{n+1}}_{\sigma,K}, \\
        B''\paraL{\varrho^{n+1}_{*,K}} \frac{\diff{\varrho_h^{n+1}}_{\sigma}^2}{2} &=& \diff{B\left(\varrho_h^{n+1}\right)}_{\sigma,L}  -  B'\left(\varrho_L^{n+1}\right) \diff{\varrho_h^{n+1}}_{\sigma,L} ,
    \end{eqnarray}
    we obtain,
    \begin{multline}
        \sumK{} \sumintKK{} \frac{1}{2} \left( \diff{B\left(\varrho_h^{n+1}\right)}_{\sigma,K}  -  B'\left(\varrho_K^{n+1}\right) \diff{\varrho_h^{n+1}}_{\sigma,K}  \right) \mathbf{u}_L^{n+1} \cdot \mathbf{n}_{\sigma,K} \\ 
        = \sumintall{} \frac{\diff{\varrho_h^{n+1}}_{\sigma}^2}{4} \diff{B''\paraL{\varrho^{n+1}_{*}}\mathbf{u}_h^{n+1}}_{\sigma} \cdot \mathbf{n}_{\sigma}. 
    \end{multline}
    Therefore, \eqref{RE} becomes
    \begin{equation}
    \label{RE ineq}
        \sumK{} \frac{B\left(\varrho_{K}^{n+1}\right)-B\left(\varrho_{K}^{n}\right)}{\Delta t^n} + \sumK{}  \paraL{\varrho_K^{n+1}B'\paraL{\varrho_K^{n+1}}-B\paraL{\varrho_K^{n+1}}} \paraL{\divh{} \paraL{\mathbf{u}_h^{n+1}}}_K \leq 0,
    \end{equation}
    if for each $\sigma = K|L$ with $\diff{\varrho_h^{n+1}}_{\sigma} \neq 0$ and $\diff{B'\paraL{\varrho_h^{n+1}}}_{\sigma} \neq 0$, $\lambda_{\sigma}$ satisfies the following condition for a convex $C^1$ function $B$:
    \begin{equation}
    \label{lambda cond RE}
        \lambda_{\sigma} \geq \frac{ \diff{\varrho_h^{n+1}}_{\sigma}^2\diff{B''\paraL{\varrho^{n+1}_{*}} \mathbf{u}_h^{n+1}}_{\sigma} \cdot \mathbf{n}_{\sigma}}{4\diff{B'\left(\varrho_h^{n+1}\right)}_{\sigma}\diff{\varrho_h^{n+1}}_{\sigma}}. 
    \end{equation} 

\begin{theorem}[Positivity of density]
\everymath{\normalfont}
\label{Thm: Positive density}
    Let the assumptions of \Cref{Lem: RE} hold, and let $\varrho_h^0>0$. For a continuously differentiable convex function $B$, let $\lambda_{\sigma}$ satisfy 
    \begin{equation}
    \label{lam: Positivity}
        \lambda_{\sigma} = \max \left\{ 1_{\diff{\varrho_h^{n+1}}_{\sigma} \neq 0,\diff{B'\paraL{\varrho_h^{n+1}}}_{\sigma} \neq 0} \frac{ \diff{\varrho_h^{n+1}}_{\sigma}^2\diff{B''\paraL{\varrho^{n+1}_{*}} \mathbf{u}_h^{n+1}}_{\sigma} \cdot \mathbf{n}_{\sigma}}{4\diff{B'\left(\varrho_h^{n+1}\right)}_{\sigma}\diff{\varrho_h^{n+1}}_{\sigma}}, \frac{\modL{\mathbf{u}_L^{n+1}\cdot \mathbf{n}_{\sigma,K}}}{2}  , \frac{\modL{\mathbf{u}_K^{n+1}\cdot \mathbf{n}_{\sigma,L}}}{2}  \right\},
    \end{equation}
    for each $\sigma = K|L$, and $1_{\diff{\varrho_h^{n+1}}_{\sigma} \neq 0,\diff{B'\paraL{\varrho_h^{n+1}}}_{\sigma} \neq 0}$ represents the characteristic function of the set of faces satisfying $\diff{\varrho_h^{n+1}}_{\sigma} \neq 0$, $\diff{B'\paraL{\varrho_h^{n+1}}}_{\sigma} \neq 0$. Then any solution of \eqref{Num mass} satisfies $\varrho_h^n>0$ for all $1\leq n \leq N$.
\end{theorem}
\begin{proof}
    We use mathematical induction, starting with the induction hypothesis $\varrho_h^{n}>0$. The goal would be to show $\varrho_h^{n+1}>0$.
    \par Since for each $\sigma = K|L$ with $\diff{\varrho_h^{n+1}}_{\sigma} \neq 0$ and $\diff{B'\paraL{\varrho_h^{n+1}}}_{\sigma} \neq 0$, $\lambda_{\sigma}$ satisfies \eqref{lambda cond RE} for a continuously differentiable convex function $B$, the renormalized continuity equation \eqref{RE} becomes the inequality \eqref{RE ineq}. 
    The nonnegativity of density could be easily shown with $B(\varrho) = \max\{ 0, - \varrho\}\geq 0$, which is $0$ if $\varrho \geq 0$, and satisfies $\varrho B'(\varrho)=B(\varrho)$. This is not a $C^1$ function, however, we can easily construct an approximate sequence
    \begin{gather*}
        B_{\delta} : \mathbb{R} \to [0,\infty), B_{\delta} \in C^1(\mathbb{R}), B_{\delta}(\varrho) = 0 \text{ for } \varrho \geq \delta, B_{\delta}'(0) = 0, \\
        B_{\delta}(\varrho) \to B(\varrho) = \max\{ 0, -\varrho\} \text{ uniformly in } \mathbb{R}, B_{\delta}'(\varrho) = -1 \text{ for all } \varrho < -\delta. 
    \end{gather*}
    Plugging $B_{\delta}$ in \eqref{RE}, and since $\lambda_{\sigma}$ satisfies \eqref{lambda cond RE} for each $B_{\delta}$, the limit $\delta \to 0$ of the renormalized continuity equation \eqref{RE} becomes,
    \begin{equation}
        \sumK{} B\left(\varrho_{K}^{n+1}\right) \leq 0.
    \end{equation}
    As $B$ is a non-negative function, we have $B\left(\varrho_{K}^{n+1}\right) = 0$ for all $K \in \mathcal{T}_h$. Thus, we have proved $\varrho_{K}^{n+1} \geq 0$ for $K \in \mathcal{T}_h$. 
    \par We use proof by contradiction to show that $\varrho_{K}^{n+1} \neq 0$ for $K \in \mathcal{T}_h$. By inserting $\varrho_{K}^{n+1} = 0$ in \eqref{Num mass}, we obtain,
    \begin{equation}
        \frac{0 - \varrho_K^n}{\Delta t^n}  =  \sumintKK{} \paraL{-\frac{1}{2}\varrho_L^{n+1}\mathbf{u}_L^{n+1} \cdot \mathbf{n}_{\sigma,K} + \lambda_{\sigma,K}\varrho_L^{n+1} } \geq 0
    \end{equation}
    since $\lambda_{\sigma,K} \geq \frac{\modL{\mathbf{u}_L^{n+1}\cdot \mathbf{n}_{\sigma,K}}}{2}$. This contradicts the hypothesis that $\varrho_K^n > 0$, and hence $\varrho_K^{n+1} \neq 0$.
\end{proof}

\subsection{Discrete energy inequality property}
In this subsection, we prove that the numerical scheme \eqref{Num mass}-\eqref{Num mom} yields discrete energy stability under suitable conditions on the numerical diffusion coefficient. In particular, we present the crucial lemmas on discrete internal and kinetic energy inequalities that lead to the proof of the formal statement of Theorem~\ref{Thm: TE-informal}. 
\label{Subsec: Discrete energy}
\begin{lemma}[Discrete internal energy inequality]
\everymath{\normalfont}
\label{Lem: IE_ineq}
     Let $\varrho_h^n \in Q_h$, $\varrho_h^n|_K$ be a solution to \eqref{Num mass} and $\mathbf{u}_h^n \in \boldsymbol{Q}_h$ for $1\leq n \leq N$. Let $\varrho_h^0$ and $\mathbf{u}_h^0$ be the initial data given by \eqref{IC} for each $K\in \mathcal{T}_h$. We suppose that periodic boundary conditions hold. Let $\lambda_{\sigma}$ satisfy the following (independently of $\varepsilon$) for each $\sigma = K|L$ with $\diff{\varrho_h^{n+1}}_{\sigma} \neq 0$: 
     \begin{equation}
     \label{lambda cond IE}
         \lambda_{\sigma} \geq   \frac{ \diff{\varrho_h^{n+1}}_{\sigma}^2\diff{P''\paraL{\varrho^{n+1}_{*}} \mathbf{u}_h^{n+1}}_{\sigma} \cdot \mathbf{n}_{\sigma}}{4\diff{P'\left(\varrho_h^{n+1}\right)}_{\sigma}\diff{\varrho_h^{n+1}}_{\sigma}},
     \end{equation}
     where $P''\paraL{\varrho^{n+1}_{*,K}}\geq0, P''\paraL{\varrho^{n+1}_{*,L}}\geq0$ with $\varrho^{n+1}_{*,K}, \varrho^{n+1}_{*,L} \in \left[ \varrho_K^{n+1}, \varrho_L^{n+1} \right]$ are such that
    \begin{eqnarray}
        \label{IE_CL}
        P''\paraL{\varrho^{n+1}_{*,L}} \frac{\diff{\varrho_h^{n+1}}_{\sigma}^2}{2} &=& \diff{P\left(\varrho_h^{n+1}\right)}_{\sigma,K}  -  P'\left(\varrho_K^{n+1}\right) \diff{\varrho_h^{n+1}}_{\sigma,K}, \\
        \label{IE_CK}
        P''\paraL{\varrho^{n+1}_{*,K}} \frac{\diff{\varrho_h^{n+1}}_{\sigma}^2}{2} &=& \diff{P\left(\varrho_h^{n+1}\right)}_{\sigma,L}  -  P'\left(\varrho_L^{n+1}\right) \diff{\varrho_h^{n+1}}_{\sigma,L}. 
    \end{eqnarray}
     Then the following discrete internal energy inequality holds:
    \begin{equation}
    \label{IE_ineq}
         \sumK{} \frac{P\left(\varrho_{K}^{n+1}\right)-P\left(\varrho_{K}^{n}\right)}{\Delta t^n} + \sumK{} p\left(\varrho_K^{n+1}\right) \paraL{\divh{} \paraL{\mathbf{u}_h^{n+1}}}_K  \leq 0.
    \end{equation}
\end{lemma}
\begin{proof}
Plugging $B=P$ in \eqref{RE}, and noting that P is a convex differentiable function, we obtain:
\begin{equation}
    \sumK{} \frac{P\left(\varrho_{K}^{n+1}\right)-P\left(\varrho_{K}^{n}\right)}{\Delta t^n} + \sumK{}  \paraL{\varrho_K^{n+1}P'\paraL{\varrho_K^{n+1}}-P\paraL{\varrho_K^{n+1}}} \paraL{\divh{} \paraL{\mathbf{u}_h^{n+1}}}_K \leq 0,
\end{equation}
if $\lambda_{\sigma}$ satisfies \eqref{lambda cond IE}. The proof is complete as $\varrho P'\paraL{\varrho}-P\paraL{\varrho} = p\left(\varrho \right)$. 
\end{proof}

\begin{remark}
    If the solution to system \eqref{ND Euler mass}-\eqref{ND Euler mom} is smooth, then we have for each $K \in \mathcal{T}_h$ and $1\leq n \leq N$,
\begin{equation}
    \diff{(\cdot)^{n}}_{\sigma,K} = d_{\sigma} D_x (\cdot)_K^{n} \approx h D_x (\cdot)_K^{n},
\end{equation}
with $D_x (\cdot)_K^{n}=\mathcal{O}(1)$. Then, the condition \eqref{lambda cond IE} on $\lambda_{\sigma}$ can be interpreted as
\begin{equation}
\label{Simple lam smooth IE}
    \lambda_{\sigma} \geq C \frac{ h D_x \varrho_K^{n+1} D_x \mathbf{u}_K^{n+1}\cdot \mathbf{n}_{\sigma}}{4 D_x P'\paraL{\varrho_K^{n+1}} }  = \mathcal{O}(h), 
\end{equation}
where $C= \max \left\{P''\paraL{\varrho^{n+1}_{*,K}}, P''\paraL{\varrho^{n+1}_{*,L}} \right\}$. 
\end{remark}

\begin{lemma}[Discrete kinetic energy inequality]
\everymath{\normalfont}
\label{Lem: KE_ineq}
    Let $\varrho_h^n>0 \in Q_h$, $\mathbf{u}_h^n \in \boldsymbol{Q}_h$, $\varrho_h^n|_K$, $\mathbf{u}_h^n|_K$ be solution to \eqref{Num mass}-\eqref{Num mom} for $1\leq n \leq N$. Let $\varrho_h^0$ and $\mathbf{u}_h^0$ be the initial data given by \eqref{IC} for each $K\in \mathcal{T}_h$. We suppose that periodic boundary conditions hold. Then, if the following condition on $\lambda_{\sigma}$ is satisfied independently of $\varepsilon$ for each $\sigma \in \mathcal{E}$ with $\diff{\mathbf{u}_h^{n+1}}_{\sigma} \neq 0$,
    \begin{equation}
    \label{lambda cond KE}
        \lambda_{\sigma} \geq  \frac{ -\avg{\varrho_h^{n+1} \mathbf{u}_h^{n+1}}_{\sigma} \cdot \mathbf{n}_{\sigma}  \diffL{\frac{\modL{\mathbf{u}_h^{n+1}}^2}{2}}_{\sigma}   + \diff{\mathbf{u}_h^{n+1}}_{\sigma} \cdot \avg{\varrho_h^n \mathbf{u}_h^n \mathbf{u}_h^n}_{\sigma} \cdot \mathbf{n}_{\sigma}}{\avg{\varrho_h^{n+1}}_{\sigma} \modL{\diff{\mathbf{u}_h^{n+1}}_{\sigma}}^2},
    \end{equation}
    the following discrete kinetic energy inequality holds:
    \begin{equation}
    \label{KE_ineq}
        \sumK{} \frac{\varrho_{K}^{n+1} \modL{\mathbf{u}_K^{n+1}}^2 - \varrho_{K}^{n} \modL{\mathbf{u}_K^{n}}^2}{2 \Delta t^n} - \sumK{} \frac{p\paraL{\varrho_K^{n+1}}}{\varepsilon^2} \paraL{\text{div}_h \mathbf{u}^{n+1}_h}_K  \leq 0.  
    \end{equation}
\end{lemma}
\begin{proof}
    See Appendix~\ref{Proof: KE}
\end{proof}
\begin{remark}
    The numerator in \eqref{lambda cond KE}  can be expanded as: 
\begin{multline}
    \diff{\mathbf{u}_h^{n+1}}_{\sigma} \cdot \left( \avg{\mathbf{u}_h^n}_{\sigma}\avg{\varrho_h^n}_{\sigma}\avg{\mathbf{u}_h^n}_{\sigma} - \avg{\mathbf{u}_h^{n+1}}_{\sigma}\avg{\varrho_h^{n+1}}_{\sigma}\avg{\mathbf{u}_h^{n+1}}_{\sigma}   \right. \\ \left. + \frac{1}{4} \paraL{\avg{\mathbf{u}_h^n}_{\sigma} \diff{\varrho_h^n}_{\sigma}\diff{\mathbf{u}_h^{n}}_{\sigma} + \diff{\mathbf{u}_h^{n}}_{\sigma} \avg{\varrho_h^n}_{\sigma} \diff{\mathbf{u}_h^{n}}_{\sigma}  + \diff{\mathbf{u}_h^{n}}_{\sigma} \diff{\varrho_h^n}_{\sigma} \avg{\mathbf{u}_h^n}_{\sigma} - \avg{\mathbf{u}_h^{n+1}}_{\sigma}\diff{\varrho_h^{n+1}}_{\sigma} \diff{\mathbf{u}_h^{n+1}}_{\sigma}} \right) \cdot \mathbf{n}_{\sigma}. 
\end{multline}
 Further, if the solution to system \eqref{ND Euler mass}, \eqref{ND Euler mom} is smooth, we have the following:
\begin{gather}
    \avg{\mathbf{u}_h^n}_{\sigma}\avg{\varrho_h^n}_{\sigma}\avg{\mathbf{u}_h^n}_{\sigma} -\avg{\mathbf{u}_h^{n+1}}_{\sigma}\avg{\varrho_h^{n+1}}_{\sigma}\avg{\mathbf{u}_h^{n+1}}_{\sigma}   = \Delta t^n \widetilde{D}_t \avg{\mathbf{u}_h}_{\sigma}\avg{\varrho_h}_{\sigma}\avg{\mathbf{u}_h}_{\sigma}    \\
    \diff{(\cdot)^{n}}_{\sigma,K} = d_{\sigma} D_x (\cdot)_K^{n} \approx h D_x (\cdot)_K^{n}, 
\end{gather}
where $\widetilde{D}_t \avg{\mathbf{u}_h}_{\sigma}\avg{\varrho_h}_{\sigma}\avg{\mathbf{u}_h}_{\sigma}$, $D_x (\cdot)_K^{n}$ are $O(1)$. Then the condition \eqref{lambda cond KE} becomes
\begin{equation}
\label{Simple lam smooth KE}
    \lambda_{\sigma} \geq \frac{D_x \mathbf{u}_K^{n+1} \cdot \widetilde{D}_t \avg{\mathbf{u}_h}_{\sigma}\avg{\varrho_h}_{\sigma}\avg{\mathbf{u}_h}_{\sigma} \Delta t^n}{\avg{\varrho_h^{n+1}} \modL{D_x \mathbf{u}_K^{n+1}}^2 h } + \mathcal{O}\paraL{h} = \mathcal{O}(1),
\end{equation}
if $\Delta t^n = \mathcal{O}(h)$.
\end{remark}

\begin{proposition}
\label{Prop: Simple lambda}
\everymath{\normalfont}
    Under the hypothesis that there exist $\underline{\varrho}, \overline{\varrho}, \overline{E} \in \mathbb{R}^+$ such that $0 < \underline{\varrho} \leq \varrho_h \leq \overline{\varrho}$ and $0 < \frac{1}{2}\varrho_h\modL{\mathbf{u}_h}^2 + P\paraL{\varrho_h} \leq \overline{E}$, the following hold:
    \begin{enumerate}
        \item There exist upper bounds $\overline{u}, \overline{m} \in \mathbb{R}^+$ such that $|\mathbf{u}_h| \leq \overline{u}$ and $|\varrho_h\mathbf{u}_h| \leq \overline{m}$,
        \item There exist bounds $\underline{P''}, \overline{P''} \in \mathbb{R}^+$ such that $\underline{P''} \leq P''(\varrho_h) \leq \overline{P''}$,
        \item If $\lambda_{\sigma} \geq \frac{\overline{P''}\overline{u}}{2\underline{P''}}$, then the condition \eqref{lambda cond IE} holds for $\diff{\varrho_h^{n+1}}_{\sigma} \neq 0$,
        \item Under the hypothesis that there exists $\underline{s} \in \mathbb{R}^+$ such that $\modL{\diff{\mathbf{u}_h^{n+1}}_{\sigma}} \geq \underline{s}$ whenever $\mathbf{u}$ is not smooth (i.e., when $D_x \mathbf{u}$ is not defined), the following is satified: \\
        If $\lambda_{\sigma} \geq 2\frac{\overline{\varrho} \ \overline{u}^2}{\underline{\varrho} \ \underline{s}}$, the condition \eqref{lambda cond KE} holds when $\mathbf{u}$ is not smooth. 
    \end{enumerate}
\end{proposition}
\begin{proof}
    \begin{enumerate}
        \item From $0 < \frac{1}{2}\varrho_h\modL{\mathbf{u}_h}^2 + P\paraL{\varrho_h} \leq \overline{E}$, it is straightforward to see that $|\mathbf{u}_h|$ has an upper bound, since the pressure potential is positive and density has a lower bound.
        \item Since $P(\varrho_h)$ is a polynomial function and since density has upper and lower bounds, there exist upper and lower bounds on $P''(\varrho_h)$. 
        \item $\lambda_{\sigma} \geq \frac{\overline{P''}\overline{u}}{2\underline{P''}} \geq \frac{ \diff{P''\paraL{\varrho^{n+1}_{*}} \mathbf{u}_h^{n+1}}_{\sigma} \cdot \mathbf{n}_{\sigma}}{4P''\left(\varrho_\dagger^{n+1}\right)} = \frac{ \diff{\varrho_h^{n+1}}_{\sigma}^2\diff{P''\paraL{\varrho^{n+1}_{*}} \mathbf{u}_h^{n+1}}_{\sigma} \cdot \mathbf{n}_{\sigma}}{4\diff{P'\left(\varrho_h^{n+1}\right)}_{\sigma}\diff{\varrho_h^{n+1}}_{\sigma}}$ for $\varrho_{*,K}^{n+1},\varrho_{*,L}^{n+1} , \varrho_\dagger^{n+1} \in [ \varrho_K^{n+1}, \varrho_L^{n+1} ]$ and $\sigma = K|L$.
        \item $\lambda_{\sigma} \geq 2\frac{\overline{\varrho} \ \overline{u}^2}{\underline{\varrho} \ \underline{s}} \geq  \frac{ \modL{\avg{\varrho_h^{n+1} \mathbf{u}_h^{n+1}}_{\sigma} \cdot \mathbf{n}_{\sigma}  \avg{\mathbf{u}_h^{n+1}}_{\sigma} }  + \modL{\avg{\varrho_h^n \mathbf{u}_h^n \mathbf{u}_h^n}_{\sigma} \cdot \mathbf{n}_{\sigma}}}{\avg{\varrho_h^{n+1}}_{\sigma} \modL{\diff{\mathbf{u}_h^{n+1}}_{\sigma}}}$ when $\mathbf{u}$ is not smooth.
    \end{enumerate}
\end{proof}

\begin{theorem}[Discrete total energy inequality]
\everymath{\normalfont}
    Let $\varrho_h^n>0 \in Q_h$, $\mathbf{u}_h^n \in \boldsymbol{Q}_h$, $\varrho_h^n|_K$, $\mathbf{u}_h^n|_K$ be solution to \eqref{Num mass}-\eqref{Num mom} for $1\leq n \leq N$. Let $\varrho_h^0$ and $\mathbf{u}_h^0$ be the initial data given by \eqref{IC} for each $K\in \mathcal{T}_h$. We suppose that periodic boundary conditions hold. Let $\lambda_{\sigma}$ satisfy 
    \begin{multline}
    \label{lambda cond TE}
        \lambda_{\sigma} = \max \left\{ 1_{\diff{\varrho_h^{n+1}}_{\sigma} \neq 0} \frac{ \diff{\varrho_h^{n+1}}_{\sigma}^2\diff{P''\paraL{\varrho^{n+1}_{*}} \mathbf{u}_h^{n+1}}_{\sigma} \cdot \mathbf{n}_{\sigma}}{4\diff{P'\left(\varrho_h^{n+1}\right)}_{\sigma}\diff{\varrho_h^{n+1}}_{\sigma}}, \right. \\ \left. 1_{\diff{\mathbf{u}_h^{n+1}}_{\sigma} \neq 0} \frac{ -\avg{\varrho_h^{n+1} \mathbf{u}_h^{n+1}}_{\sigma} \cdot \mathbf{n}_{\sigma}  \diffL{\frac{\modL{\mathbf{u}_h^{n+1}}^2}{2}}_{\sigma}   + \diff{\mathbf{u}_h^{n+1}}_{\sigma} \cdot \avg{\varrho_h^n \mathbf{u}_h^n \mathbf{u}_h^n}_{\sigma} \cdot \mathbf{n}_{\sigma}}{\avg{\varrho_h^{n+1}}_{\sigma} \modL{\diff{\mathbf{u}_h^{n+1}}_{\sigma}}^2} \right\},
    \end{multline}
    with $P''\paraL{\varrho^{n+1}_{*,K}}$, $P''\paraL{\varrho^{n+1}_{*,L}}$ given by \eqref{IE_CL},\eqref{IE_CK}, for each $\sigma = K|L$. $1_{\diff{\varrho_h^{n+1}}_{\sigma} \neq 0}$, $1_{\diff{\mathbf{u}_h^{n+1}}_{\sigma} \neq 0}$ represent the characteristic functions of the sets of faces satisfying $\diff{\varrho_h^{n+1}}_{\sigma} \neq 0$, $\diff{\mathbf{u}_h^{n+1}}_{\sigma} \neq 0$ respectively. Then, the following discrete total energy inequality holds:
    \begin{equation}
    \label{TE_ineq}
        \sumK{} \paraL{\frac{1}{2} \varrho_K^{n+1} \modL{\mathbf{u}_K^{n+1}}^2 + \frac{1}{\varepsilon^2} P\paraL{\varrho_K^{n+1}}} \leq \sumK{} \paraL{\frac{1}{2} \varrho_K^{n} \modL{\mathbf{u}_K^{n}}^2 + \frac{1}{\varepsilon^2} P\paraL{\varrho_K^{n}}}.
    \end{equation}
\end{theorem}
\begin{proof}
    By summing $\frac{1}{\varepsilon^2}$ \eqref{IE_ineq} and \eqref{KE_ineq}, we obtain \eqref{TE_ineq}. 
\end{proof}

\begin{remark}
    In view of \Cref{Prop: Simple lambda}, the condition in \eqref{lambda cond TE} is satisfied if $\lambda_{\sigma} = \max \left\{ \frac{\overline{P''}\overline{u}}{2\underline{P''}}, 2\frac{\overline{\varrho} \ \overline{u}^2}{\underline{\varrho} \ \underline{\text{d}u}} \right\}$ when $\mathbf{u}$ is not smooth. When $\mathbf{u}$ is smooth, we observed from \eqref{Simple lam smooth IE}, \eqref{Simple lam smooth KE} that $\lambda_\sigma = \mathcal{O}(1)$ if $\Delta t^n = \mathcal{O}(h)$. 
\end{remark}

\subsection{Asymptotic-preserving property}
In the previous subsections, under a specific condition on $\lambda_{\sigma}$ for each $\sigma = K|L$, we have proved the positivity of density and the discrete energy inequality for an arbitrary fixed mesh $h_K$ and time step $\Delta t^n$, $1 \leq n \leq N$ and $K \in \mathcal{T}_h$. In this subsection, we prove the asymptotic consistency of the numerical scheme \eqref{Num mass}-\eqref{Num mom}. Let us point out that the equations \eqref{Num mass}-\eqref{Num mom} are not directly used for implementation, and hence we first obtain the simplified expressions for $\varrho_K^{n+1}$ and $\varrho_K^{n+1}\mathbf{u}_K^{n+1}$ by linearizing the pressure. We consider a fixed uniform structured mesh with $h=h_K$ for all $K \in \mathcal{T}_h$ to prove the asymptotic consistency property. We also choose a fixed $\lambda$ for all $\sigma \in \mathcal{E}$ such that, $\lambda = \max_{\sigma \in \mathcal{E}} \lambda_{\sigma}$. 
\label{Subsec: AP}
\subsubsection{Implementation of the scheme}
\label{Subsec: Implementation}
In this subsection, we provide a strategy to implement \eqref{Num mass}-\eqref{Num mom} in only two steps. 
For this, we consider the following linearization of  pressure $p(\varrho_h^{n+1})$ about $p(\varrho_0)$:
\begin{equation}
\label{pressure linearization}
    p(\varrho_h^{n+1})\approx p(\varrho_0) + p'(\varrho_0) \paraL{\varrho_h^{n+1} - \varrho_0}. 
\end{equation}
Employing this linearization, \eqref{Num mass}-\eqref{Num mom} reads as:
\begin{gather}
    \label{Num mass simp}
    \frac{\varrho_{h}^{n+1}-\varrho_{h}^{n}}{\Delta t^n} + \divh{} \paraL{\varrho_h^{n+1}\mathbf{u}_h^{n+1}} - h \lambda \Delta_{h}\varrho_h^{n+1}  = 0,  \\
\label{Num mom simp}
    \frac{\varrho_{h}^{n+1} \mathbf{u}_{h}^{n+1}-\varrho_{h}^{n}\mathbf{u}_{h}^{n}}{\Delta t^n} + \divh{} \paraL{\varrho_h^{n} \mathbf{u}_h^{n} \otimes \mathbf{u}_h^{n}} -  h  \lambda \Delta_{h} \paraL{\varrho_h^{n+1}\mathbf{u}_h^{n+1}} + \frac{1}{\varepsilon^2} p'(\varrho_0)\nabla_h \varrho_h^{n+1} = \mathbf{0}.
\end{gather}
Note that \eqref{Num mass simp}$|_{K}$ is the same as \eqref{Num mass} but on a structured grid with fixed $\lambda$, whereas \eqref{Num mom simp}$|_{K}$ is an approximation of \eqref{Num mom} 
due to the linearization of pressure. From \eqref{Num mom simp}, we obtain the update equation for momentum as:
\begin{equation}
    \label{Update mom}
     \varrho_{h}^{n+1}\mathbf{u}_{h}^{n+1} = \paraL{I - \Delta t^n h \lambda\Delta_h}^{-1} \paraL{\varrho_{h}^{n}\mathbf{u}_{h}^{n}- \Delta t^n \divh{} \paraL{\varrho_h^{n} \mathbf{u}_h^{n} \otimes \mathbf{u}_h^{n}} - \frac{\Delta t^n}{\varepsilon^2} p'(\varrho_0)\nabla_h \varrho_h^{n+1} }, 
\end{equation}
 where $I$ is the identity operator. Then, we obtain the following update equation for mass by inserting \eqref{Update mom} into \eqref{Num mass simp}:
\begin{multline}
    \label{Update mass}
    \varrho_{h}^{n+1} = \paraL{I - \Delta t^n h \lambda\Delta_h - \paraL{\frac{\Delta t^n}{\varepsilon}}^2 p'(\varrho_0) \divh{} \paraL{I - \Delta t^n h \lambda\Delta_h}^{-1} \nabla_h }^{-1} \\
    \paraL{\varrho_{h}^n - \Delta t^n \divh{} \paraL{ \paraL{I - \Delta t^n h \lambda\Delta_h}^{-1} \paraL{\varrho_{h}^{n}\mathbf{u}_{h}^{n} - \Delta t^n \divh{} \paraL{\varrho_h^{n} \mathbf{u}_h^{n} \otimes \mathbf{u}_h^{n}}}}}. 
\end{multline}
Thus, the equations \eqref{Update mass}-\eqref{Update mom} are used to compute $\paraL{\varrho_h^{n+1},\varrho_h^{n+1} \mathbf{u}_h^{n+1}}$ directly in two steps. Let us note that the system \eqref{Update mass}-\eqref{Update mom} is equivalent to the system \eqref{Num mass simp}-\eqref{Num mom simp}. 

\subsubsection{Asymptotic consistency}
In this subsection, we prove that the numerical scheme \eqref{Update mass}-\eqref{Update mom} yields a consistent approximation of the limiting equations \eqref{Lim den vel exp}, \eqref{Lim div free vel}-\eqref{Lim vel}, in the limit as $\varepsilon \to 0$. 
\par Let us observe that the matrices involved in our numerical scheme are circulant (in one dimension) or block circulant with circulant blocks (BCCB) (in multiple dimensions), as the grid is uniform and the boundaries are periodic. All circulant or BCCB matrices commute (with each other) and are diagonalizable with the eigenbasis as the set of multidimensional Fourier modes dependent on the matrix size, but independent of its elements. Let us further note the following properties. 
\begin{remark}
\everymath{\normalfont}
    \begin{enumerate}
        \item The discrete Laplacian operator $\Delta_h$ is circulant or BCCB matrix of size $\mathcal{N}^d \times \mathcal{N}^d:= N_1N_2\dots N_d \times N_1N_2\dots N_d$. Here, $N_i$ represents the number of cells along direction $i$, for $i \in \{ 1,2,\dots, d\}$. 
        $\Delta_h$ is symmetric and its eigenvalues are $\mu_j \leq 0$, $j \in \{ 1,2,\dots, \mathcal{N}^d\}$.  
        \item The discrete divergence and discrete gradient operators are 
        \begin{equation}
            \divh{} = \begin{bmatrix} D_1 & D_2 & \dots D_d \end{bmatrix}, \  \nabla_h = \begin{bmatrix} D_1 & D_2 & \dots D_d \end{bmatrix}^T.
        \end{equation} 
        Here, $D_i$ for $i\in \{ 1,2,\dots, d\}$ are circulant or BCCB matrices that occur due to central differencing. Each $D_i$ is of size $\mathcal{N}^d \times \mathcal{N}^d$, is antisymmetric, and is diagonalizable. $D_i^2$ are symmetric circulant or BCCB matrices of size $\mathcal{N}^d \times \mathcal{N}^d$, and hence $\divh{} \nabla_h$ yields a symmetric circulant or BCCB matrix ${\Delta}_h = \sum_{i=1}^d D_i^2 $ of size $\mathcal{N}^d \times \mathcal{N}^d$. 
         \item The nullspaces of $D_i$ and $D_i^2$ are identical, $i\in \{ 1,2,\dots, d\}$. This is because the eigenbases of $D_i$ and $D_i^2$ are the same (as $D_i$ is circulant), and the eigenvalues of $D_i$ and $D_i^2$ are $\eta_j$ and $\eta_j^2$ respectively, for $j \in \{1,2,\dots,\mathcal{N}^d\}$. Thus,
        \begin{equation}
            Ker({\Delta}_h) = Ker(\nabla_h).
        \end{equation}
        \item  The matrix notation $\paraL{I - \Delta t^n h \lambda\Delta_h}$ occurring in the momentum update equation \eqref{Update mom} is simply a block diagonal matrix of size $d\mathcal{N}^d \times d\mathcal{N}^d := dN_1N_2\dots N_d \times dN_1N_2\dots N_d$ with identical blocks of the form $\paraL{I - \Delta t^n h \lambda\Delta_h}$, each of size $\mathcal{N}^d \times \mathcal{N}^d$. These blocks form symmetric circulant or BCCB matrices, and the eigenvalues of these blocks are $1-\Delta t^n h \lambda \mu_j \geq 1$, for $j \in \{1,2,\dots,\mathcal{N}^d\}$. Hence, these blocks are invertible, and the inverse is symmetric. Thus,
        \begin{equation}
            \divh{}\Large|\normalsize_{\mathcal{N}^d \times d \mathcal{N}^d} \paraL{I - \Delta t^n h \lambda\Delta_h}^{-1}\Huge|\normalsize_{d\mathcal{N}^d \times d\mathcal{N}^d} \nabla_h\Large|\normalsize_{d \mathcal{N}^d \times \mathcal{N}^d} = \paraL{I - \Delta t^n h \lambda\Delta_h}^{-1}\Huge|\normalsize_{\mathcal{N}^d \times \mathcal{N}^d} {\Delta}_h\Large|\normalsize_{\mathcal{N}^d \times \mathcal{N}^d}
        \end{equation}
        where ${\Delta}_h = \divh{} \nabla_h$.
        \item The matrix $\paraL{I - \Delta t^n h \lambda\Delta_h - \paraL{\frac{\Delta t^n}{\varepsilon}}^2 p'(\varrho_0) \divh{} \paraL{I - \Delta t^n h \lambda\Delta_h}^{-1} \nabla_h }$ in the mass update equation \eqref{Update mass} is a circulant or BCCB matrix of size $\mathcal{N}^d \times \mathcal{N}^d$, and it becomes: 
        \begin{equation}
        \label{Matrix mass update}
            \paraL{I - \Delta t^n h \lambda\Delta_h - \paraL{\frac{\Delta t^n}{\varepsilon}}^2 p'(\varrho_0) \paraL{I - \Delta t^n h \lambda\Delta_h}^{-1} {\Delta}_h }.
        \end{equation}
        This matrix is symmetric since $\paraL{I - \Delta t^n h \lambda\Delta_h}^{-1}$ and ${\Delta}_h$ commute (as they are circulant or BCCB). Its eigenvalues are $1-\Delta t^n h \lambda \mu_j - \paraL{\frac{\Delta t^n}{\varepsilon}}^2 p'(\varrho_0) \frac{\mu_j}{1-\Delta t^n h \lambda \mu_j } \geq 1$, for $j \in \{1,2,\dots,\mathcal{N}^d\}$. Hence, the matrix \eqref{Matrix mass update} is invertible, and the inverse is symmetric. 
    \end{enumerate}
\end{remark}

\begin{lemma}
\everymath{\normalfont}
\label{Lem: AP projection}
    The matrix \eqref{Matrix mass update} has the following projection property:
    \begin{equation}
    \label{Proj inv matrix mass}
        \paraL{I - \Delta t^n h \lambda\Delta_h - \paraL{\frac{\Delta t^n}{\varepsilon}}^2 p'(\varrho_0) \paraL{I - \Delta t^n h \lambda\Delta_h}^{-1} {\Delta}_h }^{-1} = \pi_{Ker(\Delta_h)} + \mathcal{O}\paraL{\varepsilon^2},
    \end{equation}
    where  $\pi_{Ker(\Delta_h)}$ is the projection onto $Ker(\Delta_h)$ along $\mathcal{R}(\Delta_h)$; $Ker(\Delta_h)$, $\mathcal{R}(\Delta_h)$ denote the kernel and the range of $\Delta_h$. 
\end{lemma}
\begin{proof}
    The inverse of \eqref{Matrix mass update} is diagonalizable, and its eigenvalues are
    \begin{equation}
        \frac{\varepsilon^2}{\varepsilon^2 \paraL{1 - \Delta t^n h \lambda \mu_j}-( \Delta t_n)^2 p'(\varrho_0)\frac{\mu_j}{1 - \Delta t^n h \lambda \mu_j}} = \left\{ \begin{matrix} 1 & \text{if } \mu_j=0 \\ \mathcal{O}\paraL{\varepsilon^2} & \text{ if } \mu_j\neq 0
        \end{matrix} \right..
    \end{equation}
    for $j \in \{1,2,\dots,\mathcal{N}^d\}$, and $\mu_j$ are the eigenvalues of $\Delta_h$. 
\end{proof}

\begin{lemma}
\everymath{\normalfont}
\label{Lem: AP exp}
    Let $\varrho_h^n = \varrho_0 + \mathcal{O}\paraL{\varepsilon^2}$ and $ \mathbf{u}_{h}^n=\mathcal{O}(1)$. Then the numerical solution to the linear system \eqref{Update mass}-\eqref{Update mom} satisfies $\varrho_h^{n+1} = \varrho_0 + \mathcal{O}\paraL{\varepsilon^2}$, $\mathbf{u}_{h}^{n+1} = \mathcal{O}(1)$, and $\divh{} \mathbf{u}_h^{n+1} =  \mathcal{O}\paraL{\varepsilon^2}$.
\end{lemma}
\begin{proof}
    The dimensions of the kernel of $\Delta_h$ and $\Delta_h^T$ are equal, \textit{i.e.}, $Ker\paraL{\Delta_h^T} = Ker\paraL{\Delta_h}$. Thus, $Ker\paraL{\Delta_h}^{\perp} = \mathcal{R}\paraL{\Delta_h}$, and the projection in \eqref{Proj inv matrix mass} is the orthogonal projection onto $Ker\paraL{\Delta_h}$. We also have,
    \begin{equation}
        \mathcal{R}\paraL{\Delta_h} = Ker\paraL{\Delta_h}^{\perp} = Ker\paraL{\nabla_h}^{\perp} = \mathcal{R}\paraL{\nabla_h^T} = \mathcal{R}\paraL{\divh{}}. 
    \end{equation}
    Further, inserting $\varrho_h^n = \varrho_0 + \mathcal{O}\paraL{\varepsilon^2}$ into \eqref{Update mass} and noting that $\varrho_0$ is in the kernel of \eqref{Matrix mass update}, we obtain,
    \begin{multline}
    \varrho_{h}^{n+1} - \varrho_0 = \paraL{I - \Delta t^n h \lambda\Delta_h - \paraL{\frac{\Delta t^n}{\varepsilon}}^2 p'(\varrho_0) \paraL{I - \Delta t^n h \lambda\Delta_h}^{-1} \Delta_h }^{-1} \\
    \paraL{ - \Delta t^n \divh{} \paraL{ \paraL{I - \Delta t^n h \lambda\Delta_h}^{-1} \paraL{\varrho_{h}^{n}\mathbf{u}_{h}^{n} - \Delta t^n \divh{} \paraL{\varrho_h^{n} \mathbf{u}_h^{n} \otimes \mathbf{u}_h^{n}}}} + \mathcal{O}\paraL{\varepsilon^2}}. 
    \end{multline}
    Therefore, $ \varrho_{h}^{n+1} - \varrho_0 = \mathcal{O}\paraL{\varepsilon^2}$ due to \Cref{Lem: AP projection}. 
    \par Further, inserting $\varrho_h^n = \varrho_0 + \varepsilon^2 \varrho_{h_2}^n$, $\varrho_h^{n+1} = \varrho_0 + \varepsilon^2 \varrho_{h_2}^{n+1}$ and $ \mathbf{u}_{h}^n=\mathcal{O}(1) = \mathbf{u}_{h_0}^n + \mathcal{O}(\varepsilon)$ into \eqref{Update mom}, we obtain
    \begin{equation}
        \varrho_h^{n+1}\mathbf{u}_{h}^{n+1} = \paraL{I - \Delta t^n h \lambda\Delta_h}^{-1} \paraL{ \varrho_0\mathbf{u}_{h_0}^{n}- \Delta t^n \varrho_0\divh{} \paraL{ \mathbf{u}_{h_0}^{n} \otimes \mathbf{u}_{h_0}^{n}} - \Delta t^np'(\varrho_0) \nabla_h \varrho_{h_2}^{n+1} } + \mathcal{O}(\varepsilon),
    \end{equation}
    and thus $\mathbf{u}_{h}^{n+1} = \mathcal{O}(1)$. Hence, using this fact and inserting $\varrho_h^n = \varrho_0 + \mathcal{O}\paraL{\varepsilon^2}$, $\varrho_h^{n+1} = \varrho_0 + \mathcal{O}\paraL{\varepsilon^2}$ into \eqref{Num mass simp}, we obtain $\divh{} \mathbf{u}_h^{n+1} =  \mathcal{O}\paraL{\varepsilon^2}$.
\end{proof}

\begin{theorem}[Asymptotic-preserving property]
\everymath{\normalfont}
    Let the initial data satisfy $\varrho_h^0 = \varrho_0 + \varepsilon^2 \varrho_{h_2}^0$ and $\mathbf{u}_{h}^0=\mathcal{O}(1)$. Then, the numerical solution to the linear system \eqref{Update mass}-\eqref{Update mom} satisfies $\varrho_h^n = \varrho_0 + \varepsilon^2 \varrho_{h_2}^n$, $\mathbf{u}_{h}^n=\mathcal{O}(1)$, and $\divh{} \mathbf{u}_{h}^n = \mathcal{O}\paraL{\varepsilon^2}$, for $1\leq n \leq N$. \\
    Further, in the limit as $\varepsilon \to 0$, the momentum update equation \eqref{Update mom} becomes
    \begin{equation}
        \mathbf{u}_{h_0}^{n+1} = \mathbf{u}_{h_0}^{n}- \Delta t^n \divh{} \paraL{\mathbf{u}_{h_0}^{n} \otimes \mathbf{u}_{h_0}^{n}} - \Delta t^n\frac{p'(\varrho_0)}{\varrho_0} \nabla_h \varrho_{h_2}^{n+1} + \Delta t^n h \lambda\Delta_h \mathbf{u}_{h_0}^{n+1},
    \end{equation}
    which is a consistent numerical approximation for \eqref{Lim vel}, if there is enough regularity of the solution variables. 
\end{theorem}
\begin{proof}
    Due to \Cref{Lem: AP exp} and initial data $\varrho_h^0 = \varrho_0 + \varepsilon^2 \varrho_{h_2}^0$, $\mathbf{u}_{h}^0=\mathcal{O}(1)$, we have by induction $\varrho_h^n = \varrho_0 + \varepsilon^2 \varrho_{h_2}^n$, $\mathbf{u}_{h}^n=\mathcal{O}(1)$, and $\divh{} \mathbf{u}_{h}^n = \mathcal{O}\paraL{\varepsilon^2}$, for $1\leq n \leq N$. By inserting these into \eqref{Update mom} and taking the limit $\varepsilon \to 0$, we obtain
    \begin{equation}
    \label{Num lim vel}
        \mathbf{u}_{h_0}^{n+1} = \paraL{I - \Delta t^n h \lambda\Delta_h}^{-1} \paraL{\mathbf{u}_{h_0}^{n}- \Delta t^n \divh{} \paraL{ \mathbf{u}_{h_0}^{n} \otimes \mathbf{u}_{h_0}^{n}} - \Delta t^n\frac{p'(\varrho_0)}{\varrho_0} \nabla_h \varrho_{h_2}^{n+1} }.
    \end{equation} 
\end{proof}
Let us note that, in order to get asymptotic consistency, the initial data need not be well-prepared according to \eqref{Def: Well-prep ID}, \textit{i.e.}, $\divh{}\mathbf{u}_{h}^0$ need not be $\mathcal{O}\paraL{\varepsilon^2}$.

\begin{remark}
\everymath{\normalfont}
\label{Rem: pressure linearization}
    Since $\varrho_h^{n+1} = \varrho_0 + \mathcal{O}\paraL{\varepsilon^2}$ due to \Cref{Lem: AP exp}, the pressure linearization which we considered in \eqref{pressure linearization} holds upto $\mathcal{O}\paraL{\varepsilon^4}$, as 
    \begin{equation}
        p(\varrho_h^{n+1})= p(\varrho_0) + p'(\varrho_0) \paraL{\varrho_h^{n+1} - \varrho_0} + \frac{p''(\varrho_{\xi})}{2} \paraL{\varrho_h^{n+1} - \varrho_0}^2 = p(\varrho_0) + p'(\varrho_0) \paraL{\varrho_h^{n+1} - \varrho_0} + \mathcal{O}\paraL{\varepsilon^4}
    \end{equation}
    for $\varrho_{\xi} \in \left[\varrho_0,\varrho_h^{n+1}\right]$. Furthermore, by considering $\varrho_h^{n+1} = \varrho_0 + \varepsilon^2 \varrho_{h_2}^{n+1}$, we infer that \eqref{Num mom simp}$|_{K}$ approximates \eqref{Num mom} up to $\mathcal{O}\paraL{\varepsilon^2}$, as
    \begin{equation}
    \label{rem: p linearization eqn}
        \frac{1}{\varepsilon^2}\nabla_h p(\varrho_h)^{n+1} = p'(\varrho_0) \nabla_h \varrho_{h_2}^{n+1} + \varepsilon^2 \nabla_h \paraL{ \frac{p''(\varrho_{\xi})}{2} \paraL{\varrho_{h_2}^{n+1} }^2 } = \frac{1}{\varepsilon^2} p'(\varrho_0) \nabla_h \varrho_{h}^{n+1} + \mathcal{O}\paraL{\varepsilon^2} \text{ since } \nabla_h \varrho_0 = 0. 
    \end{equation}
\end{remark}

Due to the above \Cref{Rem: pressure linearization}, we have the following proposition on the satisfaction of discrete total energy inequality up to $\mathcal{O}\paraL{\varepsilon^2}$.

 \begin{proposition}[Discrete total energy inequality up to $\mathcal{O}\paraL{\varepsilon^2}$]
 \everymath{\normalfont}
    Let $\varrho_h^n>0 \in Q_h$, $\mathbf{u}_h^n \in \boldsymbol{Q}_h$ be solution to \eqref{Num mass simp}-\eqref{Num mom simp} for $1\leq n \leq N$. Let $\varrho_h^0$ and $\mathbf{u}_h^0$ be the initial data given by \eqref{IC} for each $K\in \mathcal{T}_h$, and satisfying $\varrho_h^0 = \varrho_0 + \mathcal{O}\paraL{\varepsilon^2}$, $\varrho_0 $ a constant, $\mathbf{u}_{h}^0=\mathcal{O}(1)$. We suppose that periodic boundary conditions hold. Let $\lambda$ satisfy $\lambda = \max_{\sigma \in \mathcal{E}} \lambda_{\sigma}$, with $\lambda_{\sigma}$ given by \eqref{lambda cond TE}. Then, \eqref{TE_ineq} holds up to $\mathcal{O}\paraL{\varepsilon^2}$, that is
    \begin{equation}
        \sumK{} \paraL{\frac{1}{2} \varrho_K^{n+1} \modL{\mathbf{u}_K^{n+1}}^2 + \frac{1}{\varepsilon^2} P\paraL{\varrho_K^{n+1}}} \leq \sumK{} \paraL{\frac{1}{2} \varrho_K^{n} \modL{\mathbf{u}_K^{n}}^2 + \frac{1}{\varepsilon^2} P\paraL{\varrho_K^{n}}} + \mathcal{O}\paraL{\varepsilon^2}.
    \end{equation}
\end{proposition}
\begin{proof}
    The discrete internal energy inequality \eqref{IE_ineq} remains unchanged as \eqref{Num mass simp}$|_{K}$ is the same as \eqref{Num mass} but on a structured grid with fixed $\lambda$. 
    \par Since the pressure term in the discrete kinetic energy inequality \eqref{KE_ineq} is approximated as,
    \begin{eqnarray}
        \sumK{} \mathbf{u}_K^{n+1} \cdot  \paraL{\nabla_h \varrho_h^{n+1}}_K \frac{p'(\varrho_0)}{\varepsilon^2}  &=& \frac{1}{\varepsilon^2}\sumK{} \mathbf{u}_K^{n+1} \cdot  \paraL{\nabla_h p\paraL{\varrho_h^{n+1}}}_K + \mathcal{O}\paraL{\varepsilon^2}, \text{by \eqref{rem: p linearization eqn}}\\
        &=& - \frac{1}{\varepsilon^2}\sumK{} p\paraL{\varrho_K^{n+1}} \paraL{\text{div}_h \mathbf{u}^{n+1}_h}_K  + \mathcal{O}\paraL{\varepsilon^2},
    \end{eqnarray}
    \eqref{KE_ineq} holds up to $\mathcal{O}\paraL{\varepsilon^2}$. Summing the discrete internal energy inequality and the discrete kinetic energy inequality, we infer that the discrete total energy inequality  \eqref{TE_ineq} is satisfied up to $\mathcal{O}\paraL{\varepsilon^2}$. 
\end{proof} 

\section{Numerical validation}
\label{Sec: Num results}
In this section, we present the validation of the proposed numerical method for standard benchmark problems. We observed in the previous section that if $\lambda$ satisfies a specific condition, the method is energy-stable. We also noted that these conditions imply that $\lambda$ is $\mathcal{O}(1)$ when the solution is smooth and that $\lambda = \max \left\{ \frac{\overline{P''}\overline{u}}{2\underline{P''}}, 2\frac{\overline{\varrho} \ \overline{u}^2}{\underline{\varrho} \ \underline{\text{d}u}} \right\}$ when the solution is non-smooth. Thus, in our numerical simulations, we simply choose $\lambda$ to be a constant of $\mathcal{O}(1)$. In particular, $\lambda=1$ results in energy stability for all the test cases presented below. For the Gresho and the travelling vortex problems, it was suitable to also choose $\lambda$ by keeping all the terms at the explicit time level in the condition \eqref{lambda cond TE} that we derived, i.e. 
\begin{multline}
    \label{lambda cond numerical}
        \lambda = c \ \max_{\sigma \in \mathcal{E}} \left\{ \max \left\{ 1_{\diff{\varrho_h^{n}}_{\sigma} \neq 0} \frac{ \diff{\varrho_h^{n}}_{\sigma}^2\diff{P''\paraL{\varrho^{n}_{*}} \mathbf{u}_h^{n}}_{\sigma} \cdot \mathbf{n}_{\sigma}}{4\diff{P'\left(\varrho_h^{n}\right)}_{\sigma}\diff{\varrho_h^{n}}_{\sigma}}, \right. \right. \\ \left. \left. 1_{\diff{\mathbf{u}_h^{n}}_{\sigma} \neq 0} \frac{ -\avg{\varrho_h^{n} \mathbf{u}_h^{n}}_{\sigma} \cdot \mathbf{n}_{\sigma}  \diffL{\frac{\modL{\mathbf{u}_h^{n}}^2}{2}}_{\sigma}   + \diff{\mathbf{u}_h^{n}}_{\sigma} \cdot \avg{\varrho_h^n \mathbf{u}_h^n \mathbf{u}_h^n}_{\sigma} \cdot \mathbf{n}_{\sigma}}{\avg{\varrho_h^{n}}_{\sigma} \modL{\diff{\mathbf{u}_h^{n}}_{\sigma}}^2} \right\} \right\}.
\end{multline} 
The factor $c$ in the above expression, and the values for $\lambda$ will be detailed in the corresponding subsections for the Gresho and the travelling vortex problems. \\
Further, we also note that the proofs in the previous section hold without imposing any restriction on the time step. Thus, we simply choose the time step as: $\Delta t^n = C \ h / \max_{x \in \Omega} |\mathbf{u}^n|$ for $0\leq n \leq N-1$, where $C$ is the CFL number and $h$ is the space discretization parameter corresponding to the structured uniform grid. 

\subsection{Standard periodic problem}
\label{Subsec: SPP}
\begin{figure}[h!]
\centering
\begin{subfigure}[b]{0.32\textwidth}
\centering
\includegraphics[width=\textwidth]{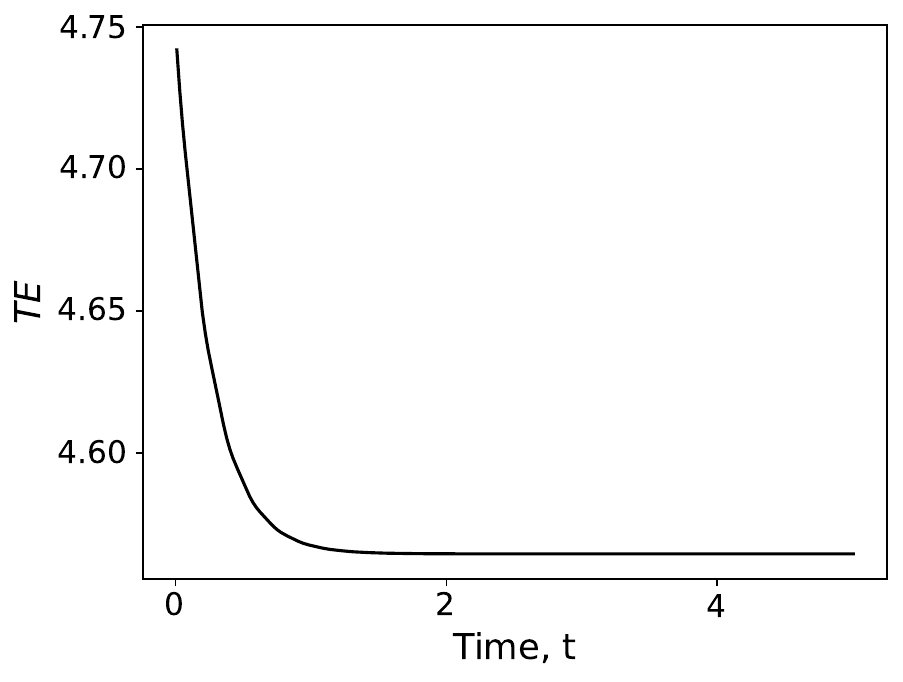}
\label{Fig:SPP_TE_eps05}
\end{subfigure}
\hspace{-0.2cm}
\begin{subfigure}[b]{0.33\textwidth}
\centering
\includegraphics[width=\textwidth]{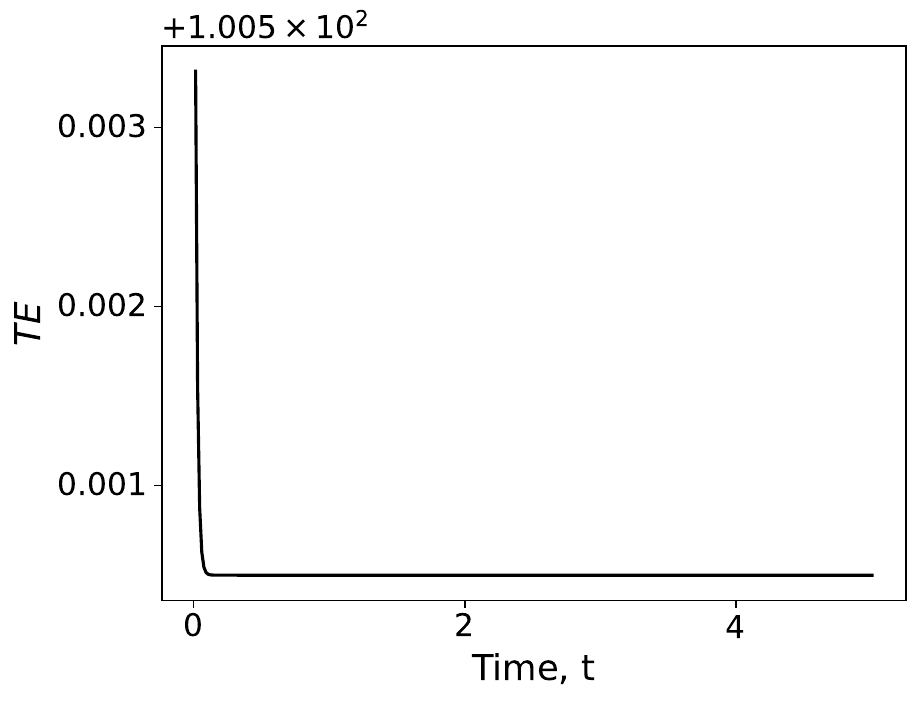}
\label{Fig:SPP_TE_eps01}
\end{subfigure}
\hspace{-0.2cm}
\begin{subfigure}[b]{0.31\textwidth}
\centering
\includegraphics[width=\textwidth]{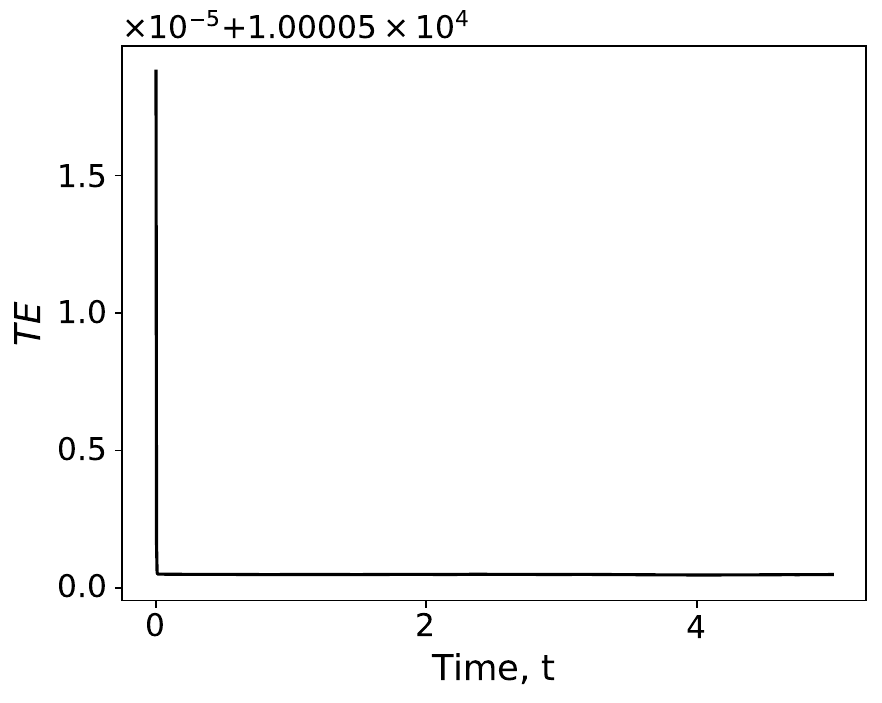}
\label{Fig:SPP_TE_eps001}
\end{subfigure}

\begin{subfigure}[b]{0.32\textwidth}
\centering
\includegraphics[width=\textwidth]{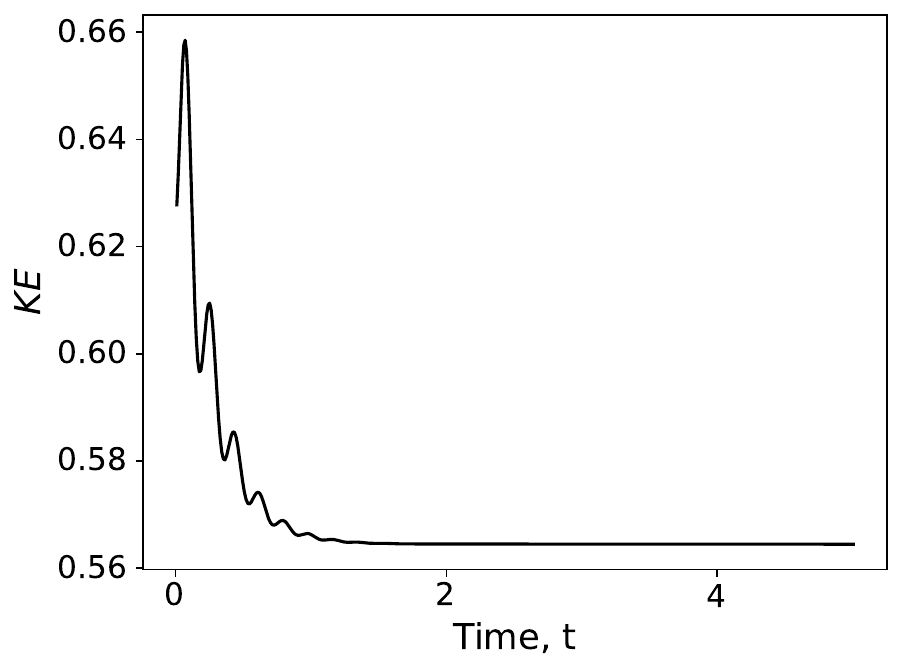}
\end{subfigure}
\hspace{-0.2cm}
\begin{subfigure}[b]{0.33\textwidth}
\centering
\includegraphics[width=\textwidth]{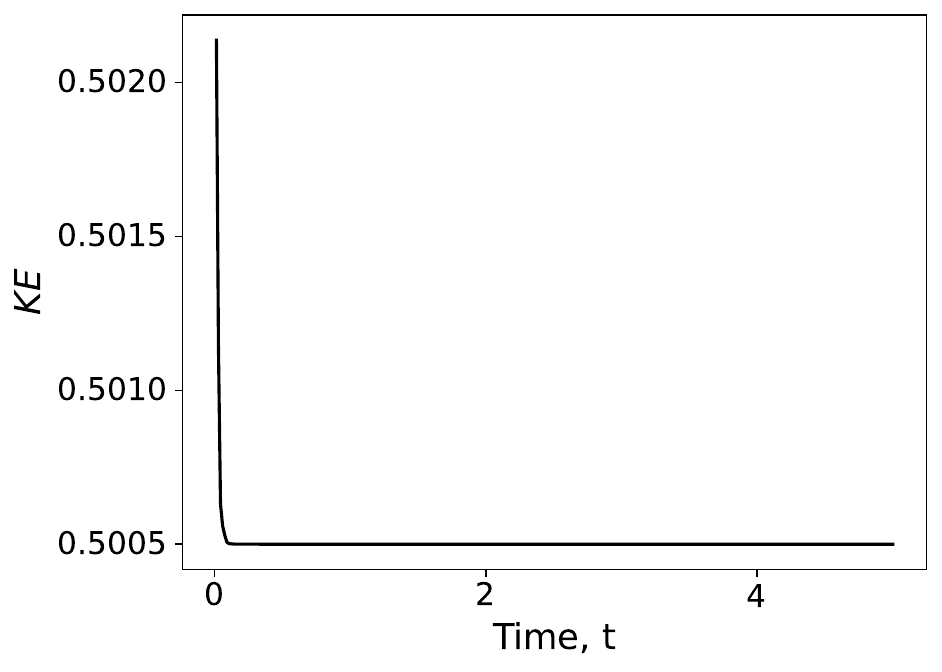}
\end{subfigure}
\hspace{-0.2cm}
\begin{subfigure}[b]{0.31\textwidth}
\centering
\includegraphics[width=\textwidth]{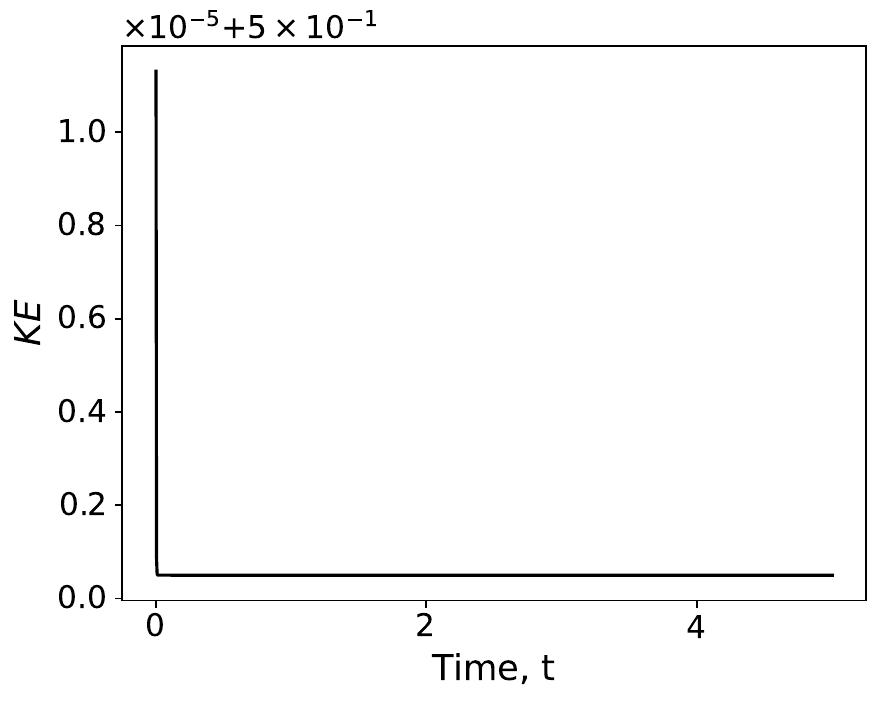}
\end{subfigure}

\begin{subfigure}[b]{0.32\textwidth}
\centering
\includegraphics[width=\textwidth]{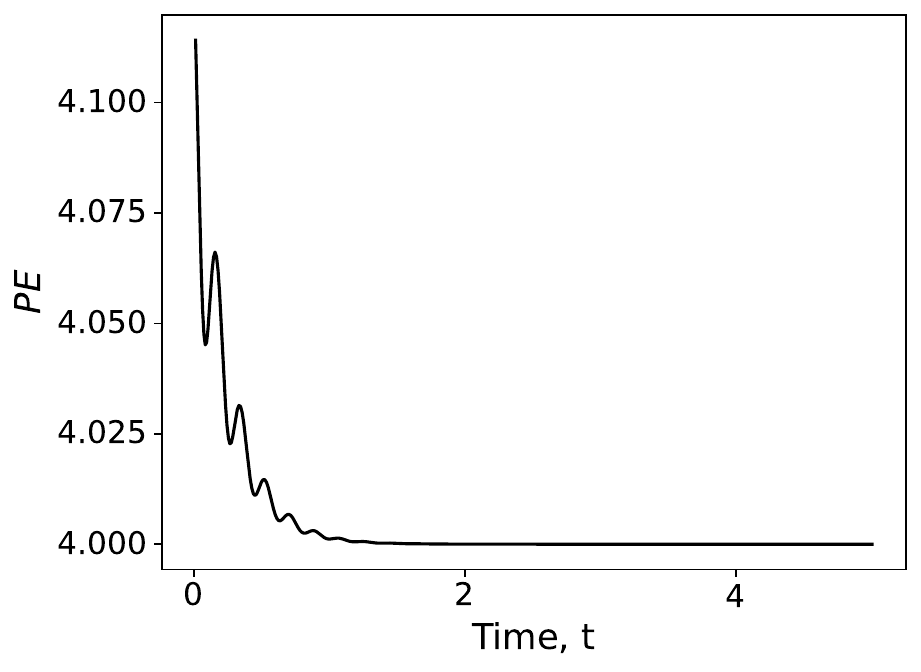}
\caption{$\varepsilon=0.5$}
\end{subfigure}
\hspace{-0.2cm}
\begin{subfigure}[b]{0.33\textwidth}
\centering
\includegraphics[width=\textwidth]{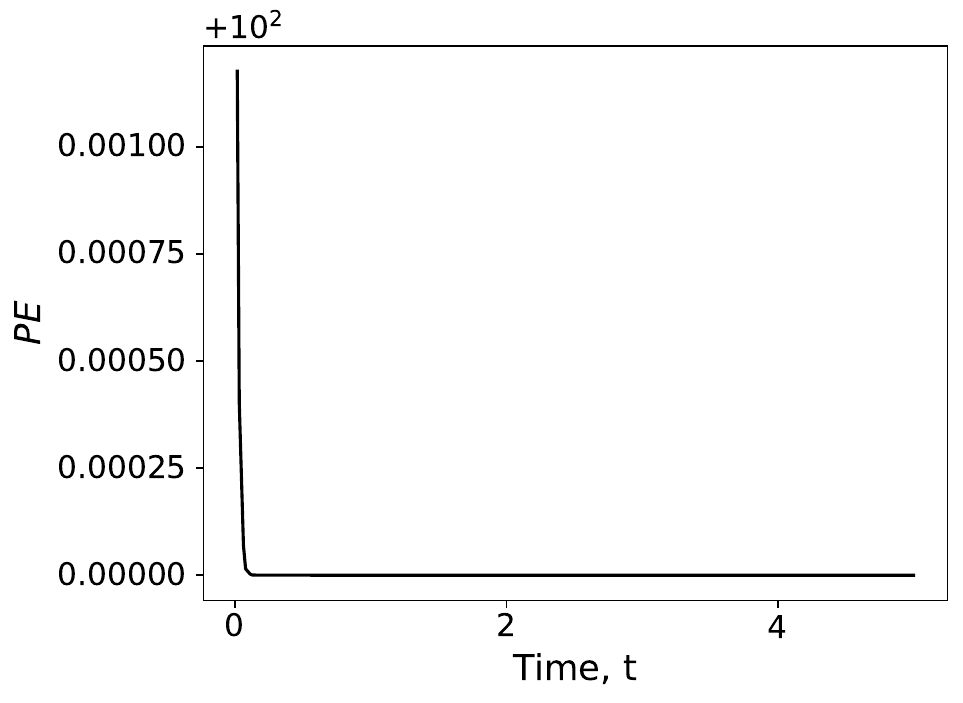}
\caption{$\varepsilon=0.1$}
\end{subfigure}
\hspace{-0.2cm}
\begin{subfigure}[b]{0.3\textwidth}
\centering
\includegraphics[width=\textwidth]{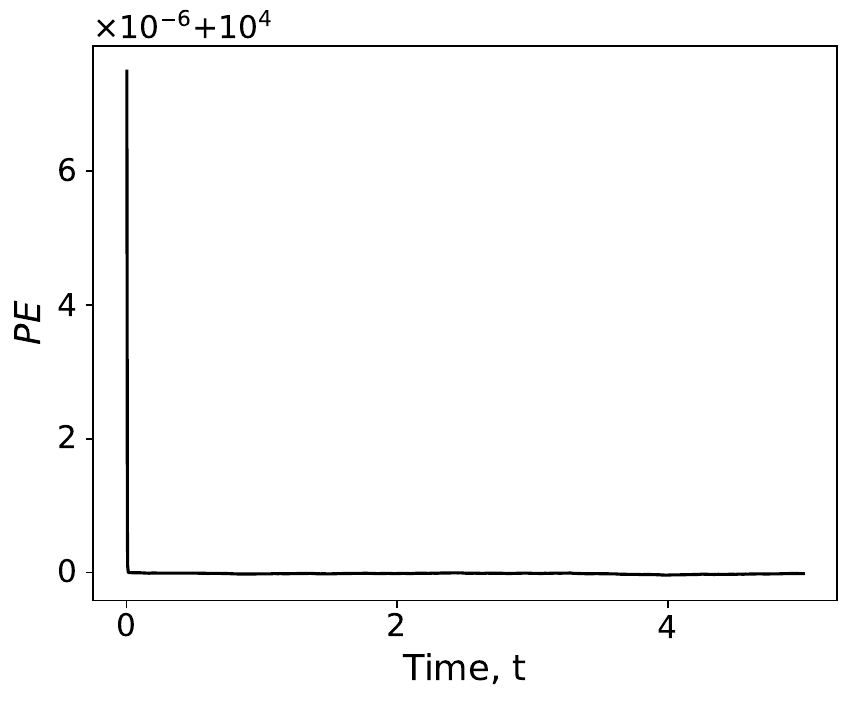}
\caption{$\varepsilon=0.01$}
\end{subfigure}

\caption{\centering \textbf{Example \ref{Subsec: SPP} - Standard periodic problem:} Top - Total energy; Middle - Kinetic energy; Bottom - Potential energy}
\label{Fig:SPP}
\end{figure}

This problem is defined on spatial domain $\Omega = [0,1]$, and the initial conditions are $\varrho^0(x)=1+\varepsilon^2 \sin \left(2\pi x\right)$, $u^0(x)=1+\varepsilon \sin \left( 2\pi x\right)$ for $x \in \Omega$. The boundary conditions are periodic. The parametric values are: $\kappa = 1$ and $\gamma = 2$, and we consider three different values of $\varepsilon$: $0.5$, $0.1$, and $0.01$. $\lambda=1$ is used as the numerical diffusion coefficient for all the values of $\varepsilon$. A CFL value of $C=0.8$ is used for $\varepsilon=0.5,0.1$, and $C=0.1$ is used for $\varepsilon=0.01$. We first present a convergence study of the proposed method \eqref{Update mass}-\eqref{Update mom} at time $t=0.1$ for $\varepsilon=0.5,0.1$, and at a smaller time $t=0.05$ for $\varepsilon=0.01$. For the discretization of $\Omega$, different numbers of grid points: N$_x = 20, 50, 100, 200, 250$ and $500$ was used. A reference of N$_x=1000$ is used in computing the $L^2$ error, and the convergence rates of $\varrho$ and $u$ are presented in \Cref{tab: SPP EOC_rho,tab: SPP EOC_u}, respectively. It can be seen that $\varrho$ converges with rates of about $2.2,1.2,1.1$, respectively, for $\varepsilon=0.5,0.1,0.01$ and $u$ converges with rates of about $1.7, 1.2, 1.3$, respectively, for $\varepsilon=0.5,0.1,0.01$. For $\varepsilon=0.01$, the errors in $\varrho$ and $u$ are of $\mathcal{O}(\varepsilon^2)$ and $\mathcal{O}(\varepsilon)$ respectively, revealing that the quantities $\varrho_0$ and $u_{0}$ in asymptotic expansion \eqref{Lim den vel exp} are constants. Furthermore, \Cref{Fig:SPP} shows the total, kinetic, and potential energy plots (with time) for N$_x=50$ with $\varepsilon=0.5,0.1$ and $0.01$, and it confirms that the fully-discrete method is energy-stable for different values of $\varepsilon$. It is interesting to observe that the oscillations in kinetic and potential energies are balancing each other to yield a decrease in total energy for $\varepsilon=0.5$. 

\begin{table}[h!]
    \centering
    \renewcommand{\arraystretch}{1.3} 
    \setlength{\tabcolsep}{10pt}      
    
    \begin{tabular}{|c|c|c|c|c|c|c|c|}
        \hline
        \multirow{2}{*}{\textbf{N$_x$}} & \multirow{2}{*}{$\mathbf{\Delta x}$} & 
        \multicolumn{2}{c|}{$\mathbf{\varepsilon=0.5}$} & 
        \multicolumn{2}{c|}{$\mathbf{\varepsilon=0.1}$} & 
        \multicolumn{2}{c|}{$\mathbf{\varepsilon=0.01}$} \\ 
        \cline{3-8}
        & & $\mathbf{||\varrho \textbf{ error}||_{L^2}}$ & \textbf{EOC} 
        & $\mathbf{||\varrho \textbf{ error}||_{L^2}}$ & \textbf{EOC} 
        & $\mathbf{||\varrho \textbf{ error}||_{L^2}}$ & \textbf{EOC} \\ 
        \hline
        20  & 0.05  & 0.04944  & -      & 0.00489  & -      & 4.25 $\times 10^{-5}$ & -      \\  
        50  & 0.02  & 0.02471  & 0.7568 & 0.00514  & -0.0539 & 4.01 $\times 10^{-5}$ & 0.0619 \\  
        100 & 0.01  & 0.01239  & 0.9956 & 0.00464  & 0.1463 & 3.30 $\times 10^{-5}$ & 0.2836 \\  
        200 & 0.005 & 0.01035  & 0.2603 & 0.00352  & 0.3987 & 2.35 $\times 10^{-5}$ & 0.4866  \\  
        250 & 0.004 & 0.00746  & 1.4644 & 0.00296  & 0.7757 & 1.81 $\times 10^{-5}$ & 1.1846  \\  
        500 & 0.002 & 0.00165  & 2.1766 & 0.00126  & 1.2346 & 8.59 $\times 10^{-6}$ & 1.0720 \\  
        \hline
    \end{tabular}
    
    \caption{\centering \textbf{Example \ref{Subsec: SPP} - Standard periodic problem:} Convergence rates of $L^2$ error in $\varrho$.}
    \label{tab: SPP EOC_rho}
\end{table}

\begin{table}[h!]
    \centering
    \renewcommand{\arraystretch}{1.3} 
    \setlength{\tabcolsep}{10pt}      
    
    \begin{tabular}{|c|c|c|c|c|c|c|c|}
        \hline
        \multirow{2}{*}{\textbf{N$_x$}} & \multirow{2}{*}{$\mathbf{\Delta x}$} & 
        \multicolumn{2}{c|}{$\mathbf{\varepsilon=0.5}$} & 
        \multicolumn{2}{c|}{$\mathbf{\varepsilon=0.1}$} & 
        \multicolumn{2}{c|}{$\mathbf{\varepsilon=0.01}$} \\ 
        \cline{3-8}
        & & $\mathbf{||u \textbf{ error}||_{L^2}}$ & \textbf{EOC} 
        & $\mathbf{||u \textbf{ error}||_{L^2}}$ & \textbf{EOC} 
        & $\mathbf{||u \textbf{ error}||_{L^2}}$ & \textbf{EOC} \\ 
        \hline
        20  & 0.05  & 0.22771  & -      & 0.05822  & -      & 8.31 $\times 10^{-3}$ & -      \\  
        50  & 0.02  & 0.10364  & 0.8590 & 0.05992  & -0.0315 & 7.98 $\times 10^{-3}$ & 0.0437 \\  
        100 & 0.01  & 0.05272  & 0.9752 & 0.05342  & 0.1657 & 6.40 $\times 10^{-3}$ & 0.3192 \\  
        200 & 0.005 & 0.02675  & 0.9786 & 0.03957  & 0.4332 & 4.37 $\times 10^{-3}$ & 0.5490  \\  
        250 & 0.004 & 0.01959  & 1.3966 & 0.03367  & 0.7230 & 3.55 $\times 10^{-3}$ & 0.9331  \\  
        500 & 0.002 & 0.00591  & 1.7296 & 0.01512  & 1.1554 & 1.46 $\times 10^{-3}$ & 1.2799 \\  
        \hline
    \end{tabular}
    
    \caption{\centering \textbf{Example \ref{Subsec: SPP} - Standard periodic problem:} Convergence rates of $L^2$ error in $u$.}
    \label{tab: SPP EOC_u}
\end{table}

\subsection{Colliding acoustic waves problem}
\label{Subsec:CAW}
\begin{figure}[h!]
\centering
\begin{subfigure}[b]{0.3\textwidth}
\centering
\includegraphics[width=\textwidth]{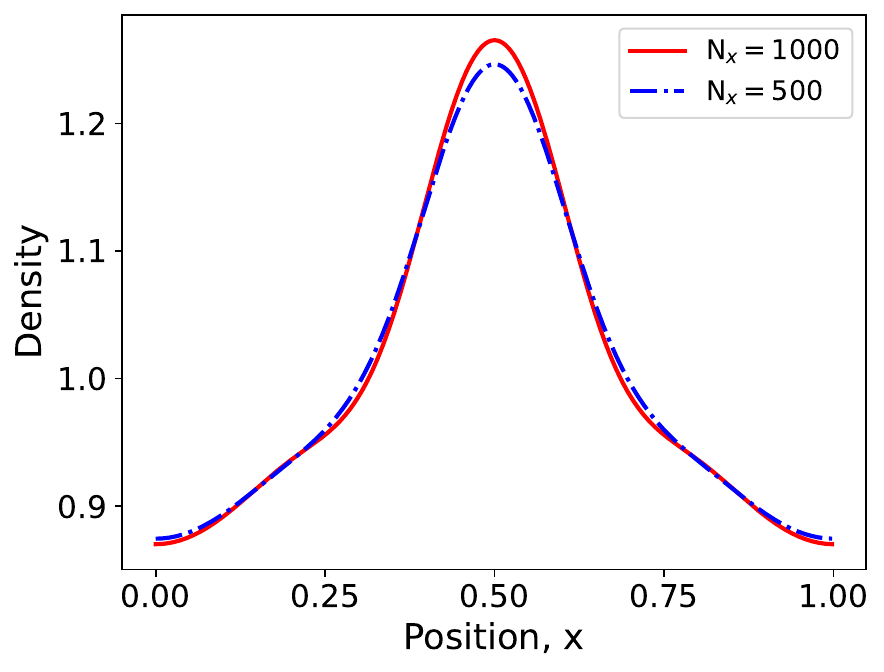}
\caption{$t=0.04$}
\label{Fig:CAW_t4_den}
\end{subfigure}
\hspace{-0.2cm}
\begin{subfigure}[b]{0.31\textwidth}
\centering
\includegraphics[width=\textwidth]{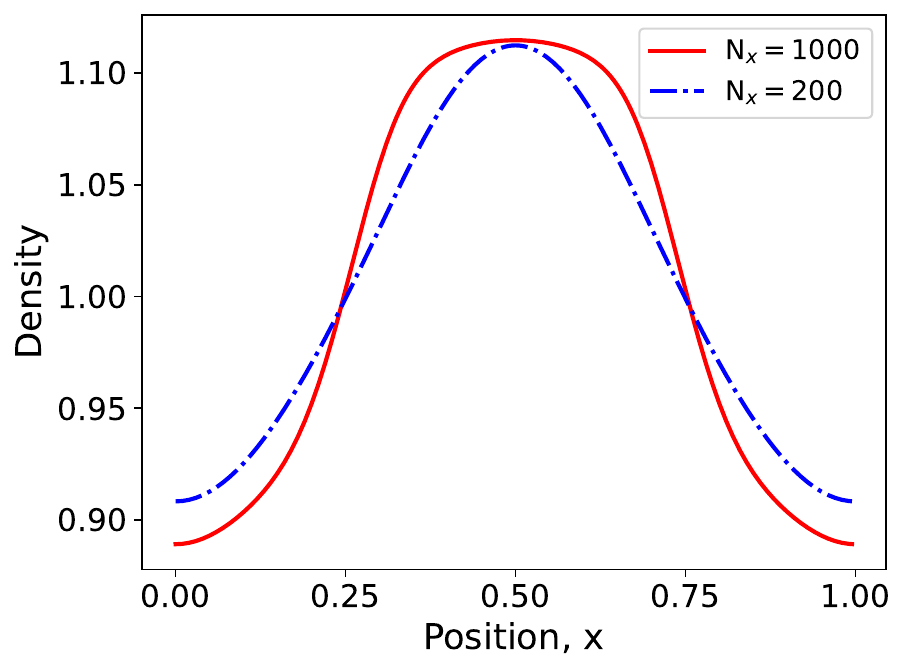}
\caption{$t=0.06$}
\label{Fig:CAW_t6_den}
\end{subfigure}
\hspace{-0.2cm} 
\begin{subfigure}[b]{0.31\textwidth}
\centering
\includegraphics[width=\textwidth]{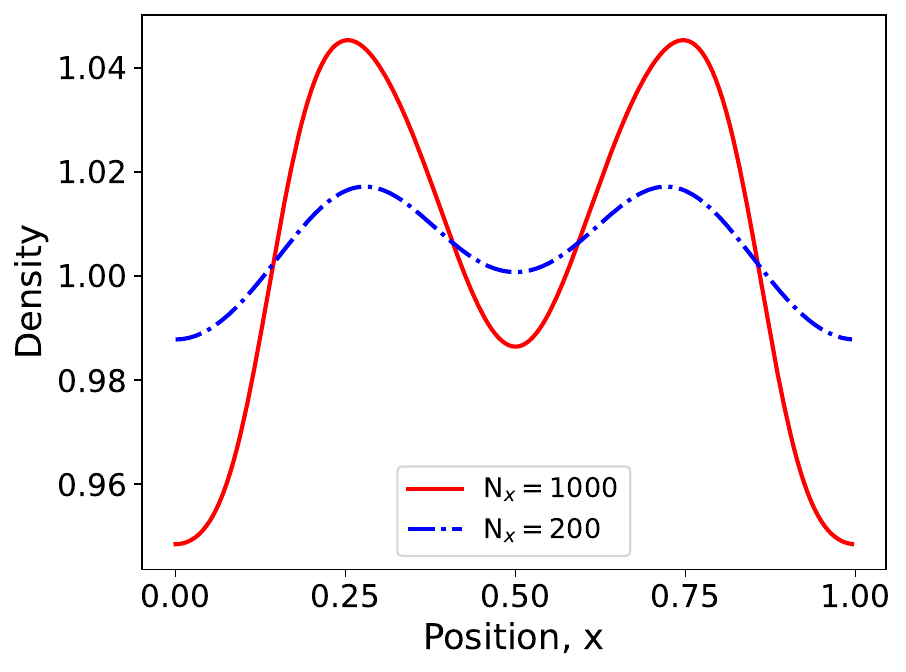}
\caption{$t=0.08$}
\label{Fig:CAW_t8_den}
\end{subfigure}
\vfill
\begin{subfigure}[b]{0.31\textwidth}
\centering
\includegraphics[width=\textwidth]{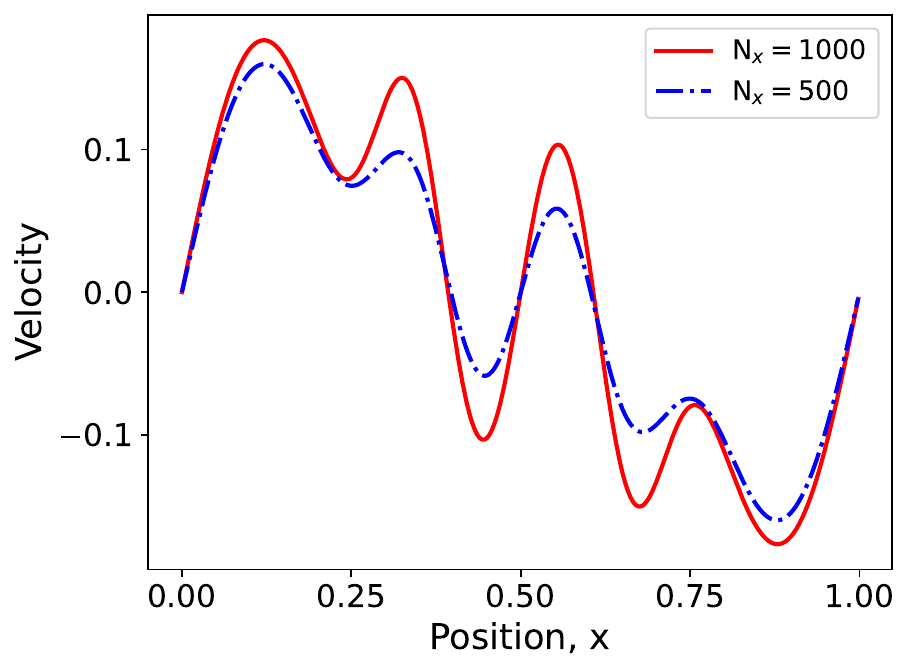}
\caption{$t=0.04$}
\label{Fig:CAW_t4_vel}
\end{subfigure}
\hspace{-0.2cm}
\begin{subfigure}[b]{0.3\textwidth}
\centering
\includegraphics[width=\textwidth]{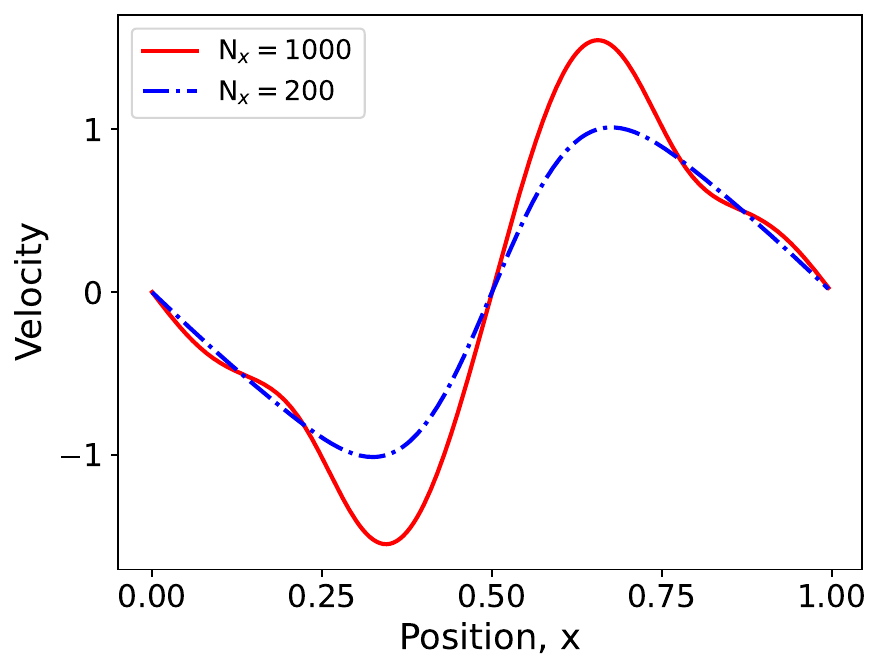}
\caption{$t=0.06$}
\label{Fig:CAW_t6_vel}
\end{subfigure}
\hspace{-0.2cm}
\begin{subfigure}[b]{0.3\textwidth}
\centering
\includegraphics[width=\textwidth]{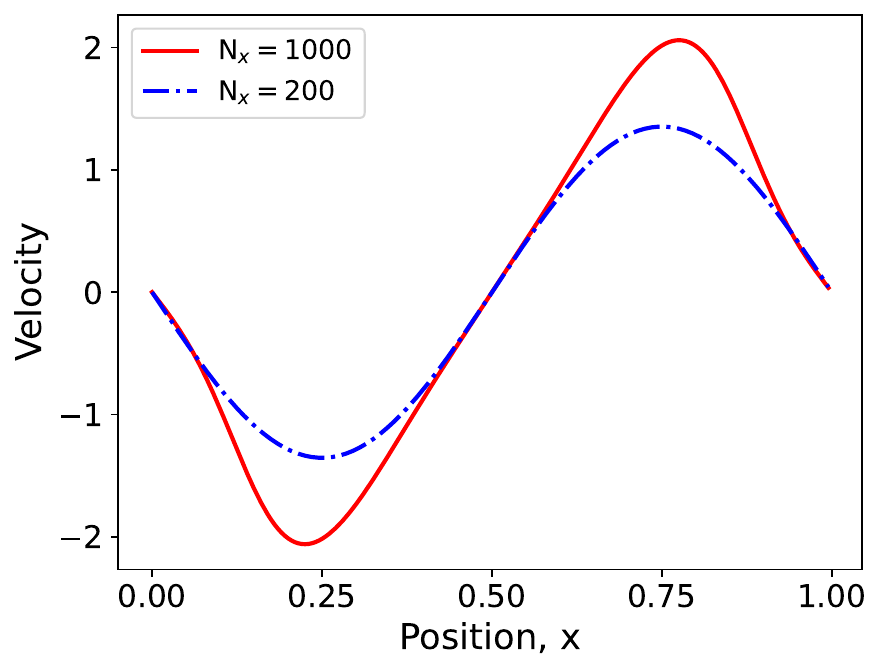}
\caption{$t=0.08$}
\label{Fig:CAW_t8_vel}
\end{subfigure}
\vfill
\begin{subfigure}[b]{0.32\textwidth}
\centering
\includegraphics[width=\textwidth]{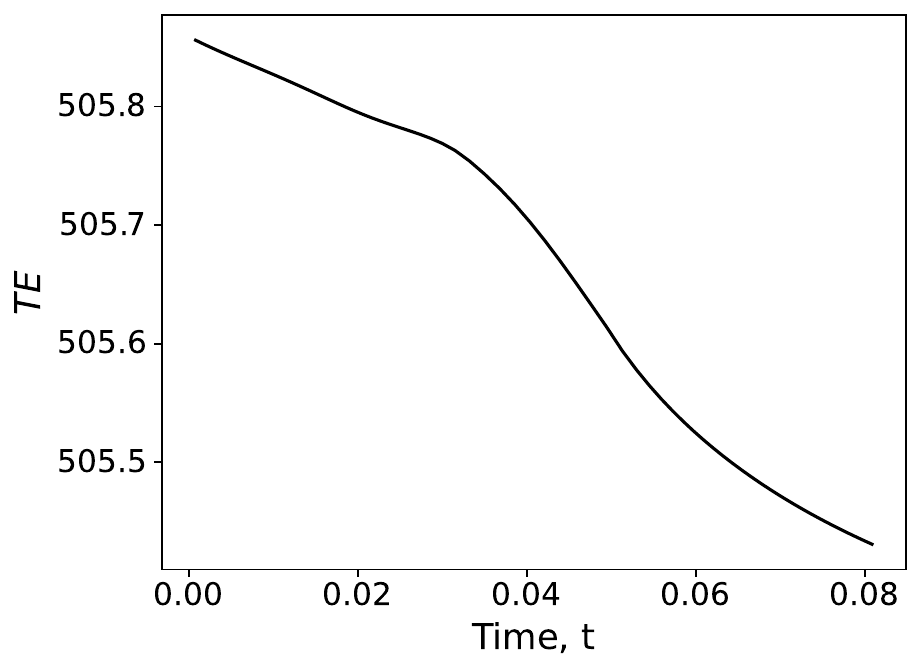}
\label{Fig:CAW_TE}
\end{subfigure}
\hspace{-0.2cm}
\begin{subfigure}[b]{0.31\textwidth}
\centering
\includegraphics[width=\textwidth]{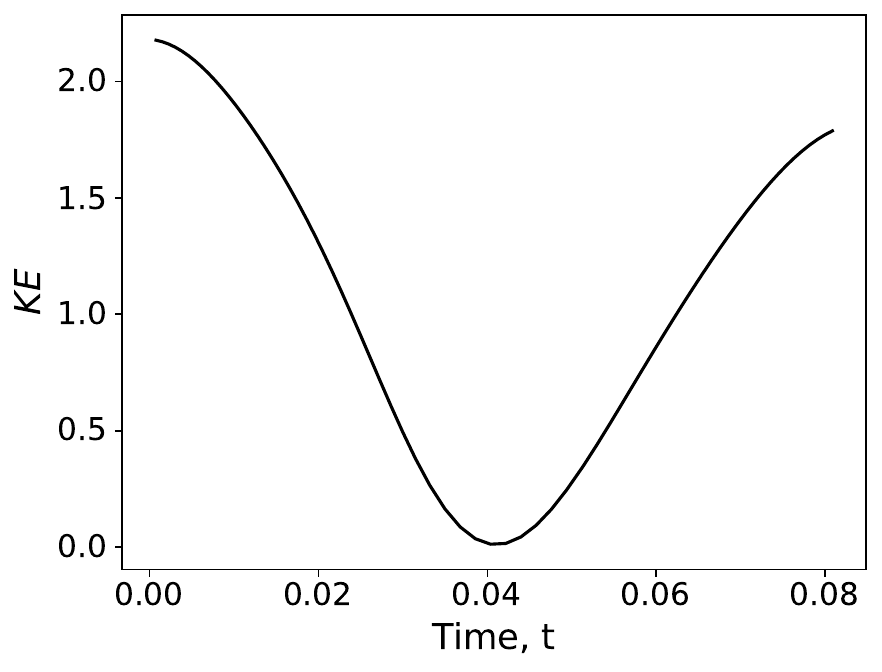}
\label{CAW_KE}
\end{subfigure}
\hspace{-0.2cm}
\begin{subfigure}[b]{0.33\textwidth}
\centering
\includegraphics[width=\textwidth]{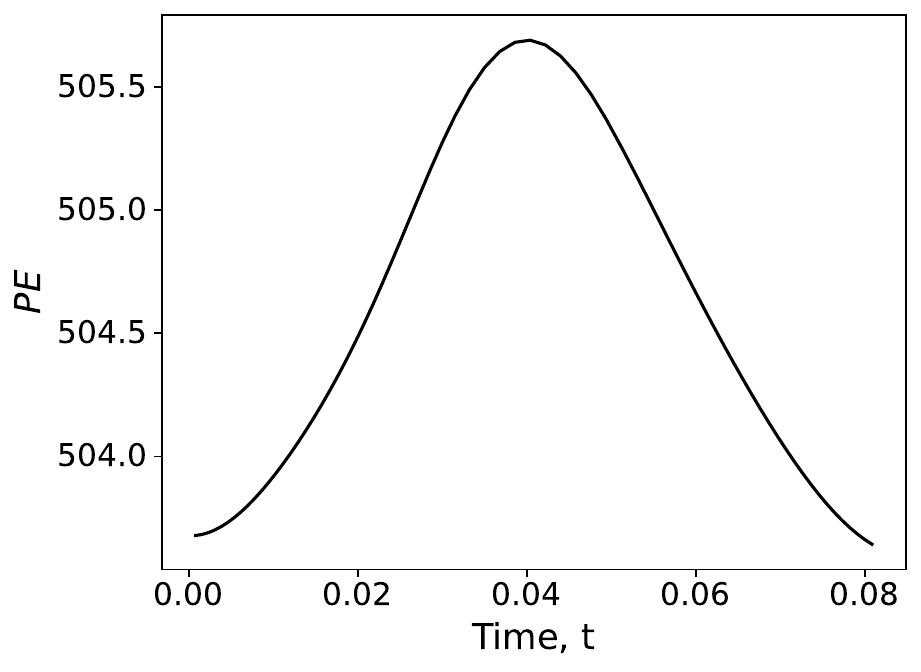}
\label{Fig:CAW_PE}
\end{subfigure}
\caption{\centering \textbf{Example \ref{Subsec:CAW} - Colliding acoustic waves problem:} Top: Density at t=0.04,0.06,0.08 (left to right); Middle - velocity at t=0.04,0.06,0.08 (left to right); Bottom - Total energy, kinetic energy, potential energy (left to right)}
\label{Fig:CAW}
\end{figure}

This problem is considered on the spatial domain, $\Omega=[-1,1]$, and the initial conditions are:
\begin{gather}
    \rho^0(x)=0.955+0.5\varepsilon \left(1-\cos \left(2\pi x\right) \right), \\
    u^0(x)=-\text{sign}(x)\sqrt{\gamma} \left(1- \cos \left( 2\pi x\right) \right).
\end{gather}
The boundary conditions are periodic, and the parametric values are: $\kappa=1$, $\gamma=1.4$, and $\varepsilon=0.1$. Further, $\lambda=1$ and $C=0.9$. A convergence study is performed at time $t=0.08$, by considering different numbers of grid points: N$_x=20,50, 100, 200, 250$ and $500$. The $L^2$ error is computed with the reference as N$_x=1000$, and the convergence rates of $\varrho,u$ are shown in \Cref{tab: CAW EOC_rho_u}. It can be seen that $\varrho$ and $u$ converge with rates of about $1.1$ and $1.4$, respectively. \Cref{Fig:CAW} shows the density, the velocity plots at time $t=0.04,0.06,0.08$, and the energy plots with respect to time for N$_x=1000$. The results indicate that the proposed fully discrete method \eqref{Update mass}-\eqref{Update mom} is energy-stable for this problem with non-well-prepared initial data. We observe that the kinetic and potential energies are balancing each other to yield decreasing total energy, which corresponds to the desired energy stability property.  

\begin{table}[h!]
    \centering
    \renewcommand{\arraystretch}{1.3} 
    \setlength{\tabcolsep}{10pt}      
    
    \begin{tabular}{|c|c|c|c|c|c|}
        \hline
        \multirow{2}{*}{\textbf{N$_x$}} & \multirow{2}{*}{$\mathbf{\Delta x}$} & 
        \multicolumn{2}{c|}{$\mathbf{\varepsilon=0.1}$} & 
        \multicolumn{2}{c|}{$\mathbf{\varepsilon=0.1}$} \\ 
        \cline{3-6}
        & & $\mathbf{||\varrho \textbf{ error}||_{L^2}}$ & \textbf{EOC} 
        & $\mathbf{||u \textbf{ error}||_{L^2}}$ & \textbf{EOC} \\ 
        \hline
        20  & 0.1  & 0.04557  & -      & 1.69840  & -  \\  
        50  & 0.04  & 0.04321  & 0.0579 & 1.37651  & 0.2293 \\  
        100 & 0.02  & 0.04762  & -0.1399 & 1.04690  & 0.3949  \\  
        200 & 0.01 & 0.03184  & 0.5804 & 0.55072  & 0.9267   \\  
        250 & 0.008 & 0.02506  & 1.0737 & 0.44314  & 0.9740   \\  
        500 & 0.004 & 0.01187  & 1.0776 & 0.17106  & 1.3733  \\  
        \hline
    \end{tabular}
    
    \caption{\centering \textbf{Example \ref{Subsec:CAW} - Colliding acoustic waves problem:} Convergence rates of $L^2$ error in $\varrho$ and $u$.}
    \label{tab: CAW EOC_rho_u}
\end{table}

\subsection{Riemann problem}
\label{Subsec:RP}
\begin{figure}[h!]
\centering
\begin{subfigure}[b]{0.3\textwidth}
\centering
\includegraphics[width=\textwidth]{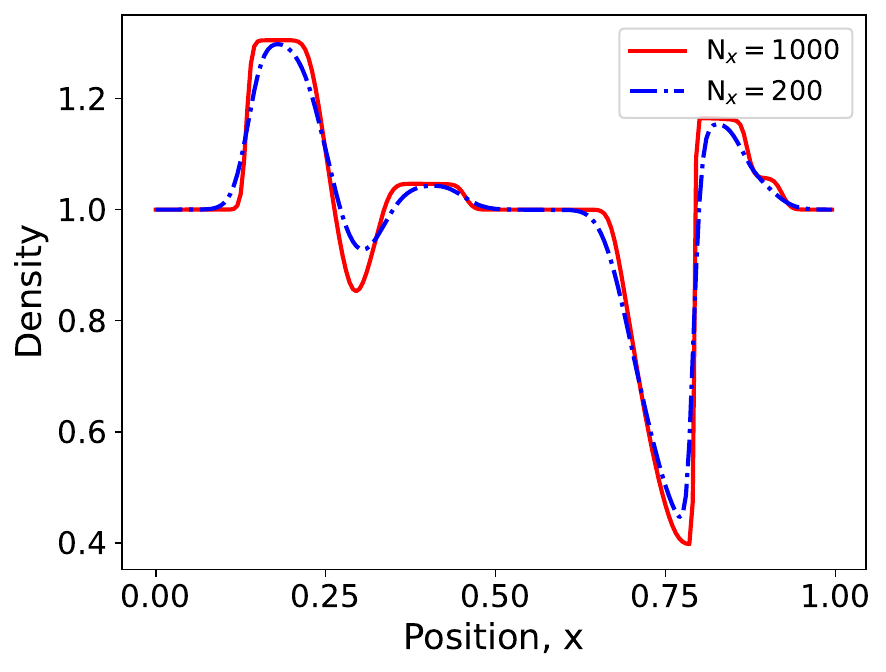}
\label{Fig:RP_eps8_den}
\end{subfigure}
\hspace{-0.2cm}
\begin{subfigure}[b]{0.31\textwidth}
\centering
\includegraphics[width=\textwidth]{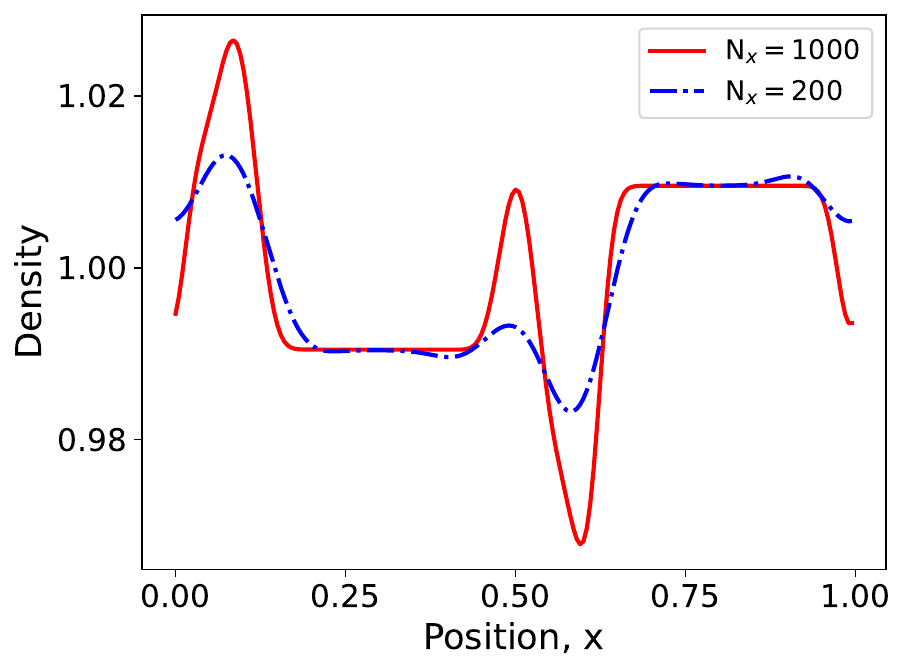}
\label{Fig:RP_eps3_den}
\end{subfigure}
\hspace{-0.2cm} 
\begin{subfigure}[b]{0.34\textwidth}
\centering
\includegraphics[width=\textwidth]{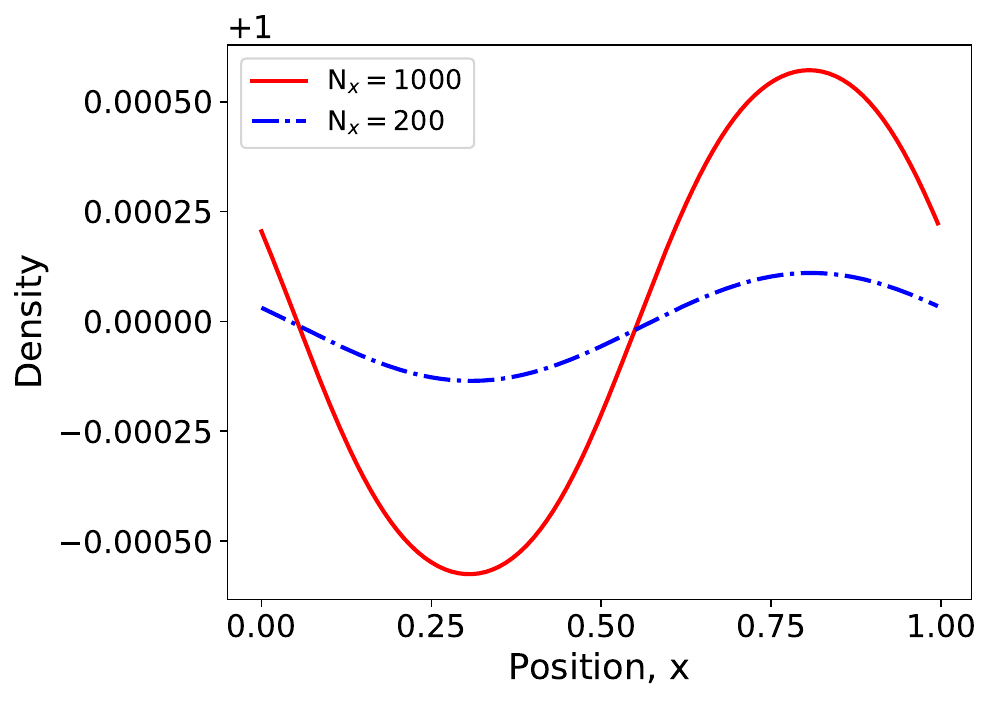}
\label{Fig:RP_eps5_den}
\end{subfigure}
\vfill
\begin{subfigure}[b]{0.3\textwidth}
\centering
\includegraphics[width=\textwidth]{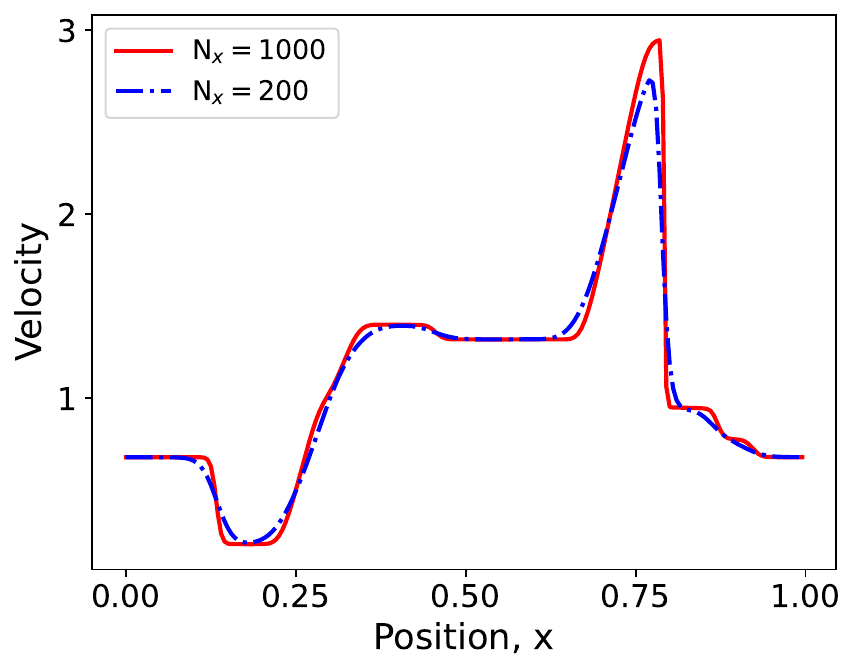}
\label{Fig:RP_eps8_vel}
\end{subfigure}
\hspace{-0.2cm}
\begin{subfigure}[b]{0.31\textwidth}
\centering
\includegraphics[width=\textwidth]{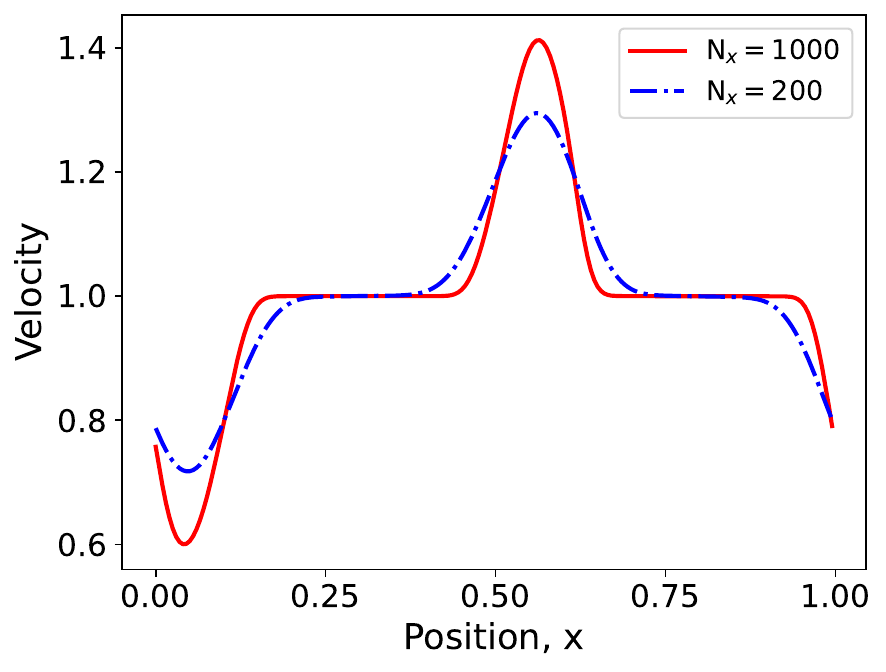}
\label{Fig:RP_eps3_vel}
\end{subfigure}
\hspace{-0.2cm}
\begin{subfigure}[b]{0.32\textwidth}
\centering
\includegraphics[width=\textwidth]{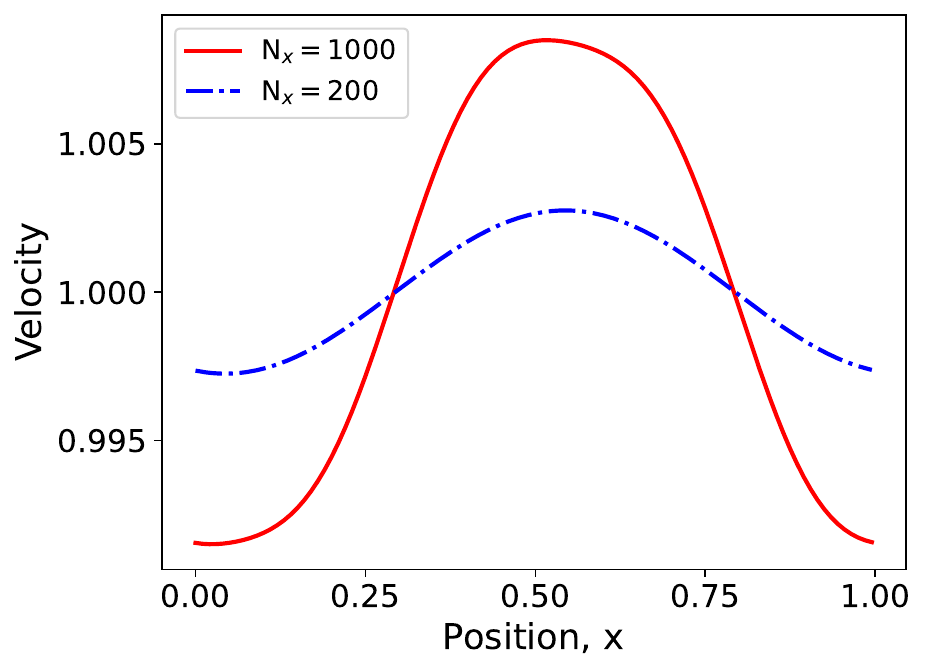}
\label{Fig:RP_eps5_vel}
\end{subfigure}
\vfill
\begin{subfigure}[b]{0.32\textwidth}
\centering
\includegraphics[width=\textwidth]{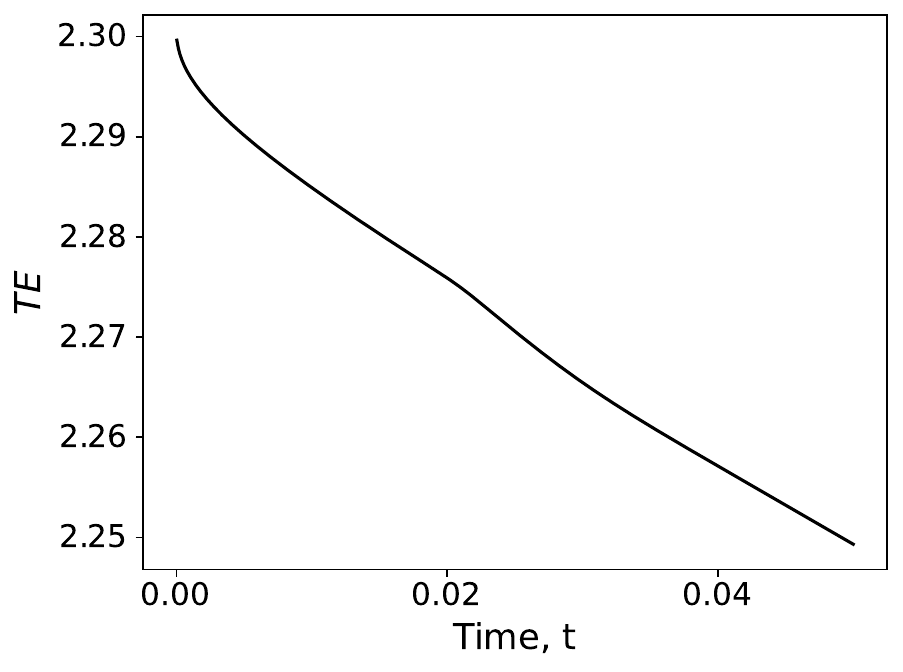}
\label{Fig:RP_TE_eps8}
\end{subfigure}
\hspace{-0.2cm}
\begin{subfigure}[b]{0.33\textwidth}
\centering
\includegraphics[width=\textwidth]{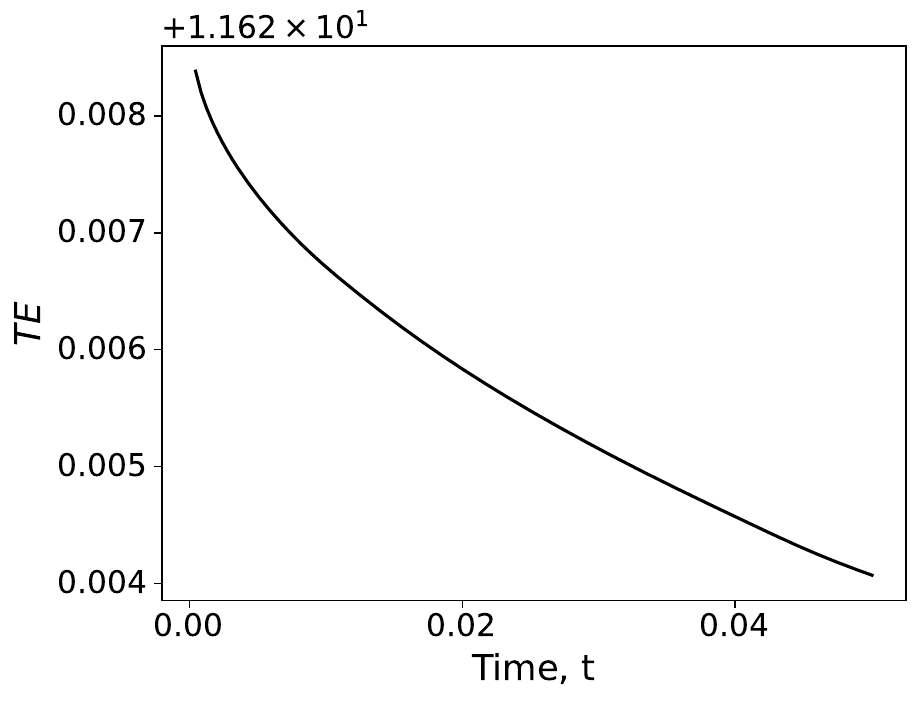}
\label{Fig:RP_TE_eps3}
\end{subfigure}
\hspace{-0.2cm}
\begin{subfigure}[b]{0.33\textwidth}
\centering
\includegraphics[width=\textwidth]{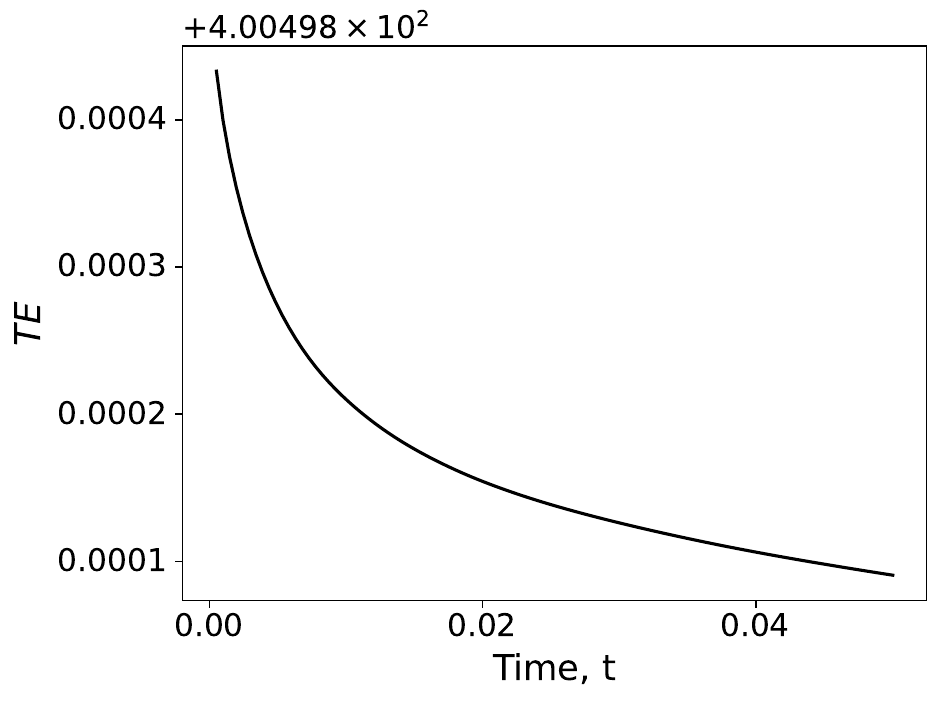}
\label{Fig:RP_TE_eps5}
\end{subfigure}

\begin{subfigure}[b]{0.32\textwidth}
\centering
\includegraphics[width=\textwidth]{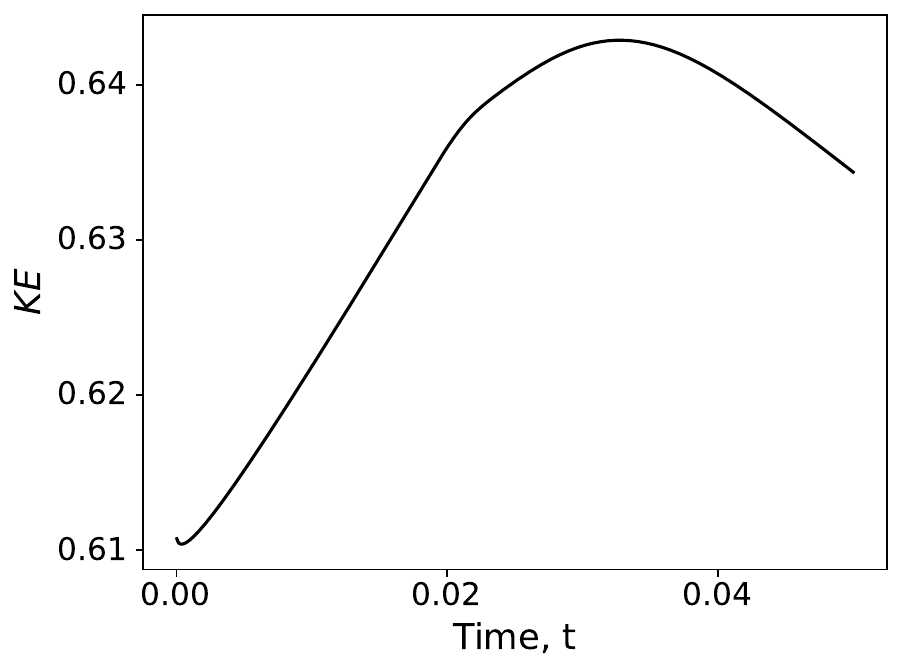}
\end{subfigure}
\hspace{-0.2cm}
\begin{subfigure}[b]{0.33\textwidth}
\centering
\includegraphics[width=\textwidth]{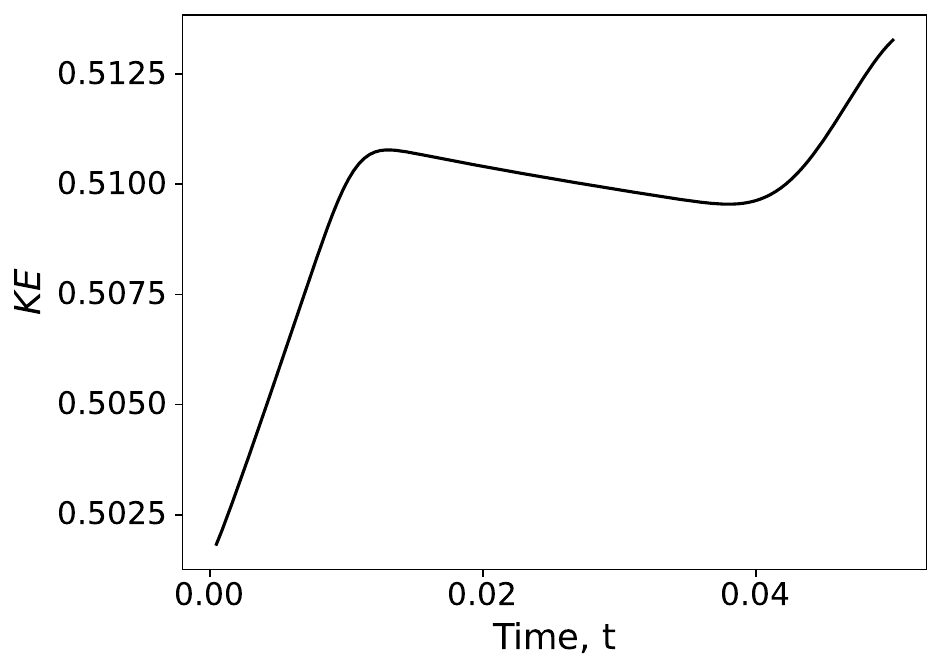}
\end{subfigure}
\hspace{-0.2cm}
\begin{subfigure}[b]{0.34\textwidth}
\centering
\includegraphics[width=\textwidth]{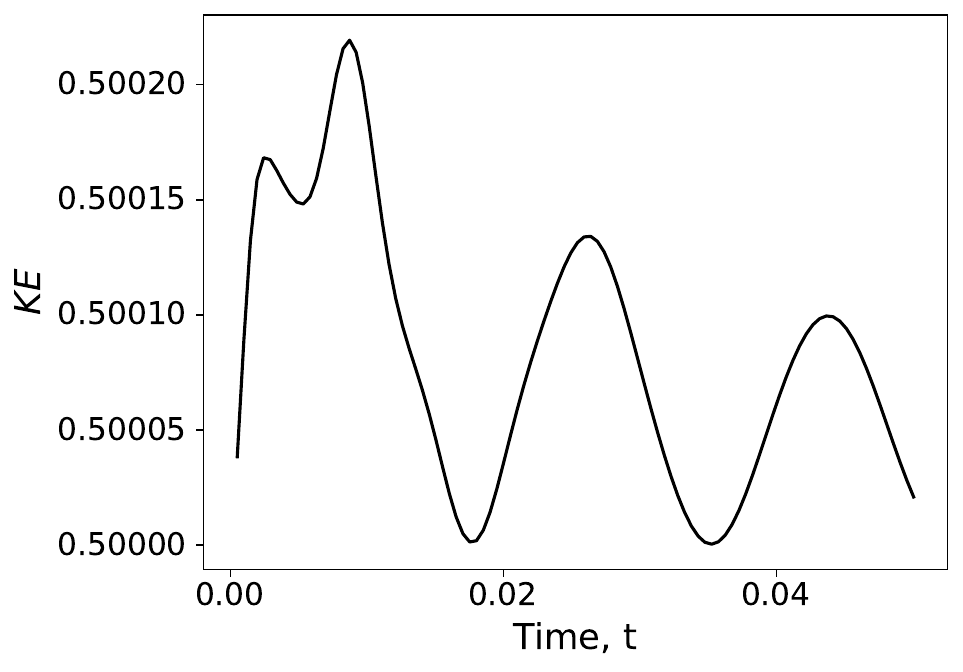}
\end{subfigure}

\begin{subfigure}[b]{0.32\textwidth}
\centering
\includegraphics[width=\textwidth]{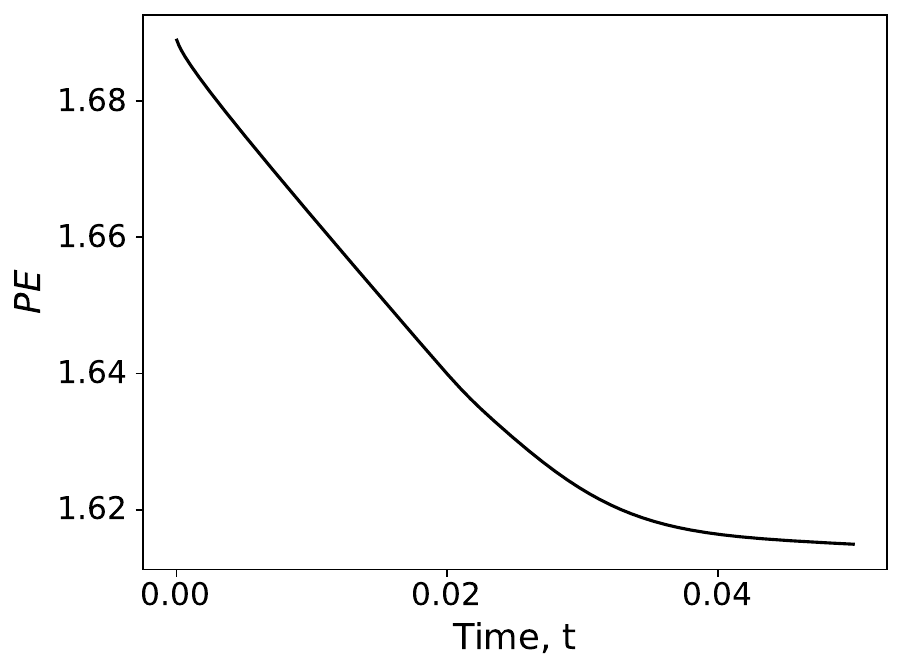}
\caption{$\varepsilon=0.8$}
\end{subfigure}
\hspace{-0.2cm}
\begin{subfigure}[b]{0.33\textwidth}
\centering
\includegraphics[width=\textwidth]{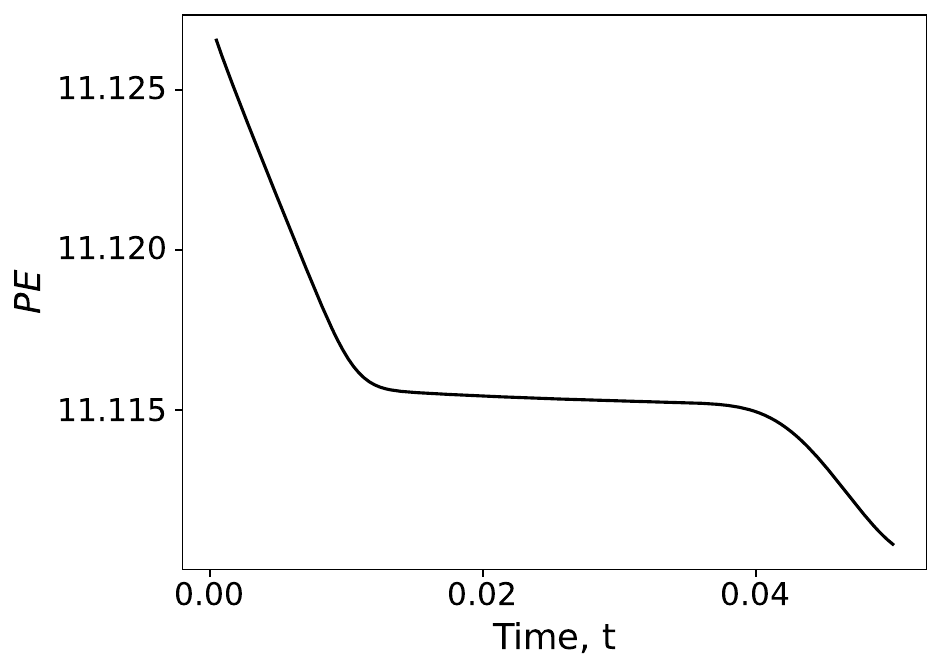}
\caption{$\varepsilon=0.3$}
\end{subfigure}
\hspace{-0.2cm}
\begin{subfigure}[b]{0.33\textwidth}
\centering
\includegraphics[width=\textwidth]{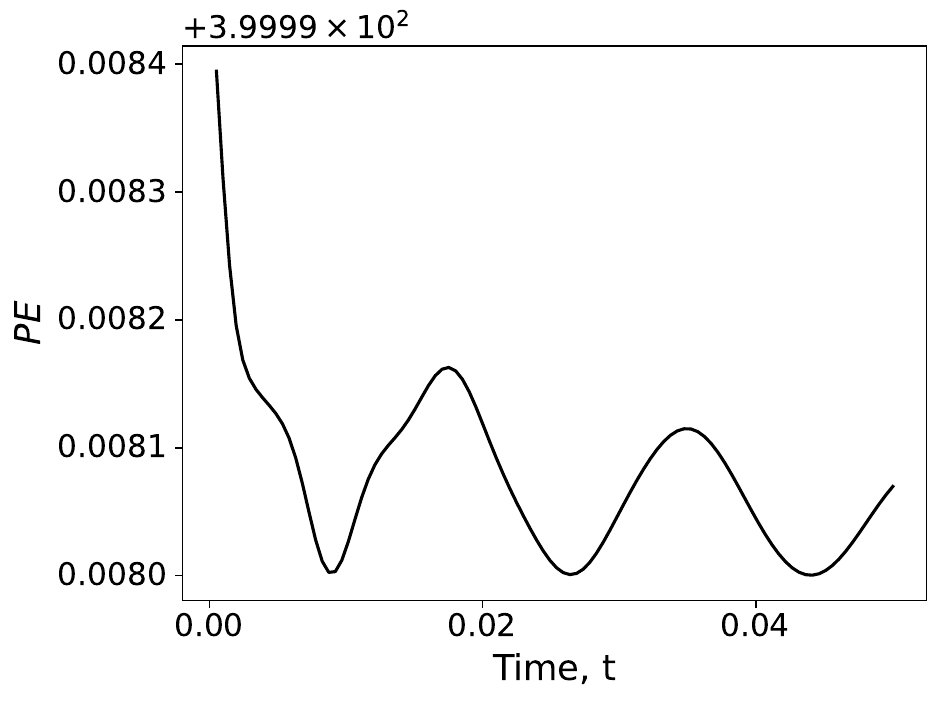}
\caption{$\varepsilon=0.05$}
\end{subfigure}

\caption{\centering \textbf{Example \ref{Subsec:RP} - Riemann problem:} Top to Bottom: Density, velocity, total energy, kinetic energy, potential energy}
\label{Fig:RP}
\end{figure}

This Riemann problem was proposed in \cite{ap_degond1}. The computational domain is taken to be $\Omega = [0,1]$, and the initial conditions are:
\begin{eqnarray}
    \varrho^0(x)=1, & \quad \left(\varrho u\right)^0(x)=1-\varepsilon^2/2, & \quad \text{if } x \in [0,0.2] \cup [0.8,1] \\
    \varrho^0(x)=1+\varepsilon^2, & \quad \left(\varrho u\right)^0(x)=1, & \quad \text{if } x \in (0.2,0.3] \\
    \varrho^0(x)=1, & \quad \left(\varrho u\right)^0(x)=1+\varepsilon^2/2, &  \quad \text{if } x \in (0.3,0.7] \\
    \varrho^0(x)=1-\varepsilon^2, & \quad \left(\varrho u\right)^0(x)=1, &  \quad \text{if } x \in (0.7,0.8]. 
\end{eqnarray}
We consider periodic boundary conditions, and choose the following parametric values: $\kappa=1$ and $\gamma=2$. Three different values of $\varepsilon$, $0.8,0.3$, and $0.05$ are considered. $\lambda=1$ is used as the numerical diffusion coefficient, and a CFL number of $C=0.5$ is used for $\varepsilon=0.3$ and $\varepsilon = 0.05$, while $C=0.1$ is used for $\varepsilon=0.8$. \Cref{Fig:RP} shows the density, the velocity plots at time $t=0.05$, and the energy plots with respect to time for N$_x=1000$. The results confirm that the proposed method \eqref{Update mass}-\eqref{Update mom} is energy-stable for different values of $\varepsilon$. We observe that the kinetic and potential energies are balancing each other to yield a decrease in total energy for different $\varepsilon$.   

\subsection{Gresho vortex problem}
\label{Subsec:GV}
This is a rotating vortex problem proposed in \cite{Num_Gresho1,Num_Gresho2}. Here, a vortex of radius $R=0.4$ centered at $(x_1,x_2)=(0.5,0.5)$ is considered at initial time $t=0$. The initial background state is: $\varrho_0=1$, $\mathbf{u}_0=(u_{1_0},0)^T$ with $u_{1_0}=0.1$ and $p_0=1$. The velocity of the vortex is given by:
\begin{equation*}
    u_{\theta}(r) = \left\{ \begin{matrix}
        2\frac{r}{R} & \text{if } 0\leq r < \frac{R}{2} \\
        2 \left( 1- \frac{r}{R} \right) & \text{if } \frac{R}{2} \leq r < R \\
        0 & \text{if } r\geq R
    \end{matrix} \right.,
\end{equation*}
and its components in Cartesian coordinates are:
\begin{equation*}
    u_1(x_1,x_2)=u_{1_0} - \frac{x_2 - x_{2_0}}{r} u_{\theta}(r), \ u_2(x_1,x_2)= \frac{x_1 - x_{1_0}}{r} u_{\theta}(r).
\end{equation*}
Here, $r=\sqrt{(x_1-0.5)^2+(x_2-0.5)^2}$. The expression for pressure is obtained by balancing the pressure gradient and the centrifugal force $\left( \textit{i.e., } \varrho_0\frac{u_{\theta}^2}{r}= \frac{1}{\varepsilon^2}\frac{\partial p}{\partial r} \right)$:
\begin{equation*}
    p(r)=p_0 + \varepsilon^2 \left\{ \begin{matrix}
        2\frac{r^2}{R^2} +2-\log 16 & \text{if } 0\leq r < \frac{R}{2} \\
        2\frac{r^2}{R^2} - 8 \frac{r}{R} +4 \log \left( \frac{r}{R} \right)+ 6 & \text{if } \frac{R}{2} \leq r < R \\
        0 & \text{if } r\geq R
    \end{matrix} \right..
\end{equation*}
 Considering the asymptotic expansion $p=p_0 + \varepsilon^2 p_2$, $\varrho=\varrho_0 + \varepsilon^2 \varrho_2$ and comparing $\frac{\partial p}{\partial \varepsilon}=2\varepsilon p_2$ with $\frac{\partial p}{\partial \varepsilon}= \frac{\partial p}{\partial \varrho} \frac{\partial \varrho}{\varepsilon} = \frac{\gamma p }{\varrho} 2\varepsilon \varrho_2$, we obtain $\varrho_2 = \frac{p_2}{\gamma}$. Note that $p=1$ and $\varrho=1$ up to the leading order. Thus, using $\varrho_2$, we can evaluate $\varrho$ as $\varrho=\varrho_0 + \varepsilon^2 \varrho_2$. Further, we consider periodic boundary conditions, and take the parametric values $\kappa$ and $\gamma$ to be $1$ and $1.4$, respectively. We consider three different values of $\varepsilon$: $\varepsilon=0.1,0.01$, and $0.001$, and use a CFL number of $C=0.5$ for $\varepsilon=0.1,0.01$, and $C=0.1$ for $\varepsilon=0.001$.
 \par
\begin{figure}[h!]
\centering
\begin{subfigure}[b]{0.3\textwidth}
\centering
\includegraphics[width=\textwidth]{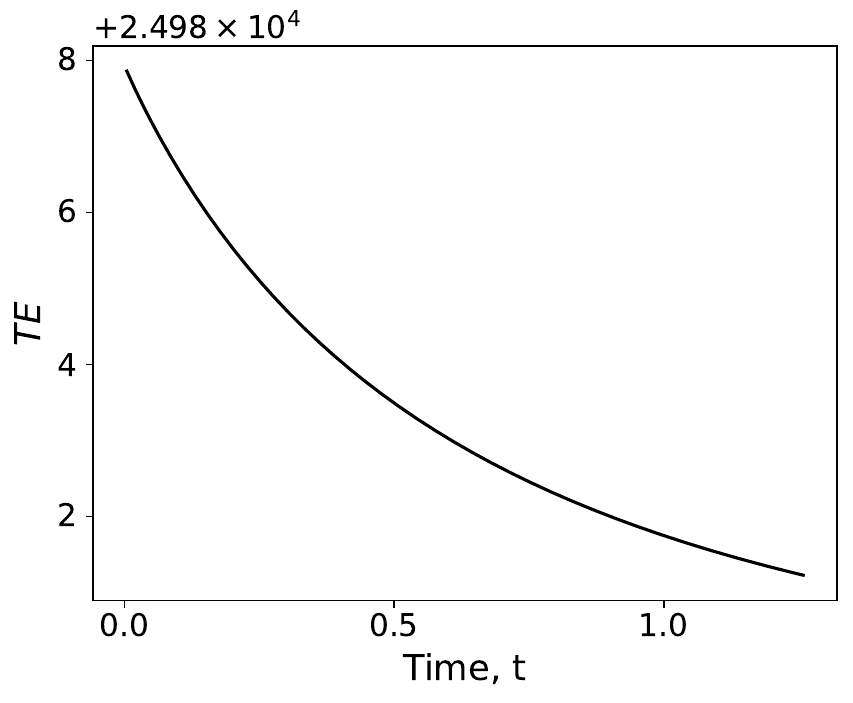}
\label{Fig:GV_TE_eps01}
\end{subfigure}
\hspace{-0.2cm}
\begin{subfigure}[b]{0.3\textwidth}
\centering
\includegraphics[width=\textwidth]{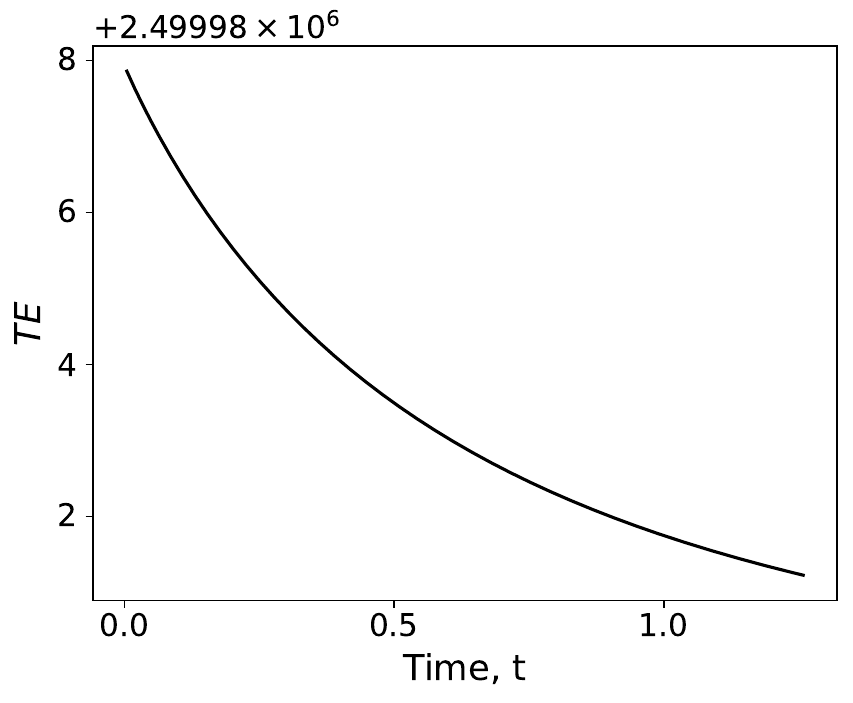}
\label{Fig:GV_TE_eps001}
\end{subfigure}
\hspace{-0.2cm}
\begin{subfigure}[b]{0.32\textwidth}
\centering
\includegraphics[width=\textwidth]{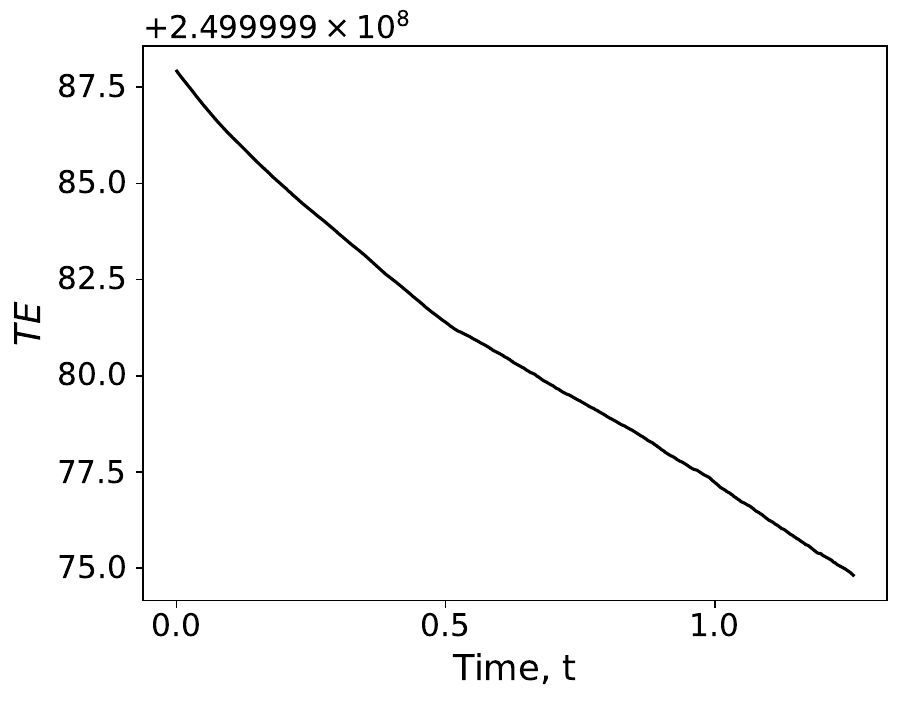}
\label{Fig:GV_TE_eps0001}
\end{subfigure}

\begin{subfigure}[b]{0.3\textwidth}
\centering
\includegraphics[width=\textwidth]{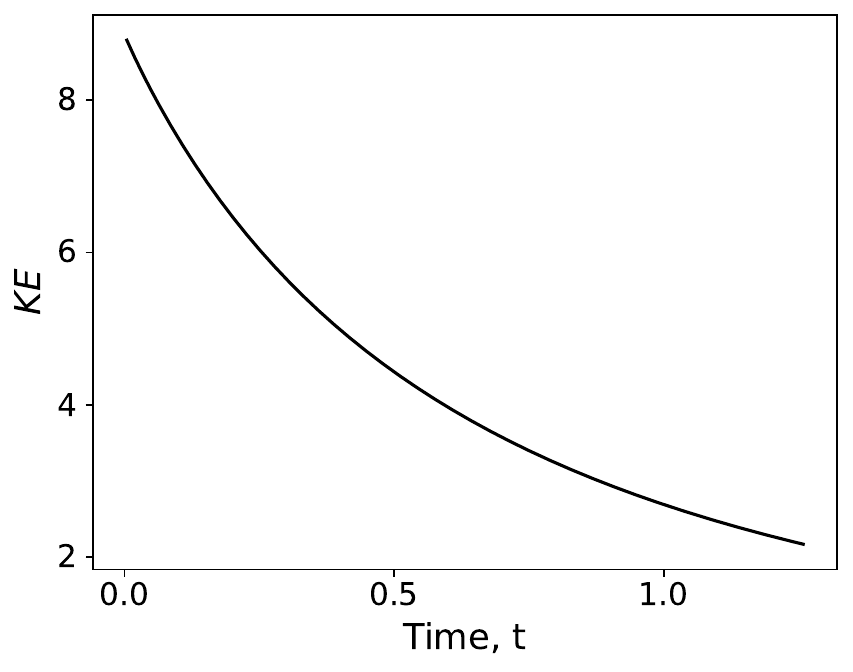}
\end{subfigure}
\hspace{-0.2cm}
\begin{subfigure}[b]{0.3\textwidth}
\centering
\includegraphics[width=\textwidth]{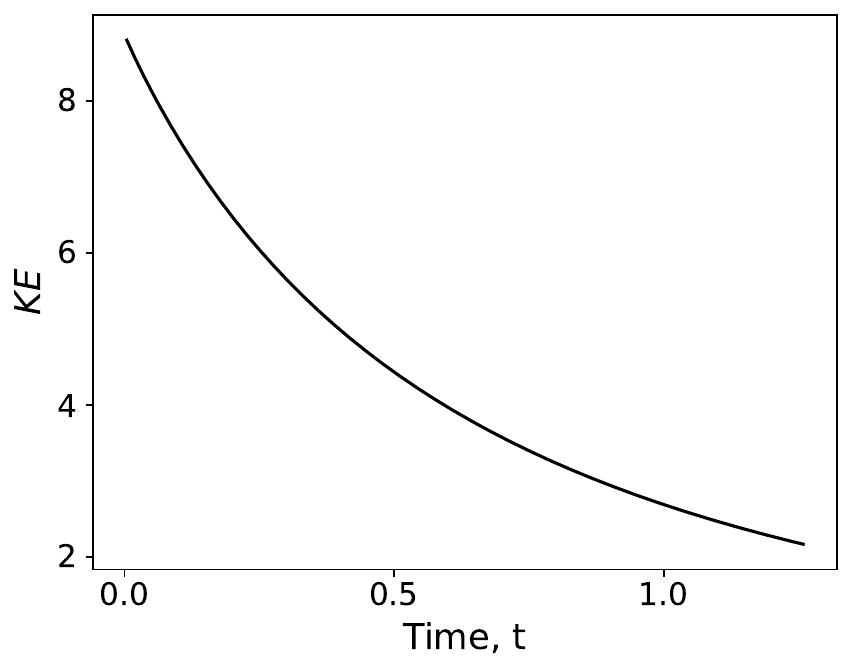}
\end{subfigure}
\hspace{-0.2cm}
\begin{subfigure}[b]{0.3\textwidth}
\centering
\includegraphics[width=\textwidth]{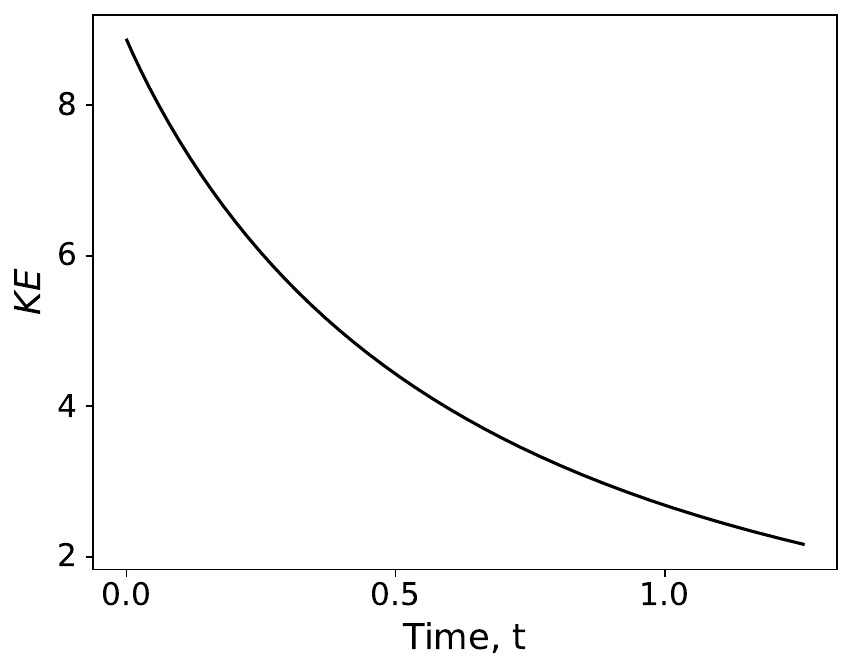}
\end{subfigure}

\begin{subfigure}[b]{0.31\textwidth}
\centering
\includegraphics[width=\textwidth]{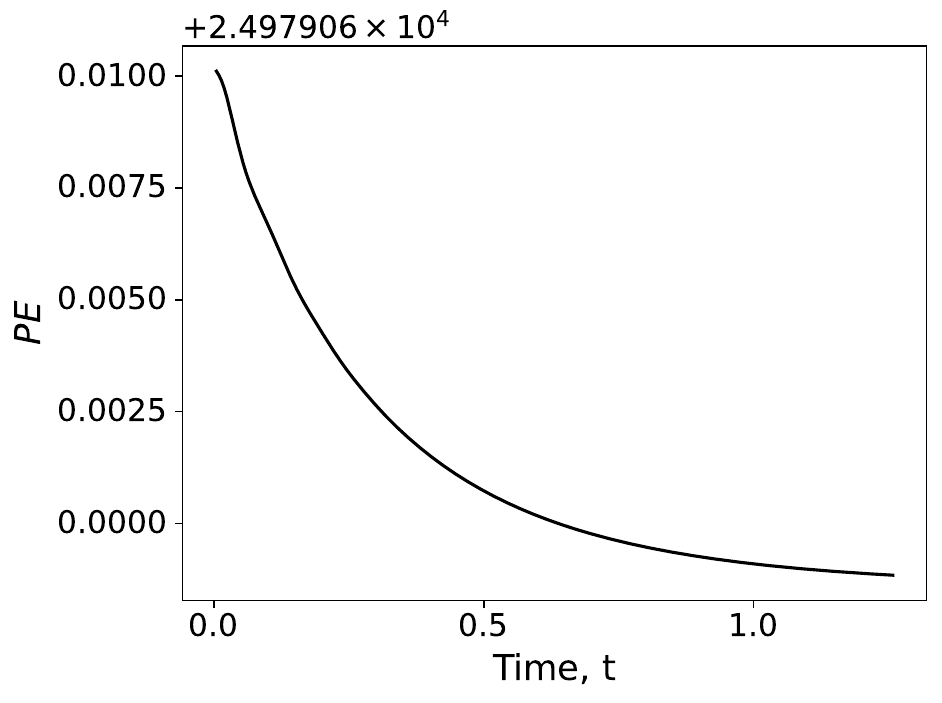}
\caption{$\varepsilon=0.1$}
\end{subfigure}
\hspace{-0.2cm}
\begin{subfigure}[b]{0.31\textwidth}
\centering
\includegraphics[width=\textwidth]{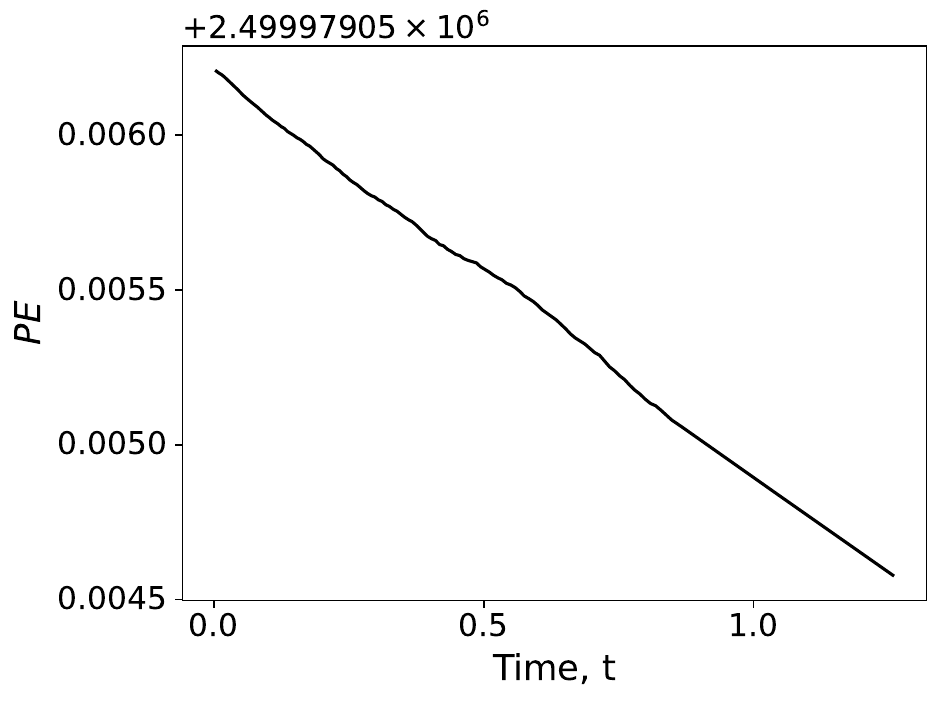}
\caption{$\varepsilon=0.01$}
\end{subfigure}
\hspace{-0.2cm}
\begin{subfigure}[b]{0.28\textwidth}
\centering
\includegraphics[width=\textwidth]{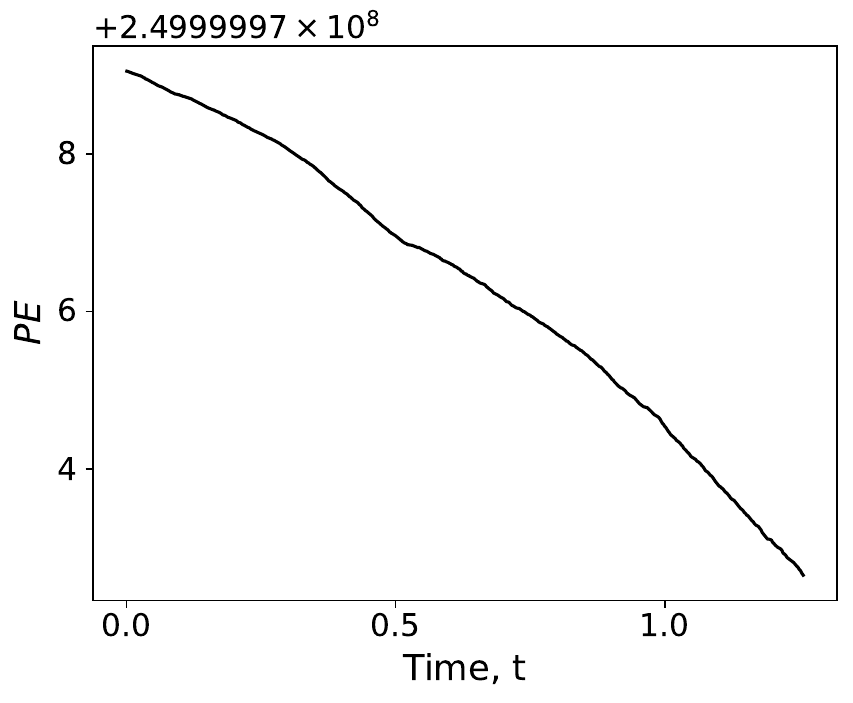}
\caption{$\varepsilon=0.001$}
\end{subfigure}

\caption{\centering \textbf{Example \ref{Subsec:GV} - Gresho vortex problem:} Top - Total energy; Middle - Kinetic energy; Bottom - Potential energy ($\lambda=1$)}
\label{Fig:GV_energy}
\end{figure}

\begin{figure}[h!]
\centering
\begin{subfigure}[b]{0.31\textwidth}
\centering
\includegraphics[width=\textwidth]{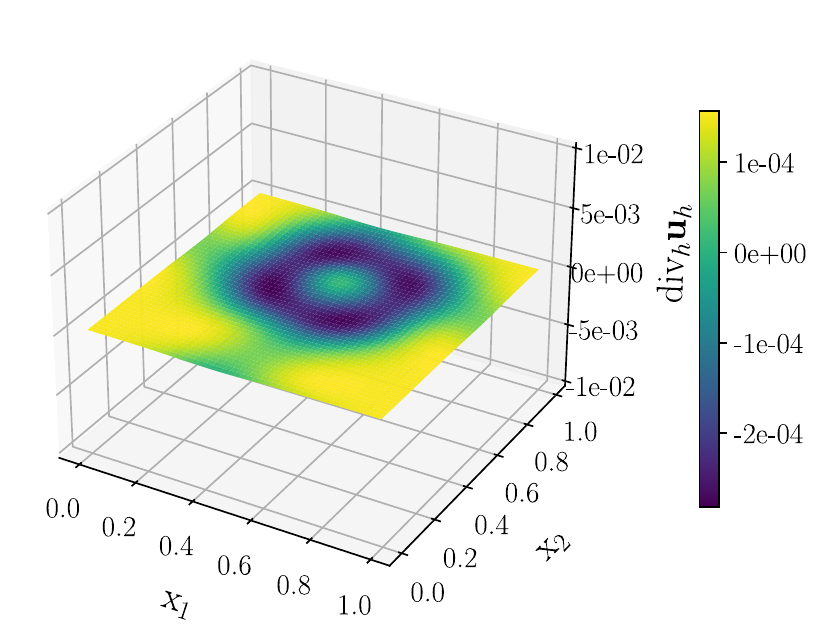}
\caption{$\varepsilon=0.1$}
\end{subfigure}
\hspace{-0.2cm}
\begin{subfigure}[b]{0.31\textwidth}
\centering
\includegraphics[width=\textwidth]{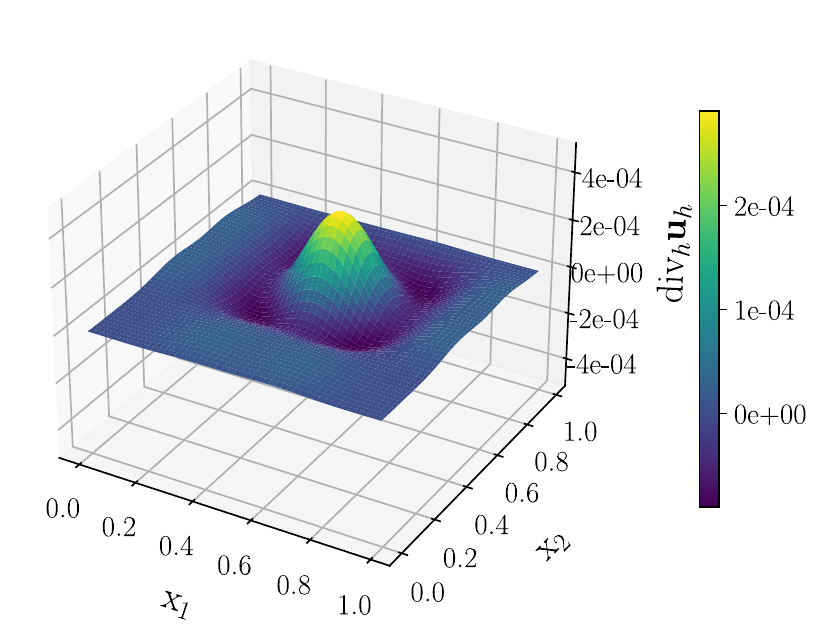}
\caption{$\varepsilon=0.01$}
\end{subfigure}
\hspace{-0.2cm}
\begin{subfigure}[b]{0.31\textwidth}
\centering
\includegraphics[width=\textwidth]{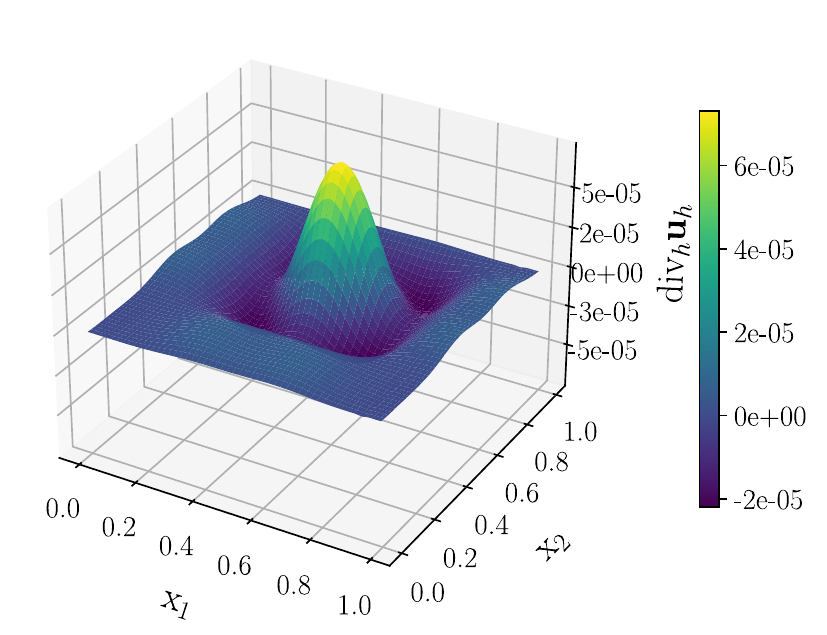}
\caption{$\varepsilon=0.001$}
\end{subfigure}

\caption{\centering \textbf{Example \ref{Subsec:GV} - Gresho vortex problem:} Surface plots of $\divh{} \mathbf{u}_h$ ($\lambda=1$)}
\label{Fig:GV_divu}
\end{figure}

\begin{figure}[h!]
\centering
\begin{subfigure}[b]{0.32\textwidth}
\centering
\includegraphics[width=\textwidth]{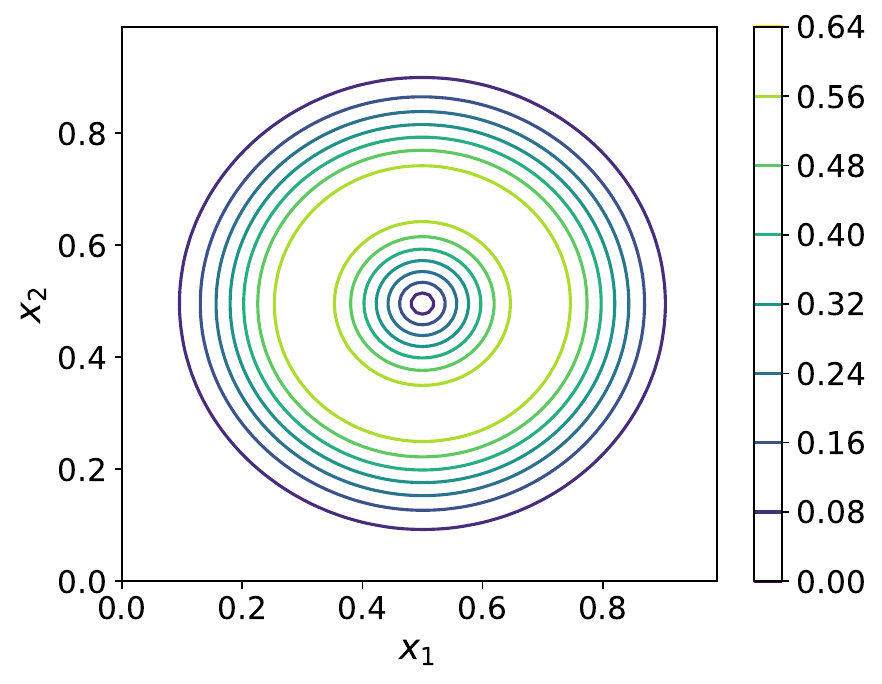}
\caption{$\varepsilon=0.1$}
\label{Fig:GV_contour_eps01}
\end{subfigure}
\hspace{-0.2cm}
\begin{subfigure}[b]{0.32\textwidth}
\centering
\includegraphics[width=\textwidth]{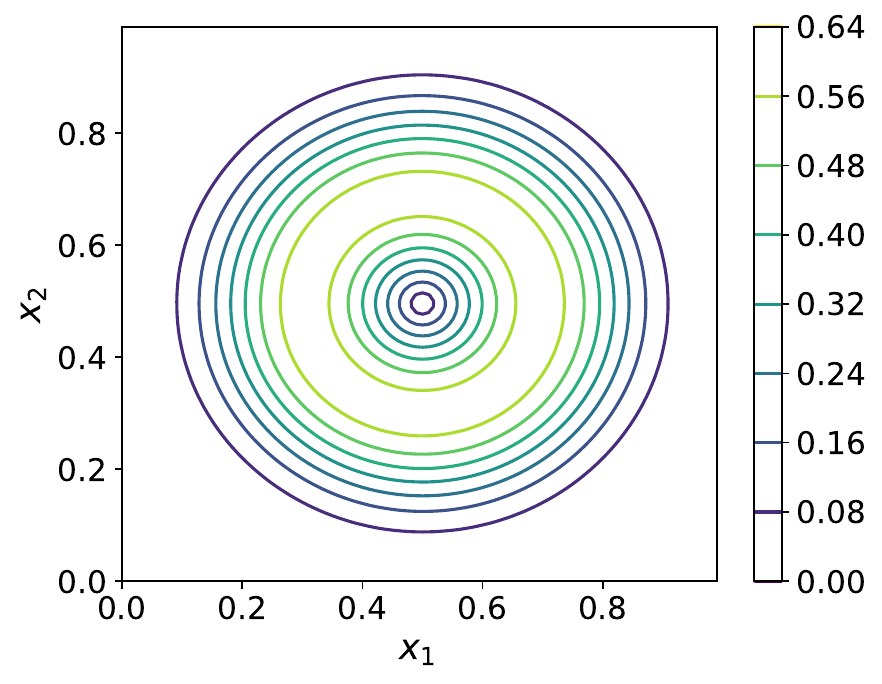}
\caption{$\varepsilon=0.01$}
\label{Fig:GV_contour_eps001}
\end{subfigure}
\hspace{-0.2cm}
\begin{subfigure}[b]{0.32\textwidth}
\centering
\includegraphics[width=\textwidth]{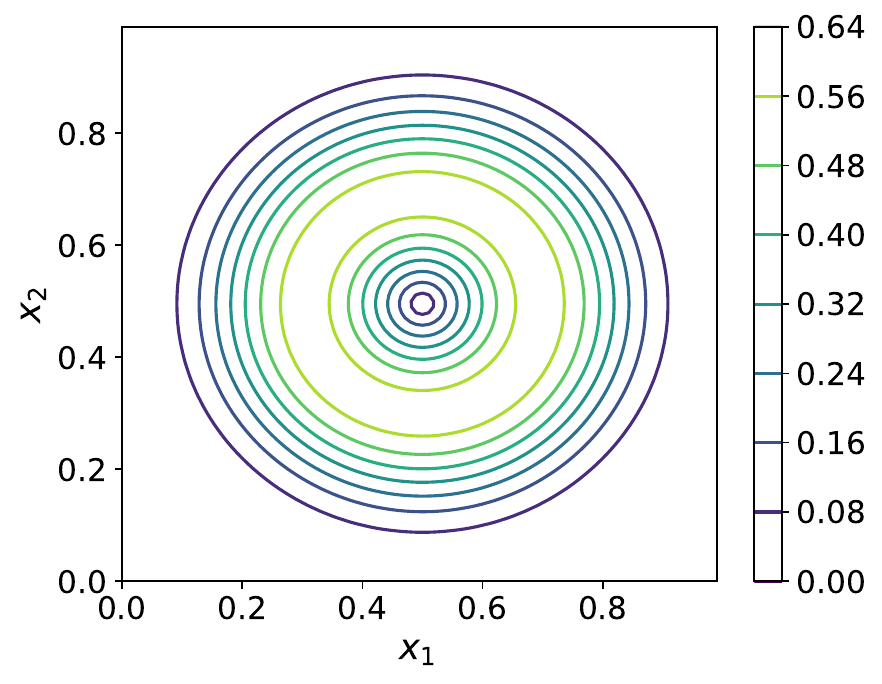}
\caption{$\varepsilon=0.001$}
\label{Fig:GV_contour_eps0001}
\end{subfigure}
\caption{\centering \textbf{Example \ref{Subsec:GV} - Gresho vortex problem:} Contour plots of Mach ratio ($\lambda$ based on \eqref{lambda cond numerical}, with $c=30,200,200$ respectively for $\varepsilon=0.1,0.01,0.001$)}
\label{Fig:GV_contour}
\end{figure}

 With $\lambda=1$, we observe energy stability for different values of $\varepsilon$, as shown in \Cref{Fig:GV_energy} for N$_x=100$. In addition to the expected decrease in total energy, we also observe that the kinetic and potential energies decrease in time. Furthermore, the $L^1$ norms of $\divh{} \mathbf{u}_h$ are presented in \Cref{tab: GV div_u}, and in \Cref{Fig:GV_divu}, we plot $\divh{} \mathbf{u}_h$ with respect to the spatial domain for a fixed grid of N$_x = 100$, and for different values of $\varepsilon$. These results confirm that the divergence of velocity is $\mathcal{O}(\varepsilon^2)$ as expected. \\   
In this example, it was suitable to compute $\lambda$ at every time step based on \eqref{lambda cond numerical} with $c=100, 200, 200$, respectively, for $\varepsilon=0.1,0.01,0.001$. For a fixed grid with N$_x=50$, this results in $\lambda$ lying in the intervals $[0.0812,1.0087]$, $[0.1423,0.3583]$, $[0.1411,0.3583]$, respectively, for $\varepsilon=0.1,0.01,0.001$. Thus, the different choices of $c$ for different $\varepsilon$ result in reasonable values of $\lambda$. A convergence study is performed at a time $t=R\pi$, by considering different numbers of grid points: N$_x=$N$_{x_1}=$N$_{x_2}=10, 20, 25$, and $50$. The $L^2$ error is computed with the reference as N$_x=100$, and the convergence rates of $\varrho$ and $\mathbf{u}$ respectively are shown in \Cref{tab: GV EOC_rho,tab: GV EOC_u}. It can be seen that the density $\varrho$ converges with rates of about $0.47,1.46,1.07$, respectively, for $\varepsilon=0.1,0.01,0.001$ and the velocity $\mathbf{u}$ converges with rates of about $0.79,1.74,1.19$, respectively, for $\varepsilon=0.1,0.01,0.001$. Thus, we obtain the expected convergence rates for $\varepsilon=0.01,0.001$ with this choice of $\lambda$. The convergence rates for $\varepsilon=0.1$ remains sub-optimal. We hypothesize that this sub-optimal convergence for $\varepsilon=0.1$ occurs due to $O\paraL{\varepsilon^2}$ error in the momentum equation due to linearization of pressure, since the $L^2$ errors in density (in \Cref{tab: GV EOC_rho}) are of $O\paraL{10^{-4}} \approx O\paraL{\varepsilon^{4}}$ for $\varepsilon=0.1$. These are orders of magnitude smaller than the $O\paraL{\varepsilon^2}$ error due to the linearization. This may indicate that for this problem, a larger modelling error is obtained due to the pressure linearization and $\varepsilon=0.1$ is on the borders of the weakly compressible regime that can be accessed using the linearization that we proposed. \\
Further, \Cref{Fig:GV_contour} shows the contour plots of Mach number ratio given by $\sqrt{\frac{(u_1 - u_{1_0})^2 + u_2^2}{\gamma p/ \varrho}}$ for N$_x=100$. 

\begin{table}[h!]
    \centering
    \renewcommand{\arraystretch}{1.3} 
    \setlength{\tabcolsep}{10pt}      
    
    \begin{tabular}{|c|c|c|c|c|c|}
        \hline
        \multirow{2}{*}{\textbf{N$_x$}} & \multirow{2}{*}{$\mathbf{\Delta x}$} & 
        \multicolumn{1}{c|}{$\mathbf{\varepsilon=0.1}$} & 
        \multicolumn{1}{c|}{$\mathbf{\varepsilon=0.01}$} & 
        \multicolumn{1}{c|}{$\mathbf{\varepsilon=0.001}$}\\ 
        \cline{3-5}
        & & $\mathbf{||\text{div}_{\textit{h}} u_{\textit{h}}||_{L^1}}$ & 
         $\mathbf{||\text{div}_{\textit{h}} u_{\textit{h}}||_{L^1}}$ & 
         $\mathbf{||\text{div}_{\textit{h}} u_{\textit{h}}||_{L^1}}$ \\ 
        \hline
        10  & 0.1  & 4.465 $\times 10^{-6}$  & 4.582 $\times 10^{-6}$ &  2.786 $\times 10^{-6}$   \\  
        20  & 0.05  & 2.060 $\times 10^{-5}$  & 2.049 $\times 10^{-5}$ & 1.619 $\times 10^{-5}$ \\  
        25 & 0.04  & 3.261 $\times 10^{-5}$  & 3.087 $\times 10^{-5}$ & 1.908 $\times 10^{-5}$ \\  
        50 & 0.02 & 6.927 $\times 10^{-5}$  & 4.194 $\times 10^{-5}$  & 1.646 $\times 10^{-5}$\\ 
        100 & 0.01 & 1.235 $\times 10^{-4}$  & 3.285 $\times 10^{-5}$  & 8.070 $\times 10^{-6}$ \\ 
        \hline
    \end{tabular}
    
    \caption{\centering \textbf{Example \ref{Subsec:GV} - Gresho vortex problem:} $L^1$ norms of $\divh{} \mathbf{u}_h$ ($\lambda = 1$).}
    \label{tab: GV div_u}
\end{table}

\begin{table}[h!]
    \centering
    \renewcommand{\arraystretch}{1.3} 
    \setlength{\tabcolsep}{10pt}      
    
    \begin{tabular}{|c|c|c|c|c|c|c|c|}
        \hline
        \multirow{2}{*}{\textbf{N$_x$}} & \multirow{2}{*}{$\mathbf{\Delta x}$} & 
        \multicolumn{2}{c|}{$\mathbf{\varepsilon=0.1}$} & 
        \multicolumn{2}{c|}{$\mathbf{\varepsilon=0.01}$} & 
        \multicolumn{2}{c|}{$\mathbf{\varepsilon=0.001}$} \\ 
        \cline{3-8}
        & & $\mathbf{||\varrho \textbf{ error}||_{L^2}}$ & \textbf{EOC} 
        & $\mathbf{||\varrho \textbf{ error}||_{L^2}}$ & \textbf{EOC} 
        & $\mathbf{||\varrho \textbf{ error}||_{L^2}}$ & \textbf{EOC} \\ 
        \hline
        10  & 0.1  & 0.000509  & -      & 8.89 $\times 10^{-6}$  & -      & 9.91 $\times 10^{-8}$ & -      \\  
        20  & 0.05  & 0.000490  & 0.0550 & 8.36 $\times 10^{-6}$  & 0.0859 & 8.24 $\times 10^{-8}$ & 0.2464 \\  
        25 & 0.04  & 0.000476  & 0.1347 & 7.00 $\times 10^{-6}$  & 0.7982 & 5.67 $\times 10^{-8}$ & 1.6028 \\  
        50 & 0.02 & 0.000344  & 0.4666 & 2.55 $\times 10^{-6}$  & 1.4566 & 2.64 $\times 10^{-8}$ & 1.0700  \\    
        \hline
    \end{tabular}
    
    \caption{\centering \textbf{Example \ref{Subsec:GV} - Gresho vortex problem:} Convergence rates of $L^2$ error in $\varrho$  ($\lambda$ based on \eqref{lambda cond numerical}, with $c=100,200,200$ respectively for $\varepsilon=0.1,0.01,0.001$).}
    \label{tab: GV EOC_rho}
\end{table}

\begin{table}[h!]
    \centering
    \renewcommand{\arraystretch}{1.3} 
    \setlength{\tabcolsep}{10pt}      
    
    \begin{tabular}{|c|c|c|c|c|c|c|c|}
        \hline
        \multirow{2}{*}{\textbf{N$_x$}} & \multirow{2}{*}{$\mathbf{\Delta x}$} & 
        \multicolumn{2}{c|}{$\mathbf{\varepsilon=0.1}$} & 
        \multicolumn{2}{c|}{$\mathbf{\varepsilon=0.01}$} & 
        \multicolumn{2}{c|}{$\mathbf{\varepsilon=0.001}$} \\ 
        \cline{3-8}
        & & $\mathbf{||u \textbf{ error}||_{L^2}}$ & \textbf{EOC} 
        & $\mathbf{||u \textbf{ error}||_{L^2}}$ & \textbf{EOC} 
        & $\mathbf{||u \textbf{ error}||_{L^2}}$ & \textbf{EOC} \\ 
        \hline
        10  & 0.1  & 0.24507  & -      & 0.33001  & -      & 0.35898 & -      \\  
        20  & 0.05  & 0.20859  & 0.2325 & 0.25962  & 0.3461 & 0.23299 & 0.5785 \\  
        25 & 0.04  & 0.19235  & 0.3632 & 0.18818  & 1.4422 & 0.13495 & 2.3376 \\  
        50 & 0.02 & 0.11150  & 0.7867 & 0.05621  & 1.7433 & 0.05786 & 1.1865  \\  
        \hline
    \end{tabular}
    
    \caption{\centering \textbf{Example \ref{Subsec:GV} Gresho vortex problem:} Convergence rates of $L^2$ error in $\mathbf{u}$  ($\lambda$ based on \eqref{lambda cond numerical}, with $c=100,200,200$ respectively for $\varepsilon=0.1,0.01,0.001$).}
    \label{tab: GV EOC_u}
\end{table}

\subsection{Travelling vortex problem}
\label{Subsec:TV}
Analogously as in \cite{ap_lukacova3}, the computational domain is taken to be $\Omega=[0,1]\times[0,1]$, and the initial conditions are given by:
\begin{eqnarray}
\varrho^0 \left( x_1,x_2\right)= 110 + \varepsilon^2\left( \frac{1.5}{4\pi}\right)^2 D\left(x_1,x_2\right) \left( k\left(r\right)-k\left( \pi\right)\right), \\
u_1^0 \left( x_1,x_2\right)=  0.6+ 1.5 \left( 1+\cos \left( r \left(x_1,x_2\right)\right)\right) D\left(x_1,x_2\right) \left( 0.5-x_2\right),  \\
u_2^0 \left( x_1,x_2\right)=  0+ 1.5 \left( 1+\cos \left( r \left(x_1,x_2\right)\right)\right) D\left(x_1,x_2\right) \left( x_1-0.5\right), 
\end{eqnarray}
with 
\begin{eqnarray}
k\left(q\right)=2\cos\left(q\right)+2q \ \sin\left(q\right) + \frac{1}{8}\cos\left(2q\right) + \frac{1}{4}q \ \sin\left(2q\right) + \frac{3}{4}q^2, \\
r \left( x_1,x_2\right) = 4 \pi \left(  \left( x_1-0.5\right)^2 + \left( x_2-0.5\right)^2 \right)^{\frac{1}{2}}, \\
D\left(x_1,x_2\right) = \left\{ \begin{matrix}  1 & \text{if } r\left( x_1,x_2\right)<\pi \\ 0 & \text{otherwise}\end{matrix} \right..
\end{eqnarray}
We use periodic boundary conditions, and the parametric values are taken as: $\kappa=1$ and $\gamma=1.4$. The problem is considered with two different values of $\varepsilon$: $0.1,0.01$. A CFL number of $C=0.5$ is used for both values of $\varepsilon$. 
\par
\begin{figure}[h!]
\centering
\begin{subfigure}[b]{0.4\textwidth}
\centering
\includegraphics[width=\textwidth]{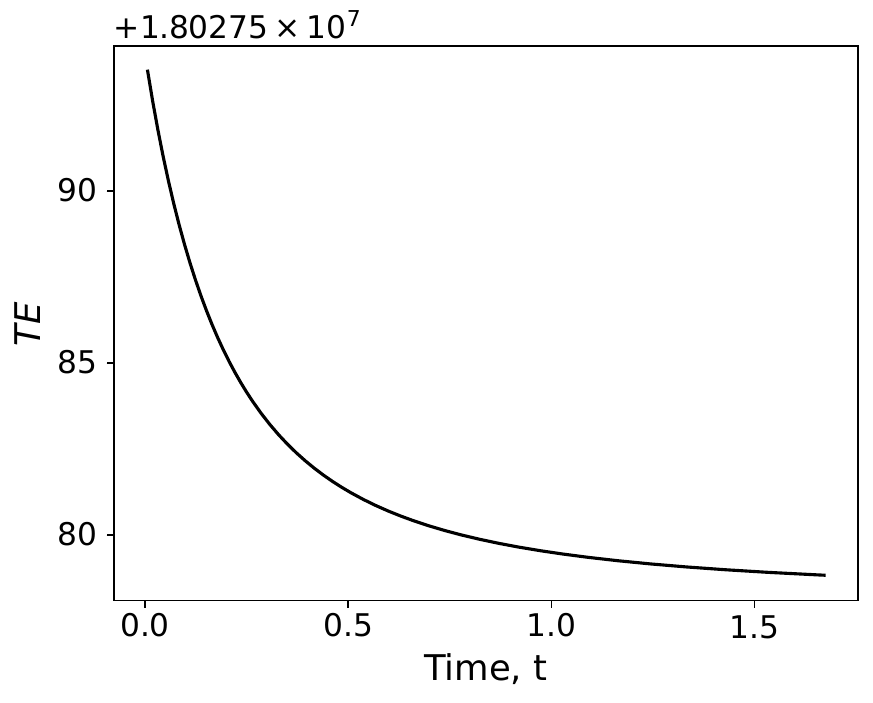}
\label{Fig:TV_TE_eps01}
\end{subfigure}
\hspace{-0.2cm}
\begin{subfigure}[b]{0.4\textwidth}
\centering
\includegraphics[width=\textwidth]{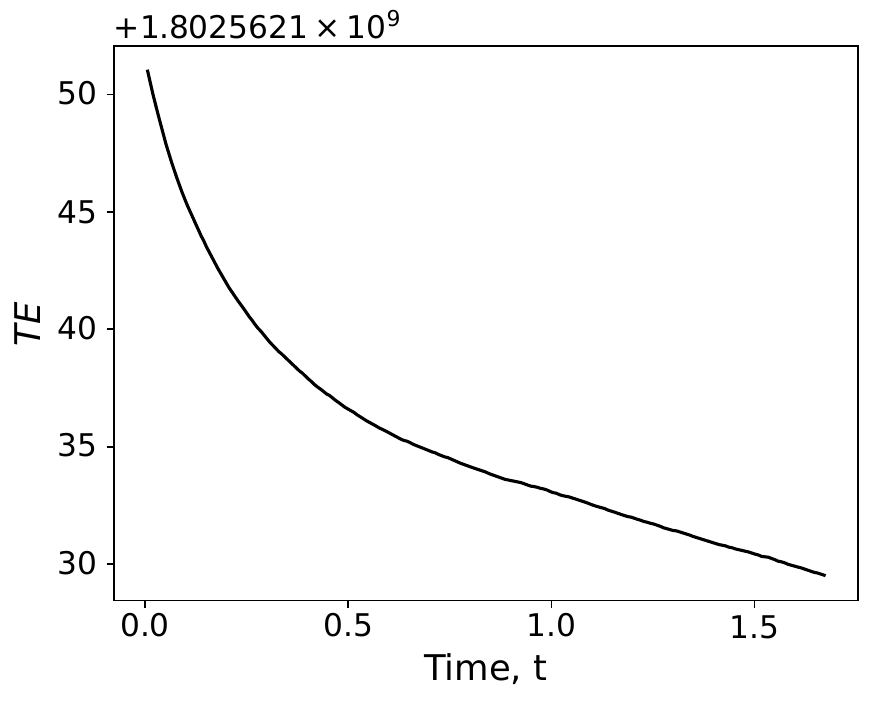}
\label{Fig:TV_TE_eps001}
\end{subfigure}

\begin{subfigure}[b]{0.4\textwidth}
\centering
\includegraphics[width=\textwidth]{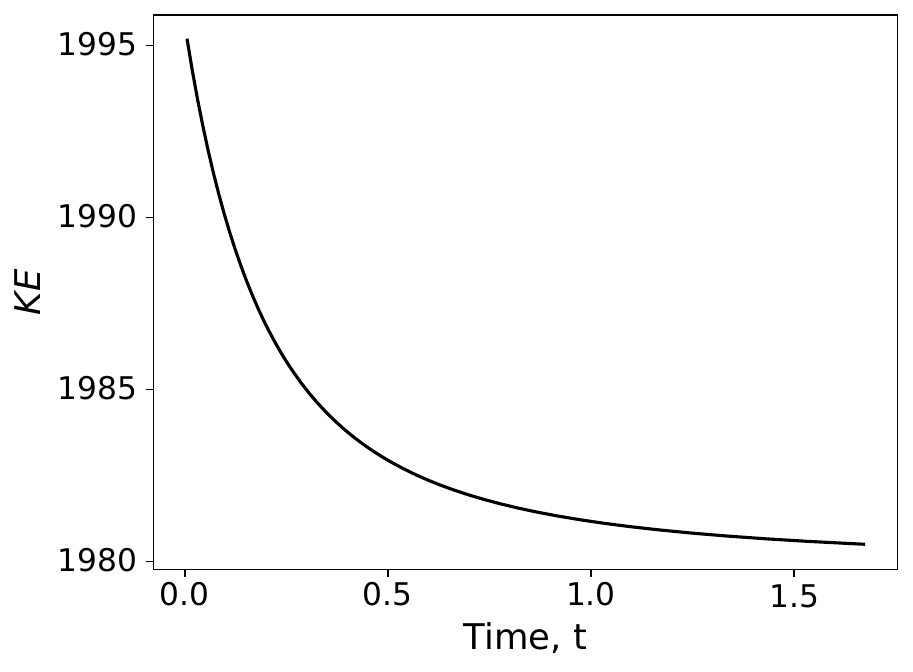}
\end{subfigure}
\hspace{-0.2cm}
\begin{subfigure}[b]{0.4\textwidth}
\centering
\includegraphics[width=\textwidth]{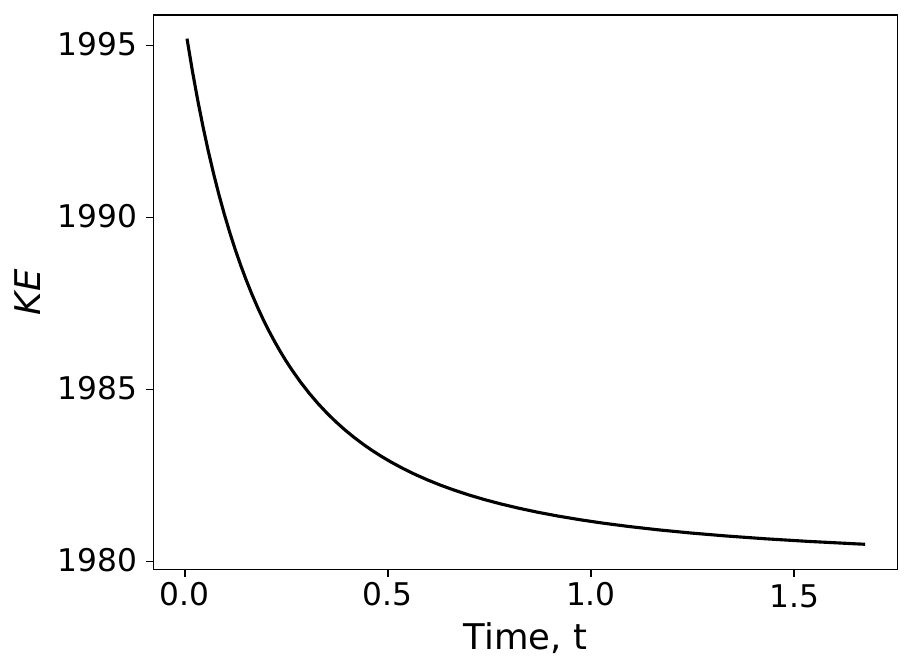}
\end{subfigure}

\begin{subfigure}[b]{0.45
\textwidth}
\centering
\includegraphics[width=\textwidth]{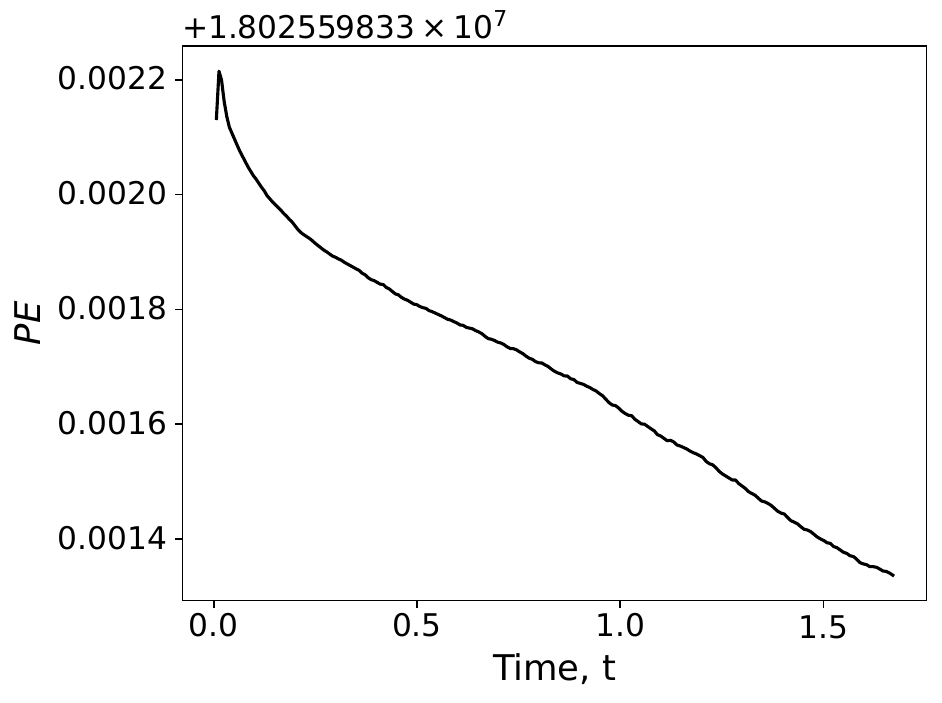}
\caption{$\varepsilon=0.1$}
\end{subfigure}
\hspace{-0.2cm}
\begin{subfigure}[b]{0.4\textwidth}
\centering
\includegraphics[width=\textwidth]{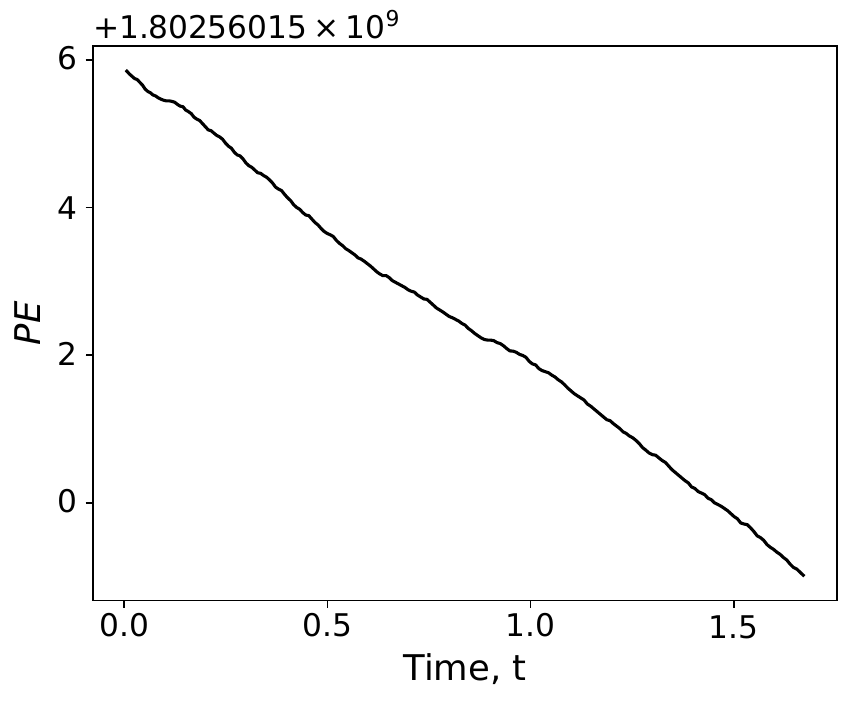}
\caption{$\varepsilon=0.01$}
\end{subfigure}

\caption{\centering \textbf{Example \ref{Subsec:TV} - Travelling vortex problem:} Top - Total energy; Middle - Kinetic energy; Bottom - Potential energy ($\lambda=1$)}
\label{Fig:TV_energy}
\end{figure}

\begin{figure}[h!]
\centering
\begin{subfigure}[b]{0.4\textwidth}
\centering
\includegraphics[width=\textwidth]{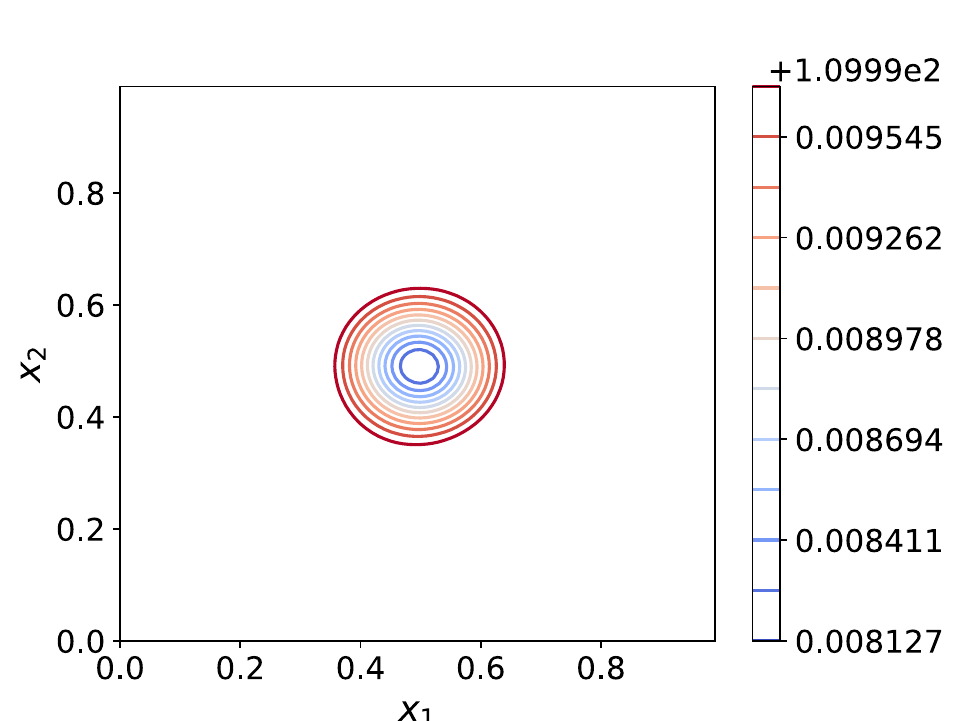}
\label{Fig:TV_contour_eps01}
\end{subfigure}
\hspace{-0.2cm}
\begin{subfigure}[b]{0.4\textwidth}
\centering
\includegraphics[width=\textwidth]{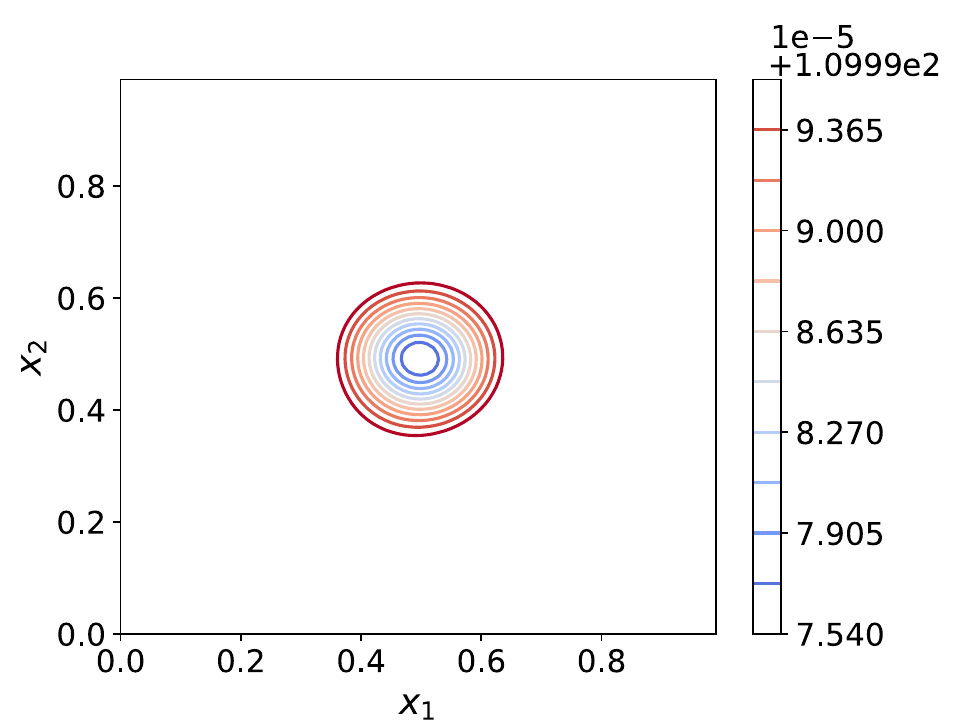}
\label{Fig:TV_contour_eps001}
\end{subfigure}
\vfill
\begin{subfigure}[b]{0.45\textwidth}
\centering
\includegraphics[width=\textwidth]{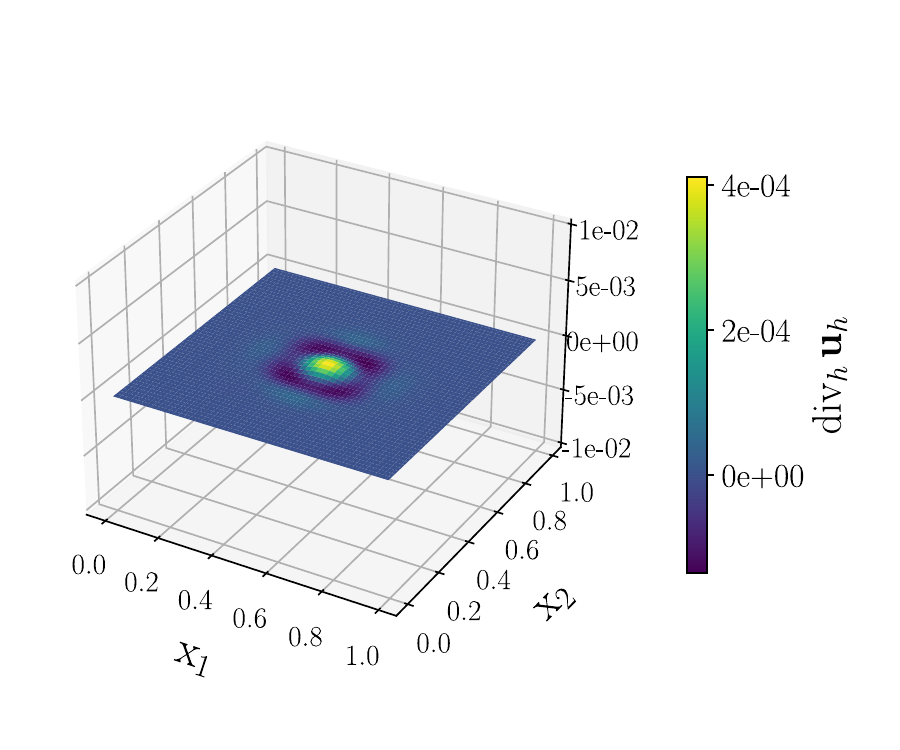}
\caption{$\varepsilon=0.1$}
\end{subfigure}
\hspace{-0.2cm}
\begin{subfigure}[b]{0.45\textwidth}
\centering
\includegraphics[width=\textwidth]{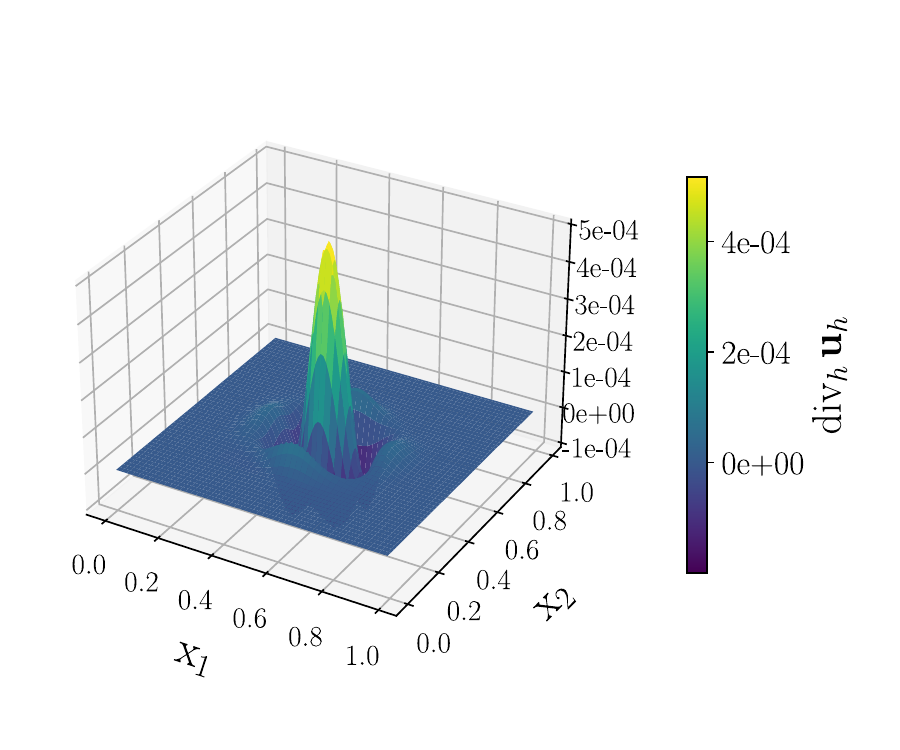}
\caption{$\varepsilon=0.01$}
\end{subfigure}

\caption{\centering \textbf{Example \ref{Subsec:TV} - Travelling vortex problem:} Top - Contour plots of $\varrho$; Bottom - Surface plots of $\divh{} \mathbf{u}_h$ ($\lambda$ based on \eqref{lambda cond numerical})} 
\label{Fig:TV_contour}
\end{figure}

 With $\lambda=1$, we observe energy stability for both the values of $\varepsilon$, as shown in \Cref{Fig:TV_energy} for N$_x=100$. In addition to total energy stability, the kinetic energy also decreases in time. In this example, it was also suitable to compute $\lambda$ at every time step based on \eqref{lambda cond numerical} with $c=30,200$, respectively, for $\varepsilon=0.1,0.01$. For a fixed grid with N$_x=50$, this condition results in $\lambda$ lying in the interval $[0.09245,0.2409]$, $[0.0239,0.7630]$, respectively, for $\varepsilon=0.1,0.01$, indicating that the different choices of $c$ for different $\varepsilon$ result in reasonable values of $\lambda$. A convergence study is performed at a time $t=1/0.6$ (one period), by considering different numbers of grid points as N$_x=$N$_{x_1}=$N$_{x_2}=10, 20, 25$, and $50$. The $L^2$ error is computed with the reference as N$_x=100$, and the convergence rates of $\varrho$ and $\mathbf{u}$ respectively are shown in \Cref{tab: TV EOC_rho,tab: TV EOC_u}. It can be seen that $\varrho$ converges with rates of about $0.5,0.87$, respectively, for $\varepsilon=0.1,0.01$, and $\mathbf{u}$ converges with rates of about $0.87,1.2$, respectively, for $\varepsilon=0.1,0.01$. Similar to the Gresho vortex problem, the $L^2$ errors in density for the travelling vortex problem (in \Cref{tab: TV EOC_rho}) are of $O\paraL{10^{-4}} \approx O\paraL{\varepsilon^{4}}$ for $\varepsilon=0.1$ and these are orders of magnitude smaller than the $O\paraL{\varepsilon^2}$ error due to the pressure linearization. \\
Further, \Cref{tab: TV div_u} shows that the $L^1$ norms of $\divh{} \mathbf{u}_h$ are $\mathcal{O}\paraL{10^{-3}}, \mathcal{O}\paraL{10^{-5}}$, respectively, for $\varepsilon=0.1,0.01$. \Cref{Fig:TV_contour} shows the contour plots of density, and the plots of $\divh{} \mathbf{u}_{h}$ for N$_x=100$ with $\varepsilon=0.1,0.01$. We observe from the $\divh{} \mathbf{u}_{h}$ plot that the method satisfies the divergence-free constraint upto $\mathcal{O}\paraL{\varepsilon^2}$. 

\begin{table}[h!]
    \centering
    \renewcommand{\arraystretch}{1.3} 
    \setlength{\tabcolsep}{10pt}      
    
    \begin{tabular}{|c|c|c|c|c|c|}
        \hline
        \multirow{2}{*}{\textbf{N$_x$}} & \multirow{2}{*}{$\mathbf{\Delta x}$} & 
        \multicolumn{2}{c|}{$\mathbf{\varepsilon=0.1}$} & 
        \multicolumn{2}{c|}{$\mathbf{\varepsilon=0.01}$} \\ 
        \cline{3-6}
        & & $\mathbf{||\varrho \textbf{ error}||_{L^2}}$ & \textbf{EOC} 
        & $\mathbf{||\varrho \textbf{ error}||_{L^2}}$ & \textbf{EOC}  \\ 
        \hline
        10  & 0.1  & 2.559 $\times 10^{-4}$ & -      & 3.265 $\times 10^{-6}$  & -          \\  
        20  & 0.05  & 2.369 $\times 10^{-4}$ & 0.1113 & 3.127 $\times 10^{-6}$  & 0.0625  \\  
        25 & 0.04  & 2.302 $\times 10^{-4}$ & 0.1275 & 2.960 $\times 10^{-6}$  & 0.2456  \\  
        50 & 0.02 & 1.634 $\times 10^{-4}$ & 0.4946 & 1.618 $\times 10^{-6}$  & 0.8716  \\    
        \hline
    \end{tabular}
    
    \caption{\centering \textbf{Example \ref{Subsec:TV} - Travelling vortex problem:} Convergence rates of $L^2$ error in $\varrho$  ($\lambda$ based on \eqref{lambda cond numerical}).}
    \label{tab: TV EOC_rho}
\end{table}

\begin{table}[h!]
    \centering
    \renewcommand{\arraystretch}{1.3} 
    \setlength{\tabcolsep}{10pt}      
    
    \begin{tabular}{|c|c|c|c|c|c|}
        \hline
        \multirow{2}{*}{\textbf{N$_x$}} & \multirow{2}{*}{$\mathbf{\Delta x}$} & 
        \multicolumn{2}{c|}{$\mathbf{\varepsilon=0.1}$} & 
        \multicolumn{2}{c|}{$\mathbf{\varepsilon=0.01}$} \\ 
        \cline{3-6}
        & & $\mathbf{||u \textbf{ error}||_{L^2}}$ & \textbf{EOC} 
        & $\mathbf{||u \textbf{ error}||_{L^2}}$ & \textbf{EOC}\\ 
        \hline
        10  & 0.1  & 4.428 $\times 10^{-2}$ & -      & 3.6566 $\times 10^{-2}$  & -           \\  
        20  & 0.05  & 2.644 $\times 10^{-2}$ & 0.7439 & 3.1727 $\times 10^{-2}$  & 0.2048  \\  
        25 & 0.04  & 2.475 $\times 10^{-2}$ & 0.2969 & 2.8092 $\times 10^{-2}$  & 0.5452  \\  
        50 & 0.02 & 1.351 $\times 10^{-2}$ & 0.8736 & 1.2186 $\times 10^{-2}$  & 1.2049  \\  
        \hline
    \end{tabular}
    
    \caption{\centering \textbf{Example \ref{Subsec:TV} - Travelling vortex problem:} Convergence rates of $L^2$ error in $\mathbf{u}$  ($\lambda$ based on \eqref{lambda cond numerical}).}
    \label{tab: TV EOC_u}
\end{table}

\begin{table}[h!]
    \centering
    \renewcommand{\arraystretch}{1.3} 
    \setlength{\tabcolsep}{10pt}      
    
    \begin{tabular}{|c|c|c|c|}
        \hline
        \multirow{2}{*}{\textbf{N$_x$}} & \multirow{2}{*}{$\mathbf{\Delta x}$} & 
        \multicolumn{1}{c|}{$\mathbf{\varepsilon=0.1}$} & 
        \multicolumn{1}{c|}{$\mathbf{\varepsilon=0.01}$} \\ 
        \cline{3-4}
        & & $\mathbf{||\text{div}_{\textit{h}} u_{\textit{h}}||_{L^1}}$ & 
         $\mathbf{||\text{div}_{\textit{h}} u_{\textit{h}}||_{L^1}}$ \\ 
        \hline
        10  & 0.1  & 2.078 $\times 10^{-3}$  & 6.954 $\times 10^{-5}$      \\  
        20  & 0.05  & 5.884 $\times 10^{-5}$  & 4.923 $\times 10^{-5}$  \\  
        25 & 0.04  & 3.953 $\times 10^{-5}$  & 5.494 $\times 10^{-5}$  \\  
        50 & 0.02 & 2.537 $\times 10^{-5}$  & 5.596 $\times 10^{-5}$  \\ 
        100 & 0.01 & 1.309 $\times 10^{-5}$  & 1.710 $\times 10^{-5}$  \\ 
        \hline
    \end{tabular}
    
    \caption{\centering \textbf{Example \ref{Subsec:TV} - Travelling vortex problem:} $L^1$ norms of $\divh{} \mathbf{u}_h$  ($\lambda$ based on \eqref{lambda cond numerical}).}
    \label{tab: TV div_u}
\end{table}

\section{Conclusion}
In this paper, we have proposed a new implicit-explicit finite volume method for the barotropic Euler system. Under a suitable numerical diffusion coefficient independent of Mach number $\varepsilon$, the method is provably positivity-preserving and energy-stable for any $\varepsilon>0$. We have also rigorously shown the asymptotic consistency of the method \eqref{Update mass}-\eqref{Update mom}, that is, it yields a consistent approximation of the incompressible system as $\varepsilon \to 0$. The numerical results validate that the proposed finite volume method is energy-stable for problems with both smooth and discontinuous solutions. The consistency and convergence analyses of the proposed asymptotic-preserving finite volume method are interesting topics for further research, and our rigorous proof of energy stability provides a necessary step in this direction. 

\section{Appendix}

\subsection{Discrete renormalization inequality}
In this section, we provide the proof of \Cref{Lem: RE}, which is crucial in showing that the density corresponding to the numerical scheme remains positive. 
\subsubsection{Proof of \Cref{Lem: RE}}
\label{Proof: RE}
    Let us consider the terms that we obtain by multiplying \eqref{Num mass} with $B'(\varrho_K^{n+1})$ one-by-one. \\
\emph{Time derivative term:}
\begin{equation}
   B'\left(\varrho_K^{n+1}\right) \frac{\varrho_{K}^{n+1}-\varrho_{K}^{n}}{\Delta t^n} = \frac{B\left(\varrho_{K}^{n+1}\right)-B\left(\varrho_{K}^{n}\right)}{\Delta t^n} + \frac{ B\left(\varrho_{K}^{n}\right) -  B\left(\varrho_{K}^{n+1}\right) - B'\left(\varrho_K^{n+1}\right) \left(\varrho_{K}^{n}-\varrho_{K}^{n+1}\right) }{\Delta t^n}.
\end{equation}
\emph{Convective term:}
\begin{eqnarray}
    B'\left(\varrho_K^{n+1}\right) \paraL{\divh{} \paraL{\varrho_h^{n+1}\mathbf{u}_h^{n+1}}}_K &=& B'\left(\varrho_K^{n+1}\right) \sumintKK{} \avg{ \varrho_h^{n+1} \mathbf{u}_h^{n+1} }_{\sigma} \cdot \mathbf{n}_{\sigma,K}  \\
   &=&  B'\left(\varrho_K^{n+1}\right)  \varrho_K^{n+1}  \sumintKK{} \avg{\mathbf{u}_h^{n+1} }_{\sigma} \cdot \mathbf{n}_{\sigma,K} \nonumber \\
    &&+ B'\left(\varrho_K^{n+1}\right) \sumintKK{} \frac{1}{2} \diff{\varrho_h^{n+1}}_{\sigma,K} \mathbf{u}_L^{n+1} \cdot \mathbf{n}_{\sigma,K} \nonumber \\
    &=& \left( B'\left(\varrho_K^{n+1}\right)  \varrho_K^{n+1} - B\left(\varrho_K^{n+1}\right) \right) \sumintKK{} \avg{\mathbf{u}_h^{n+1} }_{\sigma} \cdot \mathbf{n}_{\sigma,K} \nonumber \\ &&+ \sumintKK{} \avg{B\left(\varrho_h^{n+1}\right)\mathbf{u}_h^{n+1}}_{\sigma}\cdot \mathbf{n}_{\sigma,K} \nonumber \\
    && + B'\left(\varrho_K^{n+1}\right) \sumintKK{} \frac{1}{2} \diff{\varrho_h^{n+1}}_{\sigma,K} \mathbf{u}_L^{n+1} \cdot \mathbf{n}_{\sigma,K}  + D,
\end{eqnarray}
where 
\begin{eqnarray*}
    D &=& - \sumintKK{} \avg{B\left(\varrho_h^{n+1}\right)\mathbf{u}_h^{n+1}}_{\sigma}\cdot \mathbf{n}_{\sigma,K} + B\left(\varrho_K^{n+1}\right)  \sumintKK{} \avg{\mathbf{u}_h^{n+1} }_{\sigma} \cdot \mathbf{n}_{\sigma,K} \\
    &=& -\sumintKK{} \frac{1}{2} \diff{B\left(\varrho_h^{n+1}\right)}_{\sigma,K} \mathbf{u}_L^{n+1} \cdot \mathbf{n}_{\sigma,K}.
\end{eqnarray*}
Since the periodic boundary conditions are used, we have 
\begin{equation*}
    \sumK{} \sumintKK{}\avg{B\left(\varrho_h^{n+1}\right)\mathbf{u}_h^{n+1}}_{\sigma}\cdot \mathbf{n}_{\sigma,K} = \sumK{} \paraL{\divh{} \paraL{B\left(\varrho_h^{n+1}\right)\mathbf{u}_h^{n+1}}}_K = 0.
\end{equation*}
Hence, we obtain,
\begin{multline}
    \sumK{} B'\left(\varrho_K^{n+1}\right) \paraL{\divh{} \paraL{\varrho_h^{n+1}\mathbf{u}_h^{n+1}}}_K   = \sumK{}  \paraL{B'\left(\varrho_K^{n+1}\right)  \varrho_K^{n+1} - B\left(\varrho_K^{n+1}\right)} \paraL{\divh{} \paraL{\mathbf{u}_h^{n+1}}}_K   \\ 
     -  \sumK{} \sumintKK{} \frac{1}{2} \left( \diff{B\left(\varrho_h^{n+1}\right)}_{\sigma,K}  -  B'\left(\varrho_K^{n+1}\right) \diff{\varrho_h^{n+1}}_{\sigma,K}  \right) \mathbf{u}_L^{n+1} \cdot \mathbf{n}_{\sigma,K}.  
\end{multline}
\emph{Numerical diffusion term:}
Considering $\lambda_{\sigma,K} = \lambda_{\sigma,L} = \lambda_{\sigma}$ for each $\sigma = K | L$, we observe the following:
\begin{eqnarray}
    \sumK{} B'\left(\varrho_K^{n+1}\right) \sumintKK{} \lambda_{\sigma,K} \diff{\varrho_h^{n+1}}_{\sigma,K} = - \sumintall{} \lambda_{\sigma} \diff{B'\left(\varrho_h^{n+1}\right)}_{\sigma}\diff{\varrho_h^{n+1}}_{\sigma}. 
\end{eqnarray}
\emph{Discrete renormalization identity:}
Collecting all the terms from above, $\sumK{} B'(\varrho_K^{n+1}) $ \eqref{Num mass} yields the discrete renormalization identity \eqref{RE}.

\subsection{Discrete kinetic energy inequality}
In this section, we provide the proof of the discrete kinetic energy inequality (\Cref{Lem: KE_ineq}), which is crucial to proving the discrete energy inequality.  

\subsubsection{Proof of \Cref{Lem: KE_ineq}}
\label{Proof: KE}
Let us consider the terms that we obtain by performing $-\frac{\modL{\mathbf{u}_K^{n+1}}^2}{2}$ \eqref{Num mass} $+ \mathbf{u}_K^{n+1} \cdot $ \eqref{Num mom} one-by-one. \\ 
\emph{Time derivative terms:}
\begin{eqnarray}
    -\frac{\modL{\mathbf{u}_K^{n+1}}^2}{2} \frac{\varrho_{K}^{n+1}-\varrho_{K}^{n}}{\Delta t^n} + \mathbf{u}_K^{n+1} \cdot \frac{\varrho_{K}^{n+1} \mathbf{u}_{K}^{n+1}-\varrho_{K}^{n}\mathbf{u}_{K}^{n}}{\Delta t^n} &=& \frac{\varrho_{K}^{n+1} \modL{\mathbf{u}_K^{n+1}}^2 - \varrho_{K}^{n} \modL{\mathbf{u}_K^{n}}^2}{2 \Delta t^n} + \frac{\varrho_{K}^{n} \left( \mathbf{u}_K^{n+1} - \mathbf{u}_K^{n} \right)^2 }{2 \Delta t^n}.
\end{eqnarray}
\emph{Convective terms:}
\begin{eqnarray}
     - \frac{\modL{\mathbf{u}_K^{n+1}}^2}{2} \paraL{\divh{} \paraL{\varrho_h^{n+1}\mathbf{u}_h^{n+1}}}_K &=& - \frac{\modL{\mathbf{u}_K^{n+1}}^2}{2} \sumintKK{} \avg{\varrho_h^{n+1}\mathbf{u}_h^{n+1} }_{\sigma} \cdot \mathbf{n}_{\sigma,K} \nonumber \\ 
     &=& - \sumintKK{}{}  \avgL{\varrho_h^{n+1}\frac{\modL{ \mathbf{u}_h^{n+1} }^2}{2} \mathbf{u}_h^{n+1}}_{\sigma} \cdot \mathbf{n}_{\sigma,K} \nonumber 
    \\ && +\sumintKK{} \frac{1}{2} \varrho_L^{n+1} \mathbf{u}_L^{n+1} \cdot \mathbf{n}_{\sigma,K} \diffL{\frac{\modL{\mathbf{u}_h^{n+1}}^2}{2}}_{\sigma,K}, \\
    \mathbf{u}_K^{n+1} \cdot \paraL{\divh{} \paraL{\varrho_h^{n} \mathbf{u}_h^{n} \otimes \mathbf{u}_h^{n}}}_K &=& \mathbf{u}_K^{n+1} \cdot \sumintKK{} \avg{\varrho_h^{n}\mathbf{u}_h^{n} \otimes \mathbf{u}_h^{n} }_{\sigma} \cdot \mathbf{n}_{\sigma,K} \nonumber 
    \\ &=& 2 \sumintKK{} \avgL{\varrho_h^{n}\frac{ \mathbf{u}_h^{n} \cdot \mathbf{u}_h^{n+1}}{2} \mathbf{u}_h^{n}}_{\sigma} \cdot \mathbf{n}_{\sigma,K} \nonumber \\
     && - \sumintKK{} \varrho_L^{n} \mathbf{u}_L^{n} \cdot \mathbf{n}_{\sigma,K} \mathbf{u}_L^{n} \cdot \diffL{\frac{\mathbf{u}_h^{n+1}}{2}}_{\sigma,K}.
\end{eqnarray}
From the above, we interpret the discrete kinetic energy flux as
\begin{equation}
    \mathbf{\mathcal{F}}_{{KE}_{\sigma}} = \avgL{ - \varrho_h^{n+1}\frac{\modL{ \mathbf{u}_h^{n+1} }^2}{2} \mathbf{u}_h^{n+1} + 2 \varrho_h^{n}\frac{ \mathbf{u}_h^{n} \cdot \mathbf{u}_h^{n+1}}{2} \mathbf{u}_h^{n} }_{\sigma},
\end{equation}
and note that it is discrete divergence-free due to periodic boundaries, \textit{i.e.}, \begin{equation}
    \sumK{} \sumintKK{} \mathbf{\mathcal{F}}_{{KE}_{\sigma}} \cdot  \mathbf{n}_{\sigma,K} = 0. 
\end{equation} \\
The remaining terms become
\begin{eqnarray}
    \sumK{} \sumintKK{} \frac{1}{2} \varrho_L^{n+1} \mathbf{u}_L^{n+1} \cdot \mathbf{n}_{\sigma,K} \diffL{\frac{\modL{\mathbf{u}_h^{n+1}}^2}{2}}_{\sigma,K} &=& \sumintall{}  \avg{\varrho_h^{n+1} \mathbf{u}_h^{n+1}}_{\sigma} \cdot \mathbf{n}_{\sigma}  \diffL{\frac{\modL{\mathbf{u}_h^{n+1}}^2}{2}}_{\sigma}, \\
    - \sumK{} \sumintKK{} \varrho_L^{n} \mathbf{u}_L^{n} \cdot \mathbf{n}_{\sigma,K} \mathbf{u}_L^{n} \cdot \diffL{\frac{\mathbf{u}_h^{n+1}}{2}}_{\sigma,K} &=& - \sumintall{} \diff{\mathbf{u}_h^{n+1}}_{\sigma} \cdot \avg{\varrho_h^n \mathbf{u}_h^n \mathbf{u}_h^n}_{\sigma} \cdot \mathbf{n}_{\sigma}.
\end{eqnarray}
\emph{Numerical diffusion terms:}
Considering $\lambda_{\sigma,K} = \lambda_{\sigma,L} = \lambda_{\sigma}$ for each $\sigma = K | L$, we obtain
\begin{eqnarray} 
     \sumK{} \paraL{- \frac{\modL{\mathbf{u}_K^{n+1}}^2}{2}} \sumintKK{} \lambda_{\sigma,K} \diff{\varrho_h^{n+1}}_{\sigma,K} &=& \sumintall{} \lambda_{\sigma} \diff{\varrho_h^{n+1}}_{\sigma} \diffL{\frac{\modL{\mathbf{u}_h^{n+1}}^2}{2}}_{\sigma}, \\
     \sumK{} \mathbf{u}_K^{n+1} \cdot \sumintKK{} \lambda_{\sigma,K} \diff{\varrho_h^{n+1}\mathbf{u}_h^{n+1}}_{\sigma,K} &=& - \sumintall{} \lambda_{\sigma} \diff{\varrho_h^{n+1}\mathbf{u}_h^{n+1}}_{\sigma} \cdot \diff{\mathbf{u}_h^{n+1}}_{\sigma}.
\end{eqnarray}
Noting that $\diff{\varrho_h^{n+1}\mathbf{u}_h^{n+1}}_{\sigma} \cdot \diff{\mathbf{u}_h^{n+1}}_{\sigma} - \diff{\varrho_h^{n+1}}_{\sigma} \diffL{\frac{\modL{\mathbf{u}_h^{n+1}}^2}{2}}_{\sigma} = \avg{\varrho_h^{n+1}}_{\sigma} \modL{\diff{\mathbf{u}_h^{n+1}}_{\sigma}}^2$, we obtain the following non-positive quantity upon summing the above two equations, 
\begin{equation}
     - \sumintall{} \lambda_{\sigma} \avg{\varrho_h^{n+1}}_{\sigma} \modL{\diff{\mathbf{u}_h^{n+1}}_{\sigma}}^2.
\end{equation}
\emph{Pressure term:}
Since the boundaries are periodic, the following holds due to \Cref{Lem: Grad-div duality}:
\begin{equation}
    \sumK{} \mathbf{u}_K^{n+1} \cdot  \paraL{\nabla_h p\paraL{\varrho_h^{n+1}}}_K  = - \sumK{} p\paraL{\varrho_K^{n+1}} \paraL{\text{div}_h \mathbf{u}^{n+1}_h}_K.
\end{equation}
\emph{Discrete kinetic energy equation:}
Collecting all the terms from above, we obtain the following by performing $\sumK{} \paraL{-\frac{\modL{\mathbf{u}_K^{n+1}}^2}{2}}$ \eqref{Num mass} $+ \sumK{} \mathbf{u}_K^{n+1} \cdot $ \eqref{Num mom}:
\begin{multline}  
\sumK{} \frac{\varrho_{K}^{n+1} \modL{\mathbf{u}_K^{n+1}}^2 - \varrho_{K}^{n} \modL{\mathbf{u}_K^{n}}^2}{2 \Delta t^n} - \sumK{} \frac{p\paraL{\varrho_K^{n+1}}}{\varepsilon^2} \paraL{\text{div}_h \mathbf{u}^{n+1}_h}_K  \\ = - \sumK{}\frac{\varrho_{K}^{n} \left( \mathbf{u}_K^{n+1} - \mathbf{u}_K^{n} \right)^2 }{2 \Delta t^n}  - \sumintall{} \left(\lambda_{\sigma} \avg{\varrho_h^{n+1}}_{\sigma} \modL{\diff{\mathbf{u}_h^{n+1}}_{\sigma}}^2 \right. \\ \left. + \avg{\varrho_h^{n+1} \mathbf{u}_h^{n+1}}_{\sigma} \cdot \mathbf{n}_{\sigma}  \diffL{\frac{\modL{\mathbf{u}_h^{n+1}}^2}{2}}_{\sigma}  - \diff{\mathbf{u}_h^{n+1}}_{\sigma} \cdot \avg{\varrho_h^n \mathbf{u}_h^n \mathbf{u}_h^n}_{\sigma} \cdot \mathbf{n}_{\sigma}\right)
\end{multline}
Hence, if $\lambda_{\sigma}$ satisfies \eqref{lambda cond KE} for each $\sigma \in \mathcal{E}$ with $\diff{\mathbf{u}_h^{n+1}}_{\sigma} \neq 0$, then we have the discrete kinetic energy inequality \eqref{KE_ineq}.

\bibliographystyle{acm}
\bibliography{references}

\end{document}